# ON MANY-SERVER QUEUES IN HEAVY TRAFFIC

By Anatolii A. Puhalskii and Josh E. Reed

*University of Colorado Denver and IITP, Moscow, and New York University*

We establish a heavy-traffic limit theorem on convergence in distribution for the number of customers in a many-server queue when the number of servers tends to infinity. No critical loading condition is assumed. Generally, the limit process does not have trajectories in the Skorohod space. We give conditions for the convergence to hold in the topology of compact convergence. Some new results for an infinite server are also provided.

**1. Introduction.** Heavy-traffic limits for many-server queues is hardly a new topic. In particular, there exists a substantial body of literature on the Halfin–Whitt regime which is singled out by the requirements of critical loading and certain initial conditions. In many instances, by assuming that the service time distribution lies within a specific class of probability distributions, it has been shown that in the Halfin–Whitt regime the suitably centered and normalized processes of the number of customers in the system converge in distribution for Skorohod's $J_1$-topology to a process with continuous trajectories, which may be explicitly identified. This was first accomplished in the work of Halfin and Whitt [10] for the case of exponential service time distributions and has been continued by many additional authors as well for different classes of service time distributions, see, for instance, Jelenkovic, Mandelbaum, and Momcilovic [13], Mandelbaum and Momcilovic [18], Kaspi and Ramanan [14], Puhalskii and Reiman [21], Whitt [25] and the recent survey paper by Pang, Talreja and Whitt [19]. One exception to the above list is the work of Reed [22] in which no assumptions beyond a first moment are placed on the service time distribution. A related avenue of research concerns infinite servers in heavy traffic; see, for example,

Krichagina and Puhalskii [15], Pang, Talreja and Whitt [19] and references therein.

The purpose of this paper is to extend the aforementioned results to allow noncritical loading, generally distributed service times, and general initial conditions. We consider a sequence of $G/GI/n$ queues. The number of servers, $n$, tends to infinity and the service time distribution is held fixed. It is assumed that the centered and normalized arrival processes converge in distribution. The main result asserts the weak convergence of finite-dimensional distributions of the suitably centered and normalized number-in-the-system processes to the finite-dimensional distributions of a unique strong solution of a stochastic differential equation of convolution type. In the general case, the trajectories of the limit process have discontinuities of the second kind.

Under a certain condition on the fluid limit and the service time distribution, the limit stochastic process has trajectories in a Skorohod space and the convergence in distribution to this process holds for the $J_1$-topology. The results of Halfin and Whitt [10] and Reed [22] follow as particular cases. In fact, the convergence in distribution is stated for the stronger topology of uniform convergence on compact intervals, which we call the topology of compact convergence. This topology plays a prominent role throughout the paper so much so that the supporting results are established for this topology rather than for the $J_1$-topology. This emphasis is necessitated by the need to rely on certain continuity properties of equations describing the evolution of the system's population, which we are able to establish only for the topology of compact convergence.

Furthermore, the topology of compact convergence is almost too strong. Since it is nonseparable, the associated Borel $\sigma$-algebra is strictly larger than the Kolmogorov (i.e., cylindrical) $\sigma$-algebra, which leads to measurability problems; see, for example, Billingsley [1], Section 15. As a consequence, one needs to extend the notion of convergence in distribution. The relevant theory has been developed by Hoffman–Jørgensen and his followers and is expounded upon in van der Vaart and Wellner [23]. Its primary motivation was to provide tools to tackle convergence in distribution of empirical processes in strong topologies. Empirical processes also play an important part in the study of heavy traffic limits for many-server queues, see Louchard [17], Krichagina and Puhalskii [15]. This explains why we find methods developed in a different context useful for our setting.

One of the challenges of working with convergence in distribution in nonseparable spaces has to do with Borel probability measures not being necessarily tight. Recall that on complete spaces the tightness and separability properties are equivalent. Many conclusions of the Hoffman–Jørgensen theory rely on the assumption that the limit probability measure is separable. It is therefore important to show that the stochastic processes that appear in the limit have separable ranges, so their distributions are separable probability measures. We accomplish this by establishing that certain limit processes are Gaussian semimartingales. Since Gaussian semimartingales jump at deterministic times, their ranges are separable sets. As a side remark, there are no known examples of nonseparable Borel probability measures and the axiom that such measures do not exist can be consistently added to the Zermelo–Fränkel set of axioms; see Dudley [8], Appendix C, and van der Vaart and Wellner [23], Section 1.3, for these and other observations. Our results, however, make no use of this fact.

Another ingredient in the proof of the main theorem is a martingale argument which originates from Krichagina and Puhalskii [15]. It is instrumental in establishing tightness for the processes of interest. The key that makes an application of the methods of [15] possible is the insight of Reed [22] that the process of customers entering service in a many-server queue can be treated in analogy with the arrival process at an infinite server. There are, however, certain improvements on the derivation of the needed martingale properties as compared with the cited papers. In particular, our approach clarifies the connection with certain two-parameter processes being planar martingales, which was implicit in Krichagina and Puhalskii [15]. In addition, the construction of the limit process and a number of proofs in [15] relied on the continuity of the fluid limit. Here it is not necessarily the case, so a more subtle argument is called for. As a byproduct, in application to a $G/GI/\infty$ system our methods enable us to remove the restriction of Krichagina and Puhalskii [15] that the fluid limit for the arrival process should be continuous.

The paper is organized as follows. In Section 2, the main results are formulated and discussed and proofs are presented. We first recall basic facts about convergence in distribution for nonmeasurable mappings. We then state in Theorem 2.1 the fluid limit result of Reed [22], which is given a different proof. This proof is useful when dealing with a more general result of Theorem 2.3. Theorem 2.2 presents the stochastic approximation result. The proof relies heavily on the continuity property of convolution equations, on the one hand, and on a result of convergence in distribution of certain stochastic processes, on the other hand. The former result is proved in Appendix B and the latter result is proved in Section 4. The Gaussian limit for a $G/GI/\infty$ system is stated in Section 2.3. This subsection also contains extended versions of Theorems 2.1 and 2.2. These results are proved by similar means, so we omit the actual proofs. Appendix C contains the martingale arguments. In Appendix D, we summarize the properties of convergence in distribution for nonmeasurable mappings which are, for the most part, borrowed from van der Vaart and Wellner [23].

*Notation and conventions.* The set of real numbers is denoted by $\mathbb{R}$, the set of nonnegative reals is denoted by $\mathbb{R}_+$, the set of nonnegative rational numbers is denoted by $\mathbb{Q}_+$, the set of nonnegative integers is denoted by $\mathbb{Z}_+$, and the set of natural numbers is denoted by $\mathbb{N}$. For real numbers $x$ and $y$, $x \wedge y = \min(x, y)$, $x \vee y = \max(x, y)$, $x^+ = x \vee 0$, and $\lfloor x \rfloor$ denotes the integer part. Integrals of the form $\int_a^b$ are understood as $\int_{(a,b]}$ except when $a = 0$: we interpret $\int_0^a$ as $\int_{[0,a]}$. Similarly, $\int_{\mathbb{R}_+^2} = \int_{[0,\infty)\times[0,\infty)}$. For a function $g(x)$ of a nonnegative real-valued argument, $g(x-)$ denotes the left-hand limit at $x$, $g(x+)$ denotes the right-hand limit at $x$, and $\Delta g(x) = g(x) - g(x-)$. We always assume that $g(0-) = 0$, so $\Delta g(0) = g(0)$. If $f(x, y)$ is a function with domain $\mathbb{R}_+^2$ and $x_1 \leq x_2$ and $y_1 \leq y_2$, we define $\square f((x_1, y_1), (x_2, y_2)) = f(x_2, y_2) - f(x_2, y_1) - f(x_1, y_2) + f(x_1, y_1)$. $\mathbf{1}_A$ denotes the indicator function of set $A$. Products of topological spaces are assumed to be equipped with product topologies. All random entities are assumed to be defined on a common complete probability space $(\Omega, \mathcal{F}, \mathbf{P})$ with $\mathbf{E}$ denoting the associated expectation. Filtrations are defined as right-continuous flows of complete $\sigma$-algebras.

## 2. Main results.

2.1. *Extended convergence in distribution.* We recall the concept of convergence in distribution due to Hoffmann–Jørgensen as presented in van der Vaart and Wellner [23]. Let $\xi$ be a real-valued function on $\Omega$, which does not have to be a random variable, i.e., to be appropriately measurable. The outer expectation $\mathbf{E}^*\xi$ of $\xi$ is defined as the infimum of $\mathbf{E}\zeta$ over all random variables $\zeta$ on $(\Omega, \mathcal{F}, \mathbf{P})$ such that $\zeta \geq \xi$ a.s. and $\mathbf{E}\zeta$ is well defined. Let $\mathbb{S}$ be a metric space made into a measurable space by endowing it with the Borel $\sigma$-algebra $\mathcal{B}(\mathbb{S})$. Given a sequence $X_n$ of maps from $\Omega$ to $\mathbb{S}$ and a measurable map $X$ from $(\Omega, \mathcal{F})$ to $(\mathbb{S}, \mathcal{B}(\mathbb{S}))$, we say that the $X_n$ converge to $X$ in distribution in $\mathbb{S}$ if

$$\lim_{n\to\infty} \mathbf{E}^* f(X_n) = \mathbf{E} f(X)$$

for all bounded continuous real-valued functions $f$ on $\mathbb{S}$.

In this paper, $\mathbb{S}$ will have as its elements right-continuous with left-hand limits functions from an interval $I$ on the real line to a complete metric space $U$. We denote this set by $\mathbb{D}(I, U)$. If $U$ is, in addition, separable, then it is customary to equip $\mathbb{D}(I, U)$ with Skorohod's $J_1$-topology and a complete separable metric; see, for example, Jacod and Shiryaev [12], Ethier and Kurtz [9]. The resulting Polish space will also be denoted by $\mathbb{D}(I, U)$. We hope that no confusion will arise out of this ambiguity.

For the most part, though, we will be concerned with the stronger topology of compact convergence. This topology is compatible with a complete metric which can be defined by

$$d(\mathbf{x},\mathbf{y}) = \sum_{k=1}^{\infty} \sup_{t\in[0,k]\cap I} (1 \wedge \rho(\mathbf{x}(t),\mathbf{y}(t)))2^{-k},$$

where $\mathbf{x} = (\mathbf{x}(t), t \in I)$ and $\mathbf{y} = (\mathbf{y}(t), t \in I)$ are elements of $\mathbb{D}(I,U)$ and $\rho$ is the metric on $U$. We denote this metric space by $\mathbb{D}_c(I,U)$ and equip it with the Borel $\sigma$-algebra. One of the nice properties of the topology of compact convergence is that $\mathbb{D}_c(I,U_1 \times U_2)$ is homeomorphic to $\mathbb{D}_c(I,U_1) \times \mathbb{D}_c(I,U_2)$ for arbitrary metric spaces $U_1$ and $U_2$. [Recall that, by contrast, the topology of $\mathbb{D}(\mathbb{R}_+,\mathbb{R}^2)$ is strictly finer than the topology of $\mathbb{D}(\mathbb{R}_+,\mathbb{R})^2$].

The space $\mathbb{D}_c(I,U)$, though complete, is not separable unless either $I$ or $U$ represent a one-point set. Therefore, the Kolmogorov $\sigma$-algebra on $\mathbb{D}_c(I,U)$ is strictly smaller than the Borel $\sigma$-algebra, so a stochastic process $X$ with trajectories in $\mathbb{D}(I,U)$ need not be a random element of $\mathbb{D}_c(I,U)$. However, if the range of the mapping from $\Omega$ to $\mathbb{D}_c(I,U)$ defined by $X$ is a separable set, then $X$ is a random element of $\mathbb{D}_c(I,U)$. We call such $X$ a separable random element. In particular, if $X$ has continuous trajectories a.s., or more generally, if its jumps occur at deterministic times a.s., then $X$ is a separable random element of $\mathbb{D}_c(I,U)$. We say that $X$ is a tight random element if its distribution is a tight probability measure on $\mathbb{D}_c(I,U)$. As mentioned, $X$ is a tight random element if and only if it is a separable random element.

Suppose $X$ is a stochastic process which is a random element of $\mathbb{D}_c(I,U)$. We say that a sequence of stochastic processes $X_n$ with trajectories in $\mathbb{D}(I,U)$ converges in distribution in $\mathbb{D}_c(I,U)$ (or, for the topology of compact convergence) to $X$ if the associated maps $X_n$ converge to $X$ in distribution in $\mathbb{D}_c(I,U)$. Convergence in distribution in $\mathbb{D}(I,U)$ is defined in a standard fashion, provided $U$ is Polish. If the latter is the case, convergence in distribution in $\mathbb{D}_c(I,U)$ implies convergence in distribution in $\mathbb{D}(I,U)$ and convergence in distribution in $\mathbb{D}(I,U)$ to continuous-path processes implies convergence in distribution in $\mathbb{D}_c(I,U)$.

2.2. *Limit theorems for the many server queue.* We assume as given a sequence of many server queues indexed by $n$, where $n$ denotes the number of servers. Service is performed on a first-come–first-serve basis. We denote the number of customers present at time $0-$ by $Q_0^n$. Out of those, $Q_0^n \wedge n$ customers are in service at time $0-$ and $(Q_0^n - n)^+$ customers are in the queue. We denote the remaining service times of the customers in service at time $0-$ by $\tilde{\eta}_1, \tilde{\eta}_2, \ldots.$ The service times of the customers in the queue at time $0-$ and the service times of customers exogenously arriving after time $0-$ are denoted by $\eta_1, \eta_2, \ldots$ and come from an i.i.d. sequence of nonnegative random variables with distribution $F = (F(x), x \in \mathbb{R}_+)$. The $\eta_i$ are not equal to zero a.s. Equivalently, we require that

$$F(0) < 1. \tag{2.1}$$

We let $E^n(t)$ denote the number of exogenous arrivals by $t$. The entities $Q_0^n$, $\{\tilde{\eta}_1, \tilde{\eta}_2, \dots\}$, $\{\eta_1, \eta_2, \dots\}$, and $E^n = (E^n(t), t \in \mathbb{R}_+)$ are assumed to be independent.

Let $Q^n(t)$ denote the number of customers either in the queue or in service at time $t$, let $A^n(t)$ denote the number of customers that enter service after time $0-$ and by time $t$, and let $\tilde{Q}^n(t)$ denote the number of customers remaining in service at time $t$ out of those that were in service at time $0-$. The introduced stochastic processes are assumed to have trajectories from $\mathbb{D}(\mathbb{R}_+, \mathbb{R})$. The following equations appear in Reed [22], cf. also Krichagina and Puhalskii [15]: for $t \in \mathbb{R}_+$

$$Q^n(t) = (Q_0^n - n)^+ + \tilde{Q}^n(t) + E^n(t) - \int_0^t \int_0^t \mathbf{1}_{\{s+x \le t\}} \, d \sum_{i=1}^{A^n(s)} \mathbf{1}_{\{\eta_i \le x\}}, \tag{2.2a}$$

$$A^n(t) = (Q_0^n - n)^+ + E^n(t) - (Q^n(t) - n)^+, \tag{2.2b}$$

$$\tilde{Q}^n(t) = \sum_{i=1}^{Q_0^n \wedge n} \mathbf{1}_{\{\tilde{\eta}_i > t\}}. \tag{2.2c}$$

We note that the integral on the right-hand side of (2.2a) represents the number of departures by time $t$ of the customers that enter service after time $0-$ as seen from the representation

$$\int_0^t \int_0^t \mathbf{1}_{\{s+x \le t\}} \, d \sum_{i=1}^{A^n(s)} \mathbf{1}_{\{\eta_i \le x\}} = \sum_{i=1}^{A^n(t)} \mathbf{1}_{\{\tau_i^n + \eta_i \le t\}},$$

where

$$\tau_i^n = \inf\{t \in \mathbb{R}_+ : A^n(t) \ge i\}.$$

It is shown in Lemma A.1 below that, given $E^n$, $Q_0^n$, and the sequences $\{\eta_i\}$ and $\{\tilde{\eta}_i\}$, equations (2.2a)–(2.2c) admit solutions $A^n$, $Q^n$ and $\tilde{Q}^n$. The processes $A^n$ and $Q^n$ are specified uniquely a.s. under the additional requirement that $\eta_i > 0$ a.s. We also provide an example of nonuniqueness when the $\eta_i$ can equal zero. Since our results do not assume that the service times be positive a.s., we interpret processes $A^n$ and $Q^n$ in what follows as some solutions to (2.2a) and (2.2b) [the process $\tilde{Q}^n$ is specified uniquely by (2.2c)].

We now introduce the fluid limit equations. Let $\tilde{F} = (\tilde{F}(t), t \in \mathbb{R}_+)$ represent the distribution function of a nonnegative random variable. According to Lemma B.1, given $\mathbf{q}_0 \in \mathbb{R}_+$ and a nondecreasing function $\mathbf{e} = (\mathbf{e}(t), t \in \mathbb{R}_+) \in \mathbb{D}(\mathbb{R}_+, \mathbb{R})$, there exists a unique function $\mathbf{q} = (\mathbf{q}(t), t \in \mathbb{R}_+) \in \mathbb{D}(\mathbb{R}_+, \mathbb{R})$ which satisfies the equation

$$\begin{aligned} \mathbf{q}(t) &= \mathbf{q}_0 - (\mathbf{q}_0 - 1)^+ F(t) - \mathbf{q}_0 \wedge 1 \tilde{F}(t) + \mathbf{e}(t) \\ &\quad - \int_0^t \mathbf{e}(t-s)\, dF(s) + \int_0^t (\mathbf{q}(t-s) - 1)^+ \, dF(s). \end{aligned} \tag{2.3}$$

One can easily see that $\mathbf{q}$ is $\mathbb{R}_+$-valued. We also define a function $\mathbf{a} = (\mathbf{a}(t), t \in \mathbb{R}_+) \in \mathbb{D}(\mathbb{R}_+, \mathbb{R})$ by the relation

$$\mathbf{a}(t) = (\mathbf{q}_0 - 1)^+ + \mathbf{e}(t) - (\mathbf{q}(t) - 1)^+. \tag{2.4}$$

The next theorem is, in essence, due to Reed [22]. Let $\tilde{F}^n = (\tilde{F}^n(t), t \in \mathbb{R}_+)$ denote the empirical distribution function of $\tilde{\eta}_1, \dots, \tilde{\eta}_n$, that is,

$$\tilde{F}^n(t) = \frac{1}{n} \sum_{i=1}^n \mathbf{1}_{\{\tilde{\eta}_i \le t\}},$$

THEOREM 2.1. *Suppose that, for arbitrary $T > 0$ and $\varepsilon > 0$,*

$$\lim_{n\to\infty} \mathbf{P}\left( \left| \frac{Q_0^n}{n} - \mathbf{q}_0 \right| > \varepsilon \right) = 0,$$

$$\lim_{n\to\infty} \mathbf{P}\left( \sup_{t\in[0,T]} |\tilde{F}^n(t) - \tilde{F}(t)| > \varepsilon \right) = 0$$

*and*

$$\lim_{n\to\infty} \mathbf{P}\left( \sup_{t\in[0,T]} \left| \frac{E^n(t)}{n} - \mathbf{e}(t) \right| > \varepsilon \right) = 0.$$

*Then, for arbitrary $T > 0$ and $\varepsilon > 0$,*

$$\lim_{n\to\infty} \mathbf{P}\left( \sup_{t\in[0,T]} \left| \frac{Q^n(t)}{n} - \mathbf{q}(t) \right| > \varepsilon \right) = 0$$

*and*

$$\lim_{n\to\infty} \mathbf{P}\left( \sup_{t\in[0,T]} \left| \frac{A^n(t)}{n} - \mathbf{a}(t) \right| > \varepsilon \right) = 0.$$

Since $A^n(t)$ represents the number of customers that enter service by time $t$ out of those that are either queued at time $0-$ or exogenously arrive after time $0-$, the process $A^n$ is $\mathbb{R}_+$-valued and nondecreasing. Hence, under the hypotheses of Theorem 2.1, $\mathbf{a} = (\mathbf{a}(t), t \in \mathbb{R}_+)$ is a nondecreasing $\mathbb{R}_+$-valued function from $\mathbb{D}(\mathbb{R}_+, \mathbb{R})$. According to the Glivenko–Cantelli theorem, the compact convergence in probability of the $\tilde{F}^n$ to $\tilde{F}$ in the hypotheses holds if the $\tilde{\eta}_i$ are i.i.d. with distribution $\tilde{F}$.

An interpretation of equation (2.3) can be as follows. By (2.4), the function $\mathbf{a} = (\mathbf{a}(t), t \in \mathbb{R}_+)$ represents "fluid customers entering service after time 0−." If we write (2.3) in the form

$$\mathbf{q}(t) = \mathbf{q}_0 + \mathbf{e}(t) - \left(\mathbf{q}_0 \wedge 1\tilde{F}(t) + \int_0^t \mathbf{a}(t-s)\, dF(s)\right), \tag{2.5}$$

then the first two terms on the right have the meaning of "the number of fluid customers" seen by time $t$ and the sum in parentheses represents "fluid departures" by time $t$.

Condition (2.1) cannot be disposed of in general. If we suppose, for instance, that $F(0) = 1$, and $\tilde{F}(t_0) = 0$ for some $t_0 > 0$, then we may obtain different fluid limits for the cases where $Q_0^n = n$ and $Q_0^n = n - 1$. In the former case, $Q^n(t) = E^n(t) + n$ for all $t < t_0$ as all servers remain occupied until $t_0$, while in the latter case $Q^n(t) = n - 1$ for $t < t_0$ as the only available server will let all exogenously arriving customers out of the system immediately.

We now state the stochastic approximation result. Define processes $X^n = (X^n(t), t \in \mathbb{R}_+)$, $S^n = (S^n(t), t \in \mathbb{R}_+)$, and $Y^n = (Y^n(t), t \in \mathbb{R}_+)$ by

$$X^n(t) = \sqrt{n}\left(\frac{Q^n(t)}{n} - \mathbf{q}(t)\right), \tag{2.6}$$

$$S^n(t) = \sqrt{n}(\tilde{F}^n(t) - \tilde{F}(t)), \tag{2.7}$$

$$Y^n(t) = \sqrt{n}\left(\frac{E^n(t)}{n} - \mathbf{e}(t)\right). \tag{2.8}$$

We note that $X^n$, $S^n$ and $Y^n$ are random elements of $\mathbb{D}(\mathbb{R}_+, \mathbb{R})$. However, $X^n$ and $Y^n$ might not be random elements of $\mathbb{D}_c(\mathbb{R}_+, \mathbb{R})$. Let us also denote

$$X_0^n = \sqrt{n}\left(\frac{Q_0^n}{n} - \mathbf{q}_0\right).$$

We introduce the following condition:

$$\lim_{\varepsilon \to 0} \sup_{t \in [0,T]} \int_0^t \mathbf{1}_{\{0 < |\mathbf{q}(t-s) - 1| < \varepsilon\}}\, dF(s) = 0 \qquad \text{for all } T > 0. \tag{2.9}$$

THEOREM 2.2. *Suppose that, as $n \to \infty$, the $X_0^n$ converge in distribution in $\mathbb{R}$ to a random variable $X_0$, the $S^n$ converge in distribution in $\mathbb{D}_c(\mathbb{R}_+, \mathbb{R})$ to a process $S = (S(t), t \in \mathbb{R}_+)$, and the $Y^n$ converge in distribution in $\mathbb{D}_c(\mathbb{R}_+, \mathbb{R})$ to a process $Y = (Y(t), t \in \mathbb{R}_+)$, where $S$ and $Y$ are separable random elements of $\mathbb{D}_c(\mathbb{R}_+, \mathbb{R})$. Then finite-dimensional distributions of the processes $X^n$ weakly converge to finite-dimensional distributions* *of the process* $X = (X(t), t \in \mathbb{R}_+)$ *that is a unique strong solution to the equation*

$$X(t) = (X_0 \mathbf{1}_{\{\mathbf{q}_0 > 1\}} + X_0^+ \mathbf{1}_{\{\mathbf{q}_0 = 1\}}) \times (1 - F(t)) + (X_0 \mathbf{1}_{\{\mathbf{q}_0 < 1\}} + X_0 \wedge 0 \mathbf{1}_{\{\mathbf{q}_0 = 1\}})(1 - \tilde{F}(t)) + \sqrt{\mathbf{q}_0 \wedge 1} S(t) + Y(t) - \int_0^t Y(t-s)\, dF(s) + Z(t) + \int_0^t (X(t-s) \mathbf{1}_{\{\mathbf{q}(t-s) > 1\}} + X(t-s)^+ \mathbf{1}_{\{\mathbf{q}(t-s) = 1\}})\, dF(s),$$

*where* $Z = (Z(t), t \in \mathbb{R}_+)$ *is a zero-mean Gaussian semimartingale with trajectories from* $\mathbb{D}(\mathbb{R}_+, \mathbb{R})$ *specified by the covariance*

$$\mathbf{E} Z(s) Z(t) = \int_0^{s \wedge t} F(s \wedge t - u)(1 - F(s \vee t - u))\, d\mathbf{a}(u),$$

*the entities* $X_0$, $S$, $Z$ *and* $Y$ *being independent. The trajectories of* $X$ *are Borel measurable and locally bounded a.s. If, in addition, condition* (2.9) *holds, then* $X$ *is a separable random element of* $\mathbb{D}_c(\mathbb{R}_+, \mathbb{R})$ *and the* $X^n$ *converge in distribution in* $\mathbb{D}_c(\mathbb{R}_+, \mathbb{R})$ *to* $X$.

We now comment on this theorem. As with Theorem 2.1, the convergence of the $S^n$ in the hypotheses holds when the $\tilde{\eta}_i$ are i.i.d. with distribution $\tilde{F}$. In such a case, $S(t) = W^0(\tilde{F}(t))$, where $W^0 = (W^0(t), t \in [0,1])$ is a Brownian bridge.

Since the process $Z$ is a Gaussian semimartingale, its jumps occur at deterministic times which are the jump times of its variance, see Jacod and Shiryaev [12], Chapter II, Section 4d, and Liptser and Shiryayev [16], Chapter 4, Section 9, so this process is a separable random element of $\mathbb{D}_c(\mathbb{R}_+, \mathbb{R})$. If either $F$ or $\mathbf{a}$ is a continuous function, then $Z$ is a continuous path process. In particular, if $Y$, $F$ and $\tilde{F}$ are continuous (with $F(0) = 0$), then $X$ is continuous (see Lemma B.1).

Similarly to Krichagina and Puhalskii [15], $Z$ admits a representation as a stochastic integral, which follows from the proof of Lemma 4.4 below. Specifically, we may assume that a.s.

$$Z(t) = \int_{\mathbb{R}_+^2} \mathbf{1}_{\{s + x \le t\}}\, d\Gamma(s, x), \tag{2.10}$$

where $(\Gamma(s, x), s \in \mathbb{R}_+, x \in \mathbb{R}_+)$ is a zero-mean Gaussian process with covariance

$$\mathbf{E}\Gamma(s, x)\Gamma(t, y) = (\mathbf{a}(s) \wedge \mathbf{a}(t))(F(x) \wedge F(y) - F(x)F(y)).$$

The integral on the right of (2.10) is understood as a mean-square limit of

$$\int_{\mathbb{R}_+^2} I_{l,t}(s,x)\,d\Gamma(s,x) = \sum_{i=1}^{l}(\Gamma(s_i^l, t-s_{i-1}^l) - \Gamma(s_{i-1}^l, t-s_{i-1}^l)) + \Gamma(0,t),$$

where $0 = s_0^l < s_1^l < s_2^l < \cdots < s_l^l = t$, $I_{l,t}(s,x) = \sum_{i=1}^{l} \mathbf{1}_{\{s\in(s_{i-1}^l, s_i^l]\}} \mathbf{1}_{\{0\le x\le t-s_{i-1}^l\}} + \mathbf{1}_{\{s=0\}}\mathbf{1}_{\{0\le x\le t\}}$, and $\max_i(s_i^l - s_{i-1}^l) \to 0$ as $l\to\infty$.

As mentioned, the trajectories of $X$ may have discontinuities of the second kind. To illustrate, suppose that $F(t) = \mathbf{1}_{\{t\ge 1\}}$, $\mathbf{q}$ is continuous, $\mathbf{q}(t) < 1$ for $t < 1$, $\mathbf{q}(2) = 1$ and, as $t\uparrow 2$, $\mathbf{q}(t)$ infinitely many times changes from being less than 1 to being greater than 1 and back. Suppose also that $\mathbf{e}$, $\tilde{F}$, and $Y$ are continuous. Then the function $\mathbf{a}$ is continuous, so $Z$ is continuous. Consequently, all terms in the equation for $X(t)$ except for the last one have limits as $t\uparrow 3$. The last term equals $X(t-1)$ when $\mathbf{q}(t-1) > 1$, equals $X(t-1)^+$ when $\mathbf{q}(t-1) = 1$, and equals 0 when $\mathbf{q}(t-1) < 1$. Therefore, there are three subsequential limits as $t\uparrow 3$: $X(2)$, $X(2)^+$, and 0, which shows that $X(3-)$ may be undefined. For example, we may take

$$\mathbf{e}(t) = \begin{cases} \dfrac{1}{1+\pi}(-(t-1)^2 + \pi t + 1), & 0\le t\le 1,\\ \dfrac{1}{1+\pi}\Big(-(t-2)^2 + \pi(t-1) \\ \qquad + 1 + (t-2)^2 \sin\dfrac{\pi}{2-t}\Big) + 1, & 1\le t < 2 \end{cases}$$

and $\mathbf{q}_0 = 0$. Then $\mathbf{q}(t) = \mathbf{e}(t) \le 1$ for $t\in[0,1]$, so

$$\mathbf{q}(t) = \mathbf{e}(t) - \mathbf{e}(t-1) + (\mathbf{q}(t-1) - 1)^+ = 1 + \frac{1}{1+\pi}(t-2)^2 \sin\frac{\pi}{2-t}$$

for $t\in[1,2)$.

As follows by Lemma B.1, if any interval $[0,T]$ can be partitioned into finitely many intervals such that on each of these intervals the function $\mathbf{q}$ either stays below 1, or stays at 1, or stays above 1, then the trajectories of $X$ admit left and right-hand limits. To show that even then $X$ does not necessarily have trajectories from $\mathbb{D}(\mathbb{R}_+,\mathbb{R})$, consider the example where $Q_0^n = 0$, $E^n$ is a Poisson process of rate $n$, and $F(t) = \mathbf{1}_{\{t\ge 2\}}$. Then $\mathbf{e}(t) = t$, $Y$ is a standard Wiener process, and $Z(t) = 0$. For the fluid limit, we have $\mathbf{q}(t) = t$ for $t < 2$, $\mathbf{q}(t) = 2$ for $2\le t < 3$, $\mathbf{q}(t) = t-1$ for $3\le t < 4$, $\mathbf{q}(t) = 3$ for $4\le t < 5$, $\mathbf{q}(t) = t-2$ for $5\le t < 6$, etc. The process $X$ is of the form: $X(t) = Y(t)$ for $t < 2$, $X(t) = Y(t) - Y(t-2)$ for $2\le t < 3$, $X(3) = Y(3) - Y(1)\wedge 0$, $X(t) = Y(t)$ for $3 < t < 4$, $X(t) = Y(t) - Y(t-4)$ for $4\le t < 5$, $X(5) = Y(5) - Y(1)\wedge 0$, $X(t) = Y(t)$ for $5 < t < 6$, etc. To sum up, there are two alternating patterns: $X(t) = Y(t)$ which occurs when $0 \le t \le 2$ or $2k+1 < t \le 2k+2$ for $k \in \mathbb{N}$ and $X(t) = Y(t) - Y(t-2k)$ which occurs when $2k \le t < 2k+1$ for $k \in \mathbb{N}$. We can also see that the trajectory $X(t), t \in \mathbb{R}_+$, is right-continuous when $Y(1) \ge 0$ and is left-continuous when $Y(1) \le 0$. To see that convergence in Skorohod's $J_1$-topology to the right-continuous version of $X$ does not hold for this example, note that [cf. (3.21) in the proof of Theorem 2.2]

$$X^n(t) = Y^n(t) - Y^n(t-2)\mathbf{1}_{\{t\ge 2\}} + ((X^n(t-2) + \sqrt{n}(\mathbf{q}(t-2)-1))^+ - \sqrt{n}(\mathbf{q}(t-2)-1)^+)\mathbf{1}_{\{t\ge 2\}}.$$

It follows that $\sup_{t\in[0,4]} |\Delta X^n(t)| \le 3 \sup_{t\in[0,4]} |\Delta Y^n(t)| = 3/\sqrt{n}$. On the other hand, $X(3+) - X(3-) = Y(1)$, so the jumps of $X^n$ "do not match" the jumps of $X$, which rules out convergence in distribution of the $X^n$ to the right-continuous version of $X$ for the $J_1$-topology. One can talk of convergence in the $M_1$-topology, see Whitt [24] for the definition and basic properties.

A sketch of the proof of the latter mode of convergence is as follows. It is obvious that the $X^n$ converge in distribution to $Y$ on $[0,2]$ for the topology of compact convergence. It is also easy to see that on $[2, 3-\varepsilon]$ and on $[3+\varepsilon, 4]$, for arbitrary $\varepsilon \in (0,1)$, there is convergence for the topology of compact convergence to $Y(t) - Y(t-2)$ and to $Y(t)$, respectively. Suppose $t \in [3-\varepsilon, 3+\varepsilon]$. Then on assuming that the $X^n$ converge to $Y$ in the topology of compact convergence on $[0,2]$ a.s. and noting that

$$\begin{aligned}\sup_{t\in[3-\varepsilon,3+\varepsilon]} &|((X^n(t-2) + \sqrt{n}(\mathbf{q}(t-2)-1))^+ - \sqrt{n}(\mathbf{q}(t-2)-1)^+) \\ &\quad - ((Y(1) + \sqrt{n}(\mathbf{q}(t-2)-1))^+ - \sqrt{n}(\mathbf{q}(t-2)-1)^+)| \\ &\le \sup_{t\in[0,2]} |X^n(t) - Y(t)| + \sup_{t\in[1-\varepsilon,1+\varepsilon]} |Y(t) - Y(1)|,\end{aligned}$$

we have by continuity of $Y$ that the latter left-hand side can be made arbitrarily small if $n$ is large enough and $\varepsilon$ is small enough. Now it is an easy matter to verify that the $(Y(1) + \sqrt{n}(\mathbf{q}(t-2)-1))^+ - \sqrt{n}(\mathbf{q}(t-2)-1)^+$ converge to $Y(1)\mathbf{1}_{\{t\in[3,3+\varepsilon]\}}$ as $n \to \infty$ in the $M_1$-topology on $[3-\varepsilon, 3+\varepsilon]$, cf., Whitt [24], page 82. It follows that the $X^n$ converge to the right-continuous version of $X$ for the $M_1$-topology on $[0,4]$. Because of the periodic pattern observed above, we derive $M_1$-convergence on $\mathbb{R}_+$. We conjecture that convergence in the $M_1$-topology holds more generally.

The Halfin–Whitt regime is obtained when $\mathbf{q}_0 = 1$, $\mathbf{e}(t) = \lambda t$ for some $\lambda > 0$, and $\tilde{F}(t) = \lambda \int_0^t (1 - F(s))\, ds$. By (2.3), $\mathbf{q}(t) = 1$ for $t \in \mathbb{R}_+$. The equation for $X$ in the statement of Theorem 2.2 takes the form

$$X(t) = X_0^+(1 - F(t)) + X_0 \wedge 0(1 - \tilde{F}(t)) + S(t) + Y(t) - \int_0^t Y(t-s)\, dF(s) + Z(t) + \int_0^t X(t-s)^+\, dF(s).$$

Usually, instead of there being a limit for the processes $Y^n$ as defined in (2.8) one assumes that $\int_0^\infty s\,dF(s) < \infty$ and that there exist $\rho^n \in (0,1)$ with $\sqrt{n}(1-\rho^n)$ converging to a finite limit as $n \to \infty$ such that the $(\sqrt{n}(E^n(t)/n - \rho^n t/\int_0^\infty s\,dF(s)), t \in \mathbb{R}_+)$ converge in distribution. These hypotheses imply convergence of the $Y^n$ with $\mathbf{e}(t) = t/(\int_0^\infty s\,dF(s))$.

We discuss condition (2.9). It obviously holds for the Halfin–Whitt regime. By Lemma B.3, it also holds when $F$ is a continuous function with $F(0) = 0$. If the function $F$ has atoms, then (2.9) holds if and only if $\mathbf{q}$ does not attain level 1 "in a continuous fashion," loosely speaking. More precisely, if and only if $\lim_{\varepsilon\to 0} t_\varepsilon = \infty$, where $t_\varepsilon = \inf\{t : |\mathbf{q}(t) - 1| \in (0,\varepsilon)\}$. To see that, note that if $t_0$ represents a jump time of $F$ and $t_\varepsilon < \infty$, then the left-hand side of (2.9) is greater than or equal to $\Delta F(t_0)$ for all $T > t_0 + t_\varepsilon$.

For instance, if $F$ contains atoms, $\mathbf{e}$ is continuous, $\mathbf{q}_0 < 1$, and there exists some $t \geq 0$ such that $\mathbf{q}(t) > 1$, then (2.9) cannot hold. This follows by the fact that if $\mathbf{e}$ is continuous, then by (2.5) $\mathbf{q}$ cannot have upward jumps. Therefore, upcrossings of level 1 occur in a continuous fashion. On the other hand, $\mathbf{q}$ may have downward jumps. If it starts above 1, then it may get below this level without spending time in small neighborhoods of 1, so it is possible for (2.9) to hold. Consider the following example. Suppose that $\mathbf{q}_0 = 2$ and $\mathbf{e}(t) = 0$ for $t \geq 0$ so that (2.3) becomes

$$\mathbf{q}(t) = 2 - F(t) - \tilde{F}(t) + \int_0^t (\mathbf{q}(t-s) - 1)^+ \,dF(s).$$

Furthermore, suppose that $F(t) = \tilde{F}(t) = \mathbf{1}_{\{t\geq 1\}}$. We then have that $\mathbf{q}(t) = 2$ for $0 \leq t < 1$, $\mathbf{q}(t) = 1$ for $1 \leq t < 2$, and $\mathbf{q}(t) = 0$ for $t \geq 2$.

2.3. *Extensions.* In this subsection, we state a number of results which admit similar proofs.

We begin by extending Theorems 2.1 and 2.2 to the case where the customers in the queue at time $0-$ and the customers exogenously arriving after time $0-$ may have differently distributed service times and where there can be available servers at time $0-$ even when there are customers awaiting service. Thus, $\eta_1, \eta_2, \ldots$ represent service times of the exogenously arriving customers. Condition (2.1) is still assumed. We denote the number of customers in the queue at time $0-$ by $\hat{Q}_0^n$. Their service times are denoted $\hat{\eta}_1, \hat{\eta}_2, \ldots$ and come from an i.i.d. sequence of nonnegative random variables with distribution $\hat{F}$. As with the $\eta_i$, it is assumed that the $\hat{\eta}_i$ are not equal to zero a.s., that is, $\hat{F}(0) < 1$. We retain the rest of the previous notation and assumptions. The entities $\tilde{Q}_0^n$, $\hat{Q}_0^n$, $\{\tilde{\eta}_1, \tilde{\eta}_2, \ldots\}$, $\{\hat{\eta}_1, \hat{\eta}_2, \ldots\}$, $\{\eta_1, \eta_2, \ldots\}$, and $E^n = (E^n(t), t \in \mathbb{R}_+)$ are assumed to be independent.

Let $\breve{A}^n(t)$ denote the number of customers that enter service by time $t$ out of those that have arrived after $0-$, let $\hat{A}^n(t)$ denote the number of customers that enter service by time $t$ out of those that were in the queue at time $0-$, and let $\hat{Q}^n(t)$ denote the number of customers out of those present at time $0-$ that have not left by time $t$. The introduced quantities satisfy the following equations for $t\in\mathbb{R}_+$:

$$\begin{aligned}
Q^n(t) &= \hat{Q}^n(t) + E^n(t) - \int_0^t\int_0^t \mathbf{1}_{\{s+x\le t\}}\, d\sum_{i=1}^{\breve{A}^n(s)} \mathbf{1}_{\{\eta_i\le x\}},\\
\breve{A}^n(t) &= E^n(t) + (\hat{Q}^n(t)-n)^+ - (Q^n(t)-n)^+,\\
\hat{Q}^n(t) &= \hat{Q}_0^n + \tilde{Q}^n(t) - \int_0^t\int_0^t \mathbf{1}_{\{s+x\le t\}}\, d\sum_{i=1}^{\hat{A}^n(s)} \mathbf{1}_{\{\hat{\eta}_i\le x\}},\\
\hat{A}^n(t) &= \hat{Q}_0^n - (\hat{Q}^n(t)-n)^+,\\
\tilde{Q}^n(t) &= \sum_{i=1}^{\tilde{Q}_0^n} \mathbf{1}_{\{\tilde{\eta}_i>t\}}.
\end{aligned}$$

The fluid limit equations are of the form

$$\begin{aligned}
\mathbf{q}(t) &= \mathbf{e}(t) - \int_0^t \mathbf{e}(t-s)\,dF(s) + \hat{\mathbf{q}}(t) - \int_0^t (\hat{\mathbf{q}}(t-s)-1)^+\,dF(s)\\
&\quad + \int_0^t (\mathbf{q}(t-s)-1)^+\,dF(s),\\
\hat{\mathbf{q}}(t) &= \hat{\mathbf{q}}_0(1-\hat{F}(t)) + \tilde{\mathbf{q}}_0(1-\tilde{F}(t)) + \int_0^t (\hat{\mathbf{q}}(t-s)-1)^+\,d\hat{F}(s).
\end{aligned}$$

Existence and uniqueness of solutions to these equations are addressed in Lemma A.1 and by Lemma B.1.

We also define functions $\breve{\mathbf{a}}=(\breve{\mathbf{a}}(t), t\in\mathbb{R}_+)\in\mathbb{D}(\mathbb{R}_+,\mathbb{R})$ and $\hat{\mathbf{a}}=(\hat{\mathbf{a}}(t), t\in\mathbb{R}_+)\in\mathbb{D}(\mathbb{R}_+,\mathbb{R})$ by the equalities

$$\begin{aligned}
\breve{\mathbf{a}}(t) &= \mathbf{e}(t) + (\hat{\mathbf{q}}(t)-1)^+ - (\mathbf{q}(t)-1)^+,\\
\hat{\mathbf{a}}(t) &= \hat{\mathbf{q}}_0 - (\hat{\mathbf{q}}(t)-1)^+.
\end{aligned}$$

THEOREM 2.3. *Suppose that, for arbitrary $T>0$ and $\varepsilon>0$,*

$$\begin{aligned}
\lim_{n\to\infty}\mathbf{P}\left(\left|\frac{\hat{Q}_0^n}{n}-\hat{\mathbf{q}}_0\right|>\varepsilon\right) &= 0,\\
\lim_{n\to\infty}\mathbf{P}\left(\left|\frac{\tilde{Q}_0^n}{n}-\tilde{\mathbf{q}}_0\right|>\varepsilon\right) &= 0,\\
\lim_{n\to\infty}\mathbf{P}\Big(\sup_{t\in[0,T]}|\tilde{F}^n(t)-\tilde{F}(t)|>\varepsilon\Big) &= 0
\end{aligned}$$

*and*

$$\lim_{n\to\infty}\mathbf{P}\biggl(\sup_{t\in[0,T]}\biggl|\frac{E^n(t)}{n}-\mathbf{e}(t)\biggr|>\varepsilon\biggr)=0.$$

*Then, for arbitrary* $T>0$ *and* $\varepsilon>0$,

$$\lim_{n\to\infty}\mathbf{P}\biggl(\sup_{t\in[0,T]}\biggl|\frac{Q^n(t)}{n}-\mathbf{q}(t)\biggr|>\varepsilon\biggr)=0,$$

$$\lim_{n\to\infty}\mathbf{P}\biggl(\sup_{t\in[0,T]}\biggl|\frac{\hat{Q}^n(t)}{n}-\hat{\mathbf{q}}(t)\biggr|>\varepsilon\biggr)=0,$$

$$\lim_{n\to\infty}\mathbf{P}\biggl(\sup_{t\in[0,T]}\biggl|\frac{\breve{A}^n(t)}{n}-\breve{\mathbf{a}}(t)\biggr|>\varepsilon\biggr)=0$$

*and*

$$\lim_{n\to\infty}\mathbf{P}\biggl(\sup_{t\in[0,T]}\biggl|\frac{\hat{A}^n(t)}{n}-\hat{\mathbf{a}}(t)\biggr|>\varepsilon\biggr)=0.$$

For a stochastic approximation result, we assume that the processes $X^n=(X^n(t),t\in\mathbb{R}_+)$, $S^n=(S^n(t),t\in\mathbb{R}_+)$, and $Y^n=(Y^n(t),t\in\mathbb{R}_+)$ are defined by (2.6)–(2.8), respectively. We also denote

$$\hat{X}_0^n=\sqrt{n}\biggl(\frac{1}{n}\hat{Q}_0^n-\hat{\mathbf{q}}_0\biggr),\qquad \tilde{X}_0^n=\sqrt{n}\biggl(\frac{1}{n}\tilde{Q}_0^n-\tilde{\mathbf{q}}_0\biggr).$$

THEOREM 2.4. *Assume that the* $\hat{X}_0^n$ *and* $\tilde{X}_0^n$ *converge in distribution in* $\mathbb{R}$ *to random variables* $\hat{X}_0$ *and* $\tilde{X}_0$, *respectively, as* $n\to\infty$. *If the* $S^n$ *converge in distribution in* $\mathbb{D}_c(\mathbb{R}_+,\mathbb{R})$ *to a process* $S=(S(t),t\in\mathbb{R}_+)$ *and the* $Y^n$ *converge in distribution in* $\mathbb{D}_c(\mathbb{R}_+,\mathbb{R})$ *to a process* $Y=(Y(t),t\in\mathbb{R}_+)$ *such that* $S$ *and* $Y$ *are separable random elements of* $\mathbb{D}_c(\mathbb{R}_+,\mathbb{R})$, *then the processes* $X^n$ *converge in the sense of finite-dimensional distributions to the process* $X=(X(t),t\in\mathbb{R}_+)$ *that is a unique strong solution to the set of equations*

$$\begin{aligned}X(t)&=Y(t)-\int_0^t Y(t-s)\,dF(s)+\breve{Z}(t)+\hat{X}(t)\\&\quad-\int_0^t(\hat{X}(t-s)\mathbf{1}_{\{\hat{\mathbf{q}}(t-s)>1\}}+\hat{X}(t-s)^+\mathbf{1}_{\{\hat{\mathbf{q}}(t-s)=1\}})\,dF(s)\\&\quad+\int_0^t(X(t-s)\mathbf{1}_{\{\mathbf{q}(t-s)>1\}}+X(t-s)^+\mathbf{1}_{\{\mathbf{q}(t-s)=1\}})\,dF(s)\end{aligned}$$

*and*

$$\hat{X}(t) = \hat{X}_0(1 - \hat{F}(t)) + \tilde{X}_0(1 - \tilde{F}(t)) + \sqrt{\tilde{\mathbf{q}}_0} S(t) + \hat{Z}(t) + \int_0^t (\hat{X}(t-s)\mathbf{1}_{\{\hat{\mathbf{q}}(t-s)>1\}} + \hat{X}(t-s)^+ \mathbf{1}_{\{\hat{\mathbf{q}}(t-s)=1\}})\, d\hat{F}(s),$$

*where* $\breve{Z} = (\breve{Z}(t), t \in \mathbb{R}_+)$ *and* $\hat{Z} = (\hat{Z}(t), t \in \mathbb{R}_+)$ *are zero-mean Gaussian semimartingales with trajectories in* $\mathbb{D}(\mathbb{R}_+, \mathbb{R})$ *and with respective covariances*

$$\mathbf{E}\breve{Z}(s)\breve{Z}(t) = \int_0^{s\wedge t} F(s \wedge t - u)(1 - F(s \vee t - u))\, d\breve{\mathbf{a}}(u)$$

*and*

$$\mathbf{E}\hat{Z}(s)\hat{Z}(t) = \int_0^{s\wedge t} \hat{F}(s \wedge t - u)(1 - \hat{F}(s \vee t - u))\, d\hat{\mathbf{a}}(u),$$

*the entities* $\hat{X}_0$, $\tilde{X}_0$, $\breve{Z}$, $\hat{Z}$, $S$ *and* $Y$ *being independent. The trajectories of* $X$ *are Borel measurable and locally bounded a.s.*

*If, in addition, for all* $T > 0$,

$$\lim_{\varepsilon\to 0} \sup_{t\in[0,T]} \int_0^t (\mathbf{1}_{\{0<|\hat{\mathbf{q}}(t-s)-1|<\varepsilon\}} + \mathbf{1}_{\{0<|\mathbf{q}(t-s)-1|<\varepsilon\}})\, dF(s) = 0$$

*and*

$$\lim_{\varepsilon\to 0} \sup_{t\in[0,T]} \int_0^t \mathbf{1}_{\{0<|\hat{\mathbf{q}}(t-s)-1|<\varepsilon\}}\, d\hat{F}(s) = 0,$$

*then* $X$ *and* $\hat{X}$ *are separable random elements of* $\mathbb{D}_c(\mathbb{R}_+, \mathbb{R})$ *and the* $X^n$ *converge in distribution in* $\mathbb{D}_c(\mathbb{R}_+, \mathbb{R})$ *to* $X$.

We conclude this subsection with the Gaussian approximation result for an infinite server. Consider a sequence of $G/GI/\infty$ systems indexed by $n$. Adapting the earlier notation, we denote by $\overline{Q}^n(t)$ the number of customers present at time $t$ and we denote by $E^n(t)$ the number of exogenous arrivals by $t$. We also reuse the introduced earlier sequences $\{\tilde{\eta}_i, i \in \mathbb{N}\}$ and $\{\eta_i, i \in \mathbb{N}\}$. As above, they represent the remaining service times of customers in service at time $0-$, whose number is denoted by $\tilde{Q}_0^n$, and the service times of exogenously arriving customers, respectively. The $\eta_i$ are assumed to be i.i.d. with distribution $F$, however, by contrast with condition (2.1) we allow these random variables to equal zero a.s. The entities $E^n = (E^n(t), t \in \mathbb{R}_+)$, $\tilde{Q}_0^n$, $\{\tilde{\eta}_i, i \in \mathbb{R}_+\}$ and $\{\eta_i, i \in \mathbb{R}_+\}$ are assumed to be independent. The evolution of the process $\overline{Q}^n = (\overline{Q}^n(t), t \in \mathbb{R}_+)$ is governed by the equation

$$\overline{Q}^n(t) = \sum_{i=1}^{\tilde{Q}_0^n} \mathbf{1}_{\{\tilde{\eta}_i > t\}} + E^n(t) - \int_0^t \int_0^t \mathbf{1}_{\{s+x\le t\}}\, d \sum_{i=1}^{E^n(s)} \mathbf{1}_{\{\eta_i \le x\}}.$$

It specifies the process $\overline{Q}^n$ uniquely. In analogy with the earlier notation, given $\mathbf{q}_0 \in \mathbb{R}_+$ and a nondecreasing function $\mathbf{e} = (\mathbf{e}(t), t \in \mathbb{R}_+) \in \mathbb{D}(\mathbb{R}_+, \mathbb{R})$, we define $\overline{\mathbf{q}}(t) = \mathbf{q}_0(1 - \tilde{F}(t)) + \mathbf{e}(t) - \int_0^t \mathbf{e}(t-s)\,dF(s)$, $\overline{X}^n(t) = \sqrt{n}(\overline{Q}^n(t)/n - \overline{\mathbf{q}}(t))$, $X_0^n = \sqrt{n}(\tilde{Q}_0^n/n - \mathbf{q}_0)$, and we define $S^n(t)$ and $Y^n(t)$ by (2.7) and (2.8), respectively.

THEOREM 2.5. *Suppose that, as $n \to \infty$, the $X_0^n$ converge in distribution in $\mathbb{R}$ to a random variable $X_0$, the processes $S^n = (S^n(t), t \in \mathbb{R}_+)$ converge in distribution in $\mathbb{D}_c(\mathbb{R}_+, \mathbb{R})$ to a process $S = (S(t), t \in \mathbb{R}_+)$, and the processes $Y^n = (Y^n(t), t \in \mathbb{R}_+)$ converge in distribution in $\mathbb{D}_c(\mathbb{R}_+, \mathbb{R})$ to a process $Y = (Y(t), t \in \mathbb{R}_+)$, where $S$ and $Y$ are separable random elements of $\mathbb{D}_c(\mathbb{R}_+, \mathbb{R})$. Then the processes $\overline{X}^n = (\overline{X}^n(t), t \in \mathbb{R}_+)$ converge in distribution in $\mathbb{D}_c(\mathbb{R}_+, \mathbb{R})$ to the process $\overline{X} = (\overline{X}(t), t \in \mathbb{R}_+)$ defined by*

$$\overline{X}(t) = (1 - \tilde{F}(t))X_0 + \sqrt{\mathbf{q}_0}S(t) + Y(t) - \int_0^t Y(t-s)\,dF(s) + \overline{Z}(t),$$

*where $\overline{Z} = (\overline{Z}(t), t \in \mathbb{R}_+)$ is a zero-mean Gaussian semimartingale with trajectories in $\mathbb{D}(\mathbb{R}_+, \mathbb{R})$ and with covariance*

$$\mathbf{E}\overline{Z}(s)\overline{Z}(t) = \int_0^{s \wedge t} F(s \wedge t - u)(1 - F(s \vee t - u))\,d\mathbf{e}(u),$$

*the entities $X_0$, $S$, $\overline{Z}$ and $Y$ being independent.*

For the case of a continuous $Y$, the latter limit theorem generalizes the earlier results by Borovkov [2] (which result can also be found in Borovkov [3], Chapter 2, Section 2) and Krichagina and Puhalskii [15]. We recall that Borovkov [2] imposed a Hölder continuity condition on the function $(\int_0^t (1 - F(t-s))\,d\mathbf{e}(s), t \in \mathbb{R}_+)$ and Krichagina and Puhalskii [15] required $\mathbf{e}$ to be continuous.

**3. Proofs of the main results.** This section contains proofs of Theorems 2.1 and 2.2.

3.1. *Preliminaries.* We start by developing semimartingale representations for certain processes. On introducing

$$V^n(t,x) = \frac{1}{\sqrt{n}} \sum_{i=1}^{A^n(t)} (\mathbf{1}_{\{\eta_i \le x\}} - F(x)), \tag{3.1}$$

we have by (2.2a)

$$\frac{1}{n}Q^n(t) = \left(\frac{1}{n}Q_0^n - 1\right)^+ + \frac{1}{n}\tilde{Q}^n(t)$$

$$+\frac{1}{n}E^n(t)-\frac{1}{n}\int_0^t A^n(t-s)\,dF(s)+\frac{1}{\sqrt{n}}Z^n(t), \tag{3.2}$$

where

$$Z^n(t)=-\int_{\mathbb{R}_+^2}\mathbf{1}_{\{s+x\le t\}}\,dV^n(s,x). \tag{3.3}$$

Let

$$L^n(t,x)=\frac{1}{\sqrt{n}}\sum_{i=1}^{A^n(t)}\left(\mathbf{1}_{\{\eta_i\le x\}}-\int_0^{\eta_i\wedge x}\frac{dF(u)}{1-F(u-)}\right). \tag{3.4}$$

Then by (3.1),

$$V^n(t,x)=-\int_0^x\frac{V^n(t,u-)}{1-F(u-)}\,dF(u)+L^n(t,x). \tag{3.5}$$

Hence, by (3.3),

$$Z^n(t)=G^n(t)-M^n(t), \tag{3.6}$$

where

$$G^n(t)=\int_0^t\frac{V^n(t-x,x-)}{1-F(x-)}\,dF(x) \tag{3.7}$$

and

$$M^n(t)=\int_{\mathbb{R}_+^2}\mathbf{1}_{\{s+x\le t\}}\,dL^n(s,x). \tag{3.8}$$

We note that by (3.1), (3.4) and (3.8)

$$M^n(t)=V^n(t,0)+\int_{\mathbb{R}_+^2}\mathbf{1}_{\{x>0\}}\mathbf{1}_{\{s+x\le t\}}\,dL^n(s,x). \tag{3.9}$$

We also define, for $k\in\mathbb{N}$ and $t\in\mathbb{R}_+$,

$$M_k^n(t)=\int_{\mathbb{R}_+^2}\mathbf{1}_{\{s+x\le t\}}\,dL_k^n(s,x), \tag{3.10}$$

where

$$L_k^n(s,x)=\frac{1}{\sqrt{n}}\sum_{i=1}^{A^n(s)\wedge k}\left(\mathbf{1}_{\{0<\eta_i\le x\}}-\int_0^{\eta_i\wedge x}\mathbf{1}_{\{u>0\}}\frac{dF(u)}{1-F(u-)}\right), \tag{3.11}$$

and

$$\langle M_k^n\rangle(t)=\int_{\mathbb{R}_+^2}\mathbf{1}_{\{s+x\le t\}}\,d\langle L_k^n\rangle(s,x),$$

where

$$\langle L_k^n\rangle(s,x) = \frac{1}{n}\sum_{i=1}^{A^n(s)\wedge k}\int_0^{\eta_i\wedge x}\mathbf{1}_{\{u>0\}}\frac{1-F(u)}{(1-F(u-))^2}\,dF(u).$$

For $t\in\mathbb{R}_+$, let $\hat{\mathcal{G}}_t^n$ denote the complete $\sigma$-algebra generated by the random variables $\mathbf{1}_{\{\tau_i^n\le s\}}\mathbf{1}_{\{\eta_i\le x\}}$, where $x+s\le t$ and $i\in\mathbb{N}$, and by the $A^n(s)$ (or, equivalently, by the $\mathbf{1}_{\{\tau_i^n\le s\}}$ for $i\in\mathbb{N}$), where $s\le t$. We introduce "the right-continuous modification" by $\mathcal{G}_t^n=\bigcap_{\varepsilon>0}\hat{\mathcal{G}}_{t+\varepsilon}^n$. Then $\mathbf{G}^n=(\mathcal{G}_t^n,t\in\mathbb{R}_+)$ is a filtration.

The next lemma originates in Krichagina and Puhalskii [15], Lemma 3.5, see also Reed [22], Lemma 1, where larger filtrations are used. Our conditions are closer to minimal ones. The proof is also different from those in the cited papers and is more direct.

LEMMA 3.1. *For each $k\in\mathbb{N}$, the process $M_k^n=(M_k^n(t),t\in\mathbb{R}_+)$ is a $\mathbf{G}^n$-square integrable martingale starting at $0$ with predictable quadratic variation process $\langle M_k^n\rangle=(\langle M_k^n\rangle(t),t\in\mathbb{R}_+)$.*

PROOF. According to Lemma C.3 it suffices to prove that, given $i\in\mathbb{N}$, $s\in\mathbb{R}_+$ and $x\in\mathbb{R}_+$, the random variable $\eta_i$ is independent of the random variables $\eta_j$ for $j\ne i$, of the $\tau_j^n$ for $j\le i$, and of the $\tau_j^n\wedge(s+x)$ for $j>i$, when conditioned on the event $\{\tau_i^n\ge s,\eta_i>x\}$ and that $\eta_i$ is independent of the $\tau_j^n$ for $j\le i$. By (2.2a)–(2.2c), the $\tau_j^n$ for $j\le i$ are measurable with respect to the $\sigma$-algebra generated by $\tilde{Q}_0^n$, $\hat{Q}_0^n$, $\tilde{\eta}_j,j\in\mathbb{N}$, $\eta_j,j<i$, and $E^n(t),t\in\mathbb{R}_+$. By the independence assumptions in the hypotheses, this $\sigma$-algebra and $\eta_i$ are independent, hence, $\tau_j^n$ for $j\le i$ and $\eta_i$ are independent.

Let $A^{n,i}=(A^{n,i}(t),t\in\mathbb{R}_+)$ denote the process of exogenous arrivals entering service that would occur if the $i$th exogenously arriving customer had an infinite service time and let $\tau_j^{n,i}=\inf\{t:A^{n,i}(t)\ge j\}$. Then by (2.2a)–(2.2c) the $\tau_j^{n,i}$ for $j\in\mathbb{N}$ are measurable with respect to the $\sigma$-algebra generated by $\tilde{Q}_0^n$, $\hat{Q}_0^n$, $\tilde{\eta}_j,j\in\mathbb{N}$, $\eta_j,j<i$, and $E^n(t),t\in\mathbb{R}_+$. Hence, they are independent of $\eta_i$. Therefore, $\eta_i$, on the one hand, and $\eta_j$ for $j\ne i$, $\tau_j^n$ for $j\le i$, and $\tau_j^{n,i}$ for $j>i$, on the other hand, are independent. In addition, as it follows by (2.2a)–(2.2c), on the event $\{\tau_i^n\ge s,\eta_i>x\}$ the $i$th exogenously arriving customer does not depart until after time $s+x$, which means that she has the same effect on the epochs before time $s+x$ when exogenously arriving customers enter service as if she never completed service. To put it precisely, $\tau_j^{n,i}\wedge(s+x)=\tau_j^n\wedge(s+x)$ for $j\in\mathbb{N}$ on the event $\{\tau_i^n\ge s,\eta_i>x\}$.

By the established independence property, we have for natural numbers $r_1,r_2,\dots,r_l$ none of which equals $i$ and $i<p_1<p_2<\cdots<p_m$, and for bounded Borel functions $f_1,\dots,f_l$, $g_1,\dots,g_i$, $h_1,\dots,h_m$, and $f$, on assuming that conditional probabilities given events of probability zero equal zero and adopting the convention that $0/0=0$,

$$\begin{aligned}
&\mathbf{E}\Biggl(\prod_{j=1}^{l} f_j(\eta_{r_j})\prod_{j=1}^{i} g_j(\tau_j^n)\prod_{j=1}^{m} h_j(\tau_{p_j}^n\wedge(s+x))f(\eta_i)\Bigm| \tau_i^n\ge s,\eta_i>x\Biggr)\\
&\quad=\mathbf{E}\Biggl(\prod_{j=1}^{l} f_j(\eta_{r_j})\prod_{j=1}^{i} g_j(\tau_j^n)\\
&\qquad\quad\times\prod_{j=1}^{m} h_j(\tau_{p_j}^n\wedge(s+x))f(\eta_i)\mathbf{1}_{\{\tau_i^n\ge s,\eta_i>x\}}\Biggr)\Big/\mathbf{P}(\tau_i^n\ge s,\eta_i>x)\\
&\quad=\mathbf{E}\Biggl(\prod_{j=1}^{l} f_j(\eta_{r_j})\prod_{j=1}^{i} g_j(\tau_j^n)\\
&\qquad\quad\times\prod_{j=1}^{m} h_j(\tau_{p_j}^{n,i}\wedge(s+x))f(\eta_i)\mathbf{1}_{\{\tau_i^n\ge s,\eta_i>x\}}\Biggr)\Big/\mathbf{P}(\tau_i^n\ge s,\eta_i>x)\\
&\quad=\mathbf{E}\Biggl(\prod_{j=1}^{l} f_j(\eta_{r_j})\prod_{j=1}^{i} g_j(\tau_j^n)\prod_{j=1}^{m} h_j(\tau_{p_j}^{n,i}\wedge(s+x))\mathbf{1}_{\{\tau_i^n\ge s\}}\Biggr)\\
&\qquad\times\mathbf{E}(f(\eta_i)\mathbf{1}_{\{\eta_i>x\}})\Big/\mathbf{P}(\tau_i^n\ge s)\mathbf{P}(\eta_i>x)\\
&\quad=\frac{\mathbf{E}(\prod_{j=1}^{l} f_j(\eta_{r_j})\prod_{j=1}^{i} g_j(\tau_j^n)\prod_{j=1}^{m} h_j(\tau_{p_j}^{n,i}\wedge(s+x))\mathbf{1}_{\{\tau_i^n\ge s\}}\mathbf{1}_{\{\eta_i>x\}})}{\mathbf{P}(\tau_i^n\ge s,\eta_i>x)}\\
&\qquad\times\frac{\mathbf{E}(f(\eta_i)\mathbf{1}_{\{\eta_i>x\}}\mathbf{1}_{\{\tau_i^n\ge s\}})}{\mathbf{P}(\tau_i^n\ge s,\eta_i>x)}\\
&\quad=\frac{\mathbf{E}(\prod_{j=1}^{l} f_j(\eta_{r_j})\prod_{j=1}^{i} g_j(\tau_j^n)\prod_{j=1}^{m} h_j(\tau_{p_j}^{n}\wedge(s+x))\mathbf{1}_{\{\tau_i^n\ge s\}}\mathbf{1}_{\{\eta_i>x\}})}{\mathbf{P}(\tau_i^n\ge s,\eta_i>x)}\\
&\qquad\times\frac{\mathbf{E}(f(\eta_i)\mathbf{1}_{\{\eta_i>x\}}\mathbf{1}_{\{\tau_i^n\ge s\}})}{\mathbf{P}(\tau_i^n\ge s,\eta_i>x)}\\
&\quad=\mathbf{E}\Biggl(\prod_{j=1}^{l} f_j(\eta_{r_j})\prod_{j=1}^{i} g_j(\tau_j^n)\prod_{j=1}^{m} h_j(\tau_{p_j}^n\wedge(s+x))\Bigm| \tau_i^n\ge s,\eta_i>x\Biggr)\\
&\qquad\times\mathbf{E}(f(\eta_i)\mid\tau_i^n\ge s,\eta_i>x).
\end{aligned}$$

An application of Lemma C.3 completes the proof. □

REMARK 3.1. We stop short of claiming that $M^n$ is a $\mathbf{G}^n$-locally square integrable martingale because there might not be a sequence of stopping times such that the "stopped" processes $A^n$ are bounded. However, if, in addition, the jumps of $A^n$ (including the "jump" at zero) are uniformly bounded, then $M^n$ is a $\mathbf{G}^n$-locally square integrable martingale.

3.2. *Proof of Theorem 2.1.*

LEMMA 3.2. *Under the hypotheses of Theorem 2.1, for $T > 0$,*

$$\lim_{b\to\infty} \limsup_{n\to\infty} \mathbf{P}\Big(\sup_{t\in[0,T]} |Z^n(t)| > b\Big) = 0.$$

PROOF. In view of (3.6), it suffices to prove that

$$\lim_{b\to\infty} \limsup_{n\to\infty} \mathbf{P}\Big(\sup_{t\in[0,T]} |G^n(t)| > b\Big) = 0 \tag{3.12}$$

and

$$\lim_{b\to\infty} \limsup_{n\to\infty} \mathbf{P}\Big(\sup_{t\in[0,T]} |M^n(t)| > b\Big) = 0. \tag{3.13}$$

Let us firstly note that, for $t > 0$ and $b$ large enough,

$$\lim_{n\to\infty} \mathbf{P}\bigg(\frac{1}{n}A^n(t) > b\bigg) = 0, \tag{3.14}$$

which follows by the bound $A^n(t) \le (Q_0^n - n)^+ + E^n(t)$ [see (2.2b)] and the convergence in probability of $Q_0^n/n$ to $\mathbf{q}_0$ and of $E^n(t)/n$ to $\mathbf{e}(t)$ in the hypotheses of Theorem 2.1.

The proof of (3.12) is analogous to the proof of (3.23) in Krichagina and Puhalskii [15]. Let

$$\breve{V}^n(t,x) = \frac{1}{\sqrt{n}} \sum_{i=1}^{\lfloor nt\rfloor} (\mathbf{1}_{\{\eta_i \le x\}} - F(x)). \tag{3.15}$$

We have, for $c > 0$, on taking into account (3.1), (3.7) and (3.15),

$$\mathbf{P}\Big(\sup_{t\in[0,T]} |G^n(t)| > b\Big) \le \mathbf{P}\bigg(\frac{1}{n}A^n(T) > cT\bigg) + \mathbf{P}\bigg(\int_0^\infty \frac{\sup_{t\in[0,cT]} |\breve{V}^n(t,x-)|}{1-F(x-)}\, dF(x) > b\bigg).$$

By (3.15), for fixed $x$, $\breve{V}^n(t,x)$ is a locally square-integrable martingale in $t$ relative to the natural filtration with predictable quadratic variation process $(\lfloor nt\rfloor/n)F(x)(1-F(x))$. Theorem 1.9.5 in Liptser and Shiryayev [16] yields the bound $\mathbf{E}\sup_{t\in[0,cT]}|\breve{V}^n(t,x-)|\le 3\sqrt{cT(1-F(x-))}$. Hence,

$$\int_0^\infty \frac{\mathbf{E}\sup_{t\in[0,cT]}|\breve{V}^n(t,F(x-))|}{1-F(x-)}\,dF(x)\le 6\sqrt{cT},$$

so (3.12) follows by an application of Markov's inequality and (3.14).

Next, on noting that $M^n(t)=\breve{V}^n(A^n(t)/n,0)+M_k^n(t)$ when $A^n(t)\le k$, we have by (3.9), Lemma 3.1, Kolmogorov's inequality, and the Lenglart–Rebolledo inequality; see, for example, Liptser and Shiryayev [16], Theorem 1.9.3, for $b>0$ and $c>0$,

$$\begin{aligned}
\mathbf{P}\Bigl(\sup_{t\in[0,T]}|M^n(t)|>b\Bigr)&\le \mathbf{P}(A^n(T)>\lfloor nbT\rfloor)+\mathbf{P}\Bigl(\sup_{t\in[0,bT]}|\breve{V}^n(t,0)|>\frac{b}{2}\Bigr)\\
&\quad+\mathbf{P}\Bigl(\sup_{t\in[0,T]}|M^n_{\lfloor nbT\rfloor}(t)|>\frac{b}{2}\Bigr)\\
&\le \mathbf{P}(A^n(T)>\lfloor nbT\rfloor)+\frac{4(T+c)}{b^2}+\mathbf{P}(\langle M^n_{\lfloor nbT\rfloor}\rangle(T)>c)\\
&\le \mathbf{P}(A^n(T)>\lfloor nbT\rfloor)+\frac{4(T+c)}{b^2}\\
&\quad+\mathbf{P}\Biggl(\frac{1}{n}\sum_{i=1}^{\lfloor nbT\rfloor}\int_0^\infty \mathbf{1}_{\{u\le\eta_i\}}\frac{1-F(u)}{(1-F(u-))^2}\,dF(u)>c\Biggr).
\end{aligned}$$

By Markov's inequality,

$$\mathbf{P}\Biggl(\frac{1}{n}\sum_{i=1}^{\lfloor nbT\rfloor}\int_0^\infty \mathbf{1}_{\{u\le\eta_i\}}\frac{1-F(u)}{(1-F(u-))^2}\,dF(u)>c\Biggr)\le\frac{bT}{c},$$

so, on picking $c=b^{3/2}T^{1/2}$,

$$\mathbf{P}\Bigl(\sup_{t\in[0,T]}|M^n(t)|>b\Bigr)\le \mathbf{P}(A^n(T)>\lfloor nbT\rfloor)+\frac{4T}{b^2}+5\sqrt{\frac{T}{b}}.$$

The convergence in (3.13) now follows by (3.14). □

Proof of Theorem 2.1. Denote $\mathbf{e}^n(t)=E^n(t)/n$, $\mathbf{q}_0^n=Q_0^n/n$, $\mathbf{q}^n(t)=Q^n(t)/n$, and $\tilde{\mathbf{q}}^n(t)=\tilde{Q}^n(t)/n$ so that by (2.2b) and (3.2),

$$\begin{aligned}
\mathbf{q}^n(t)&=(\mathbf{q}_0^n-1)^+ +\tilde{\mathbf{q}}^n(t)+\mathbf{e}^n(t)-\int_0^t \mathbf{e}^n(t-s)\,dF(s)\\
&\quad+\int_0^t(\mathbf{q}^n(t-s)-1)^+\,dF(s)+\frac{1}{\sqrt{n}}Z^n(t).
\end{aligned}\tag{3.16}$$

By Lemma 3.2, the $\sup_{t\in[0,T]}|Z^n(t)|/\sqrt{n}$ converge to zero in probability for every $T>0$ as $n\to\infty$. By hypotheses, the $(\mathbf{q}_0^n-1)^+$ and $\mathbf{q}_0^n\wedge 1$ converge in probability to $(\mathbf{q}_0-1)^+$ and $\mathbf{q}_0\wedge 1$, respectively. By (2.2c) and the hypotheses, the $\tilde{\mathbf{q}}^n(t)$ converge in probability uniformly over compact intervals of $t$ to $\mathbf{q}_0\wedge 1(1-\tilde{F}(t))$. Also, the compact convergence in probability of the $\mathbf{e}^n$ to $\mathbf{e}$ implies that

$$\sup_{t\in[0,T]}\left|\mathbf{e}^n(t)-\int_0^t \mathbf{e}^n(t-s)\,dF(s)-\left(\mathbf{e}(t)-\int_0^t \mathbf{e}(t-s)\,dF(s)\right)\right|$$

converges to zero in probability. Hence, on applying part 1 of Lemma B.2 to (3.17), we conclude by the fact that a sequence of random variables converges in probability if and only if its every subsequence contains a further subsequence that converges a.s. that the sequence $\sup_{t\in[0,T]}|\mathbf{q}^n(t)-\mathbf{q}(t)|$ converges in probability to zero as required.

The convergence

$$\lim_{n\to\infty}\mathbf{P}\left(\sup_{t\in[0,T]}\left|\frac{A^n(t)}{n}-\mathbf{a}(t)\right|>\varepsilon\right)=0$$

follows by (2.2b) and the part of the theorem already proven. □

3.3. *Proof of Theorem 2.2.* The key to the proof of Theorem 2.2 is the following result whose proof is deferred until the next section. Denote

$$H(t)=Y(t)-\int_0^t Y(t-s)\,dF(s), \tag{3.17}$$

$$H^n(t)=Y^n(t)-\int_0^t Y^n(t-s)\,dF(s), \tag{3.18}$$

$$\tilde{X}^n(t)=\sqrt{n}\left(\frac{1}{n}\tilde{Q}^n(t)-\mathbf{q}_0\wedge 1(1-\tilde{F}(t))\right), \tag{3.19}$$

$$\tilde{X}(t)=(X_0\mathbf{1}_{\{\mathbf{q}_0<1\}}+X_0\wedge 0\mathbf{1}_{\{\mathbf{q}_0=1\}})(1-\tilde{F}(t))+\sqrt{\mathbf{q}_0\wedge 1}S(t). \tag{3.20}$$

Let $H=(H(t),t\in\mathbb{R}_+)$, $H^n=(H^n(t),t\in\mathbb{R}_+)$, $\tilde{X}=(\tilde{X}(t),t\in\mathbb{R}_+)$, $\tilde{X}^n=(\tilde{X}^n(t),t\in\mathbb{R}_+)$, and $Z^n=(Z^n(t),t\in\mathbb{R}_+)$. We also denote

$$\hat{X}_0^n=\sqrt{n}\left(\frac{1}{n}(Q_0^n-n)^+-(\mathbf{q}_0-1)^+\right),$$

$$\hat{X}_0=X_0\mathbf{1}_{\{\mathbf{q}_0>1\}}+X_0^+\mathbf{1}_{\{\mathbf{q}_0=1\}}.$$

THEOREM 3.1. *Under the hypotheses of Theorem 2.2, the processes $\tilde{X}$, $H$ and $Z$ are separable random elements of $\mathbb{D}_c(\mathbb{R}_+,\mathbb{R})$. As $n\to\infty$, the $(\hat{X}_0^n,\tilde{X}^n,H^n,Z^n)$ converge in distribution in $\mathbb{R}\times\mathbb{D}_c(\mathbb{R}_+,\mathbb{R})^3$ to $(\hat{X}_0,\tilde{X},H,Z)$.*

Given Theorem 3.1 and Lemmas B.1 and B.3, the proof of Theorem 2.2 is now routine. On recalling that $X^n(t) = \sqrt{n}(Q^n(t)/n - \mathbf{q}(t))$, we have by (2.2b), (2.3), (2.4), (3.2), (3.18) and (3.19), that

$$\begin{aligned} X^n(t) &= \hat{X}_0^n(1 - F(t)) + \tilde{X}^n(t) + H^n(t) + Z^n(t) \\ &\quad + \int_0^t ((X^n(t-s) + \sqrt{n}(\mathbf{q}(t-s) - 1))^+ \\ &\qquad - \sqrt{n}(\mathbf{q}(t-s) - 1)^+ )\, dF(s). \end{aligned} \tag{3.21}$$

On writing (3.21) as $X^n = \Psi^n(\hat{X}_0^n, \tilde{X}^n, H^n, Z^n)$, we have by applying part 2 of Lemma B.2 with $f^n(y,t) = (y + \sqrt{n}(\mathbf{q}(t) - 1))^+ - \sqrt{n}(\mathbf{q}(t) - 1)^+$ and $f(y,t) = y\mathbf{1}_{\{\mathbf{q}(t)>1\}} + y^+\mathbf{1}_{\{\mathbf{q}(t)=1\}}$ that if $\mathbf{x}^n \to \mathbf{x}$ in $\mathbb{R} \times \mathbb{D}_c(\mathbb{R}_+, \mathbb{R})^3$, then $\Psi^n(\mathbf{x}^n)(t) \to \mathbf{y}(t)$ for all $t$, where $\mathbf{y} = (\mathbf{y}(t), t \in \mathbb{R}_+)$ is determined by the following equations assuming that $\mathbf{x} = (\mathbf{x}_1, \mathbf{x}_2, \mathbf{x}_3, \mathbf{x}_4)$ with $\mathbf{x}_1 \in \mathbb{R}$ and $\mathbf{x}_i = (\mathbf{x}_i(t), t \in \mathbb{R}_+)$ for $i = 2, 3, 4$:

$$\begin{aligned} \mathbf{y}(t) &= \mathbf{x}_1(1 - F(t)) + \mathbf{x}_2(t) + \mathbf{x}_3(t) + \mathbf{x}_4(t) \\ &\quad + \int_0^t (\mathbf{y}(t-s)\mathbf{1}_{\{\mathbf{q}(t-s)>1\}} + \mathbf{y}(t-s)^+\mathbf{1}_{\{\mathbf{q}(t-s)=1\}})\, dF(s). \end{aligned}$$

These equations specify $\mathbf{y}$ uniquely by Lemma B.1. The hypotheses and Theorem 3.1 imply that, as $n \to \infty$, the $(\hat{X}_0^n, \tilde{X}^n, H^n, Z^n)$ converge in distribution in $\mathbb{R} \times \mathbb{D}_c(\mathbb{R}_+, \mathbb{R})^3$ to $(\hat{X}_0, \tilde{X}, H, Z)$, which is a separable random element. Therefore, by the continuous mapping principle (see Theorem D.2), the $X^n$ converge in the sense of weak convergence of finite-dimensional distributions to $X = \Psi(\hat{X}_0, \tilde{X}, H, Z)$, where $\Psi(\mathbf{x}_1, \mathbf{x}_2, \mathbf{x}_3, \mathbf{x}_4) = \mathbf{y}$. The trajectories of $X$ are a.s. Borel measurable and locally bounded by Lemma B.1.

If, in addition, condition (2.9) holds, then, by Lemmas B.2 and B.3, $X$ has trajectories in $\mathbb{D}(\mathbb{R}_+, \mathbb{R})$ and the convergence $\Psi^n(\mathbf{x}^n) \to \mathbf{y}$ holds in $\mathbb{D}_c(\mathbb{R}_+, \mathbb{R})$. An application of the continuous mapping principle implies convergence in distribution in $\mathbb{D}_c(\mathbb{R}_+, \mathbb{R})$ of the $X^n$ to $X$. By Lemma B.2, $\Psi$ is a continuous mapping from $\mathbb{R} \times \mathbb{D}_c(\mathbb{R}_+, \mathbb{R})^3$ to $\mathbb{D}_c(\mathbb{R}_+, \mathbb{R})$, so $X$ is a tight random element. Hence, it is a separable random element.

## 4. Proof of Theorem 3.1.

In Section 4.1, we derive certain properties of the processes that appear in the limit, emphasizing the property of being a separable random element. In Section 4.2, results on joint convergence in distribution are established. These developments culminate in the proof of Theorem 3.1 in Section 4.3. The hypotheses of Theorem 2.2 are assumed.

### 4.1. *Properties of processes associated with the Kiefer process.*

Let $K = ((K(t,x), x \in [0,1]), t \in \mathbb{R}_+)$ be a Kiefer process such that $K$, $Y$, $\hat{X}_0$ and $\tilde{X}$ are independent. Recall that this means that $K$ is a zero-mean Gaussian process with $\mathbf{E}K(t,x)K(s,y)=(t\wedge s)(x\wedge y-xy)$. In particular, $K(t,x)$ is a Brownian motion in $t$ for fixed $x$. We choose a modification of $K$ which is continuous in both variables a.s. Introduce for $t\in\mathbb{R}_+$ and $x\in\mathbb{R}_+$,

$$V(t,x)=K(\mathbf{a}(t),F(x)), \tag{4.1}$$

$$U(t,x)=K(t,F(x)), \tag{4.2}$$

and define

$$L(t,x)=V(t,x)+\int_0^x \frac{V(t,u-)}{1-F(u-)}\,dF(u), \tag{4.3}$$

$$L'(t,x)=U(t,x)+\int_0^x \frac{U(t,u-)}{1-F(u-)}\,dF(u), \tag{4.4}$$

$$G(t)=\int_0^t \frac{V(t-u,u-)}{1-F(u-)}\,dF(u). \tag{4.5}$$

The integrals on the right of (4.3)–(4.5) converge absolutely a.s. To see this, note that $V(t,x)$ is a locally square integrable martingale in $t$ for fixed $x$ relative to the natural filtration with predictable quadratic variation process $\mathbf{a}(t)F(x)(1-F(x))$. Therefore, by Theorem 1.9.5 in Liptser and Shiryayev [16], for $T>0$, $\mathbf{E}\sup_{t\in[0,T]}|V(t,x-)|\le 3\sqrt{\mathbf{a}(T)(1-F(x-))}$. We thus obtain

$$\begin{aligned}\int_0^\infty \frac{\mathbf{E}\sup_{t\in[0,T]}|V(t,u-)|}{1-F(u-)}\,dF(u)&\le\int_0^\infty \frac{3\sqrt{\mathbf{a}(t)}}{\sqrt{1-F(u-)}}\,dF(u)\\&\le 6\sqrt{\mathbf{a}(t)}\end{aligned} \tag{4.6}$$

so that

$$\int_0^\infty \frac{\sup_{t\in[0,T]}|V(t,u-)|}{1-F(u-)}\,dF(u)<\infty \qquad \text{a.s. for } T>0. \tag{4.7}$$

A similar argument shows that for arbitrary $T>0$ and $\delta>0$

$$\lim_{\varepsilon\to 0}\mathbf{P}\biggl(\sup_{t\in[0,T]}\sup_{x\in\mathbb{R}_+}\biggl|\int_0^x \frac{V(t,u-)}{1-F(u-)}\mathbf{1}_{\{F(u-)>1-\varepsilon\}}\,dF(u)\biggr|>\delta\biggr)=0 \tag{4.8}$$

and

$$\lim_{\varepsilon\to 0}\mathbf{P}\biggl(\sup_{t\in[0,T]}\sup_{x\in\mathbb{R}_+}\biggl|\int_0^x \frac{U(t,u-)}{1-F(u-)}\mathbf{1}_{\{F(u-)>1-\varepsilon\}}\,dF(u)\biggr|>\delta\biggr)=0. \tag{4.9}$$

Denote

$$F'(x)=\int_0^x \frac{1-F(u)}{1-F(u-)}\,dF(u). \tag{4.10}$$

Let $U=((U(t,x),x\in\mathbb{R}_+),t\in\mathbb{R}_+)$, $V=((V(t,x),x\in\mathbb{R}_+),t\in\mathbb{R}_+)$, $L=((L(t,x),x\in\mathbb{R}_+),t\in\mathbb{R}_+)$, $L'=((L'(t,x),x\in\mathbb{R}_+),t\in\mathbb{R}_+)$, and $G=(G(t),t\in\mathbb{R}_+)$.

LEMMA 4.1. *The process $L$ is a zero-mean Gaussian process with trajectories in $\mathbb{D}(\mathbb{R}_+, \mathbb{D}_c(\mathbb{R}_+, \mathbb{R}))$. Its covariance is given by*

$$\mathbf{E}L(t,x)L(s,y) = \mathbf{a}(t \wedge s)F'(x \wedge y).$$

*Furthermore, the pair $(L, V)$ is Gaussian.*

PROOF. By (4.1), (4.7) and Lebesgue's dominated convergence theorem, the definition of $L(t,x)$ in (4.3) implies that $L(t,x)$ is right-continuous in $x$ with left-hand limits for a given $t$. Similarly,

$$\lim_{t \downarrow s} \int_0^\infty \frac{|V(t,u-) - V(s,u-)|}{1 - F(u-)} \, dF(u) = 0,$$
$$\lim_{t \uparrow s} \int_0^\infty \frac{|V(t,u-) - V(s-,u-)|}{1 - F(u-)} \, dF(u) = 0.$$

It follows that $L$ has trajectories in $\mathbb{D}(\mathbb{R}_+, \mathbb{D}_c(\mathbb{R}_+, \mathbb{R}))$ a.s.

Since the integrand in the integral on the right of (4.3) is left-continuous in $u$, Lebesgue's dominated convergence theorem shows that this integral can be assumed to be a Stjeltjes integral, that is, to be a limit of Riemann sums. Since finite-dimensional distributions of $V$ are Gaussian, it follows that finite-dimensional distributions of $(L, V)$ are Gaussian too. The formula for the covariance of $L$ is obtained by a direct calculation. $\square$

Lemma 4.1 implies that $L$ defines an orthogonal random measure on $\mathbb{R}_+^2$ in that $\mathbf{E}\square L((t_1,x_1),(t_2,x_2))\square L((s_1,y_1),(s_2,y_2)) = 0$, where $t_1 \le t_2$, $x_1 \le x_2$, $s_1 \le s_2$, and $y_1 \le y_2$, whenever the rectangles with the vertices $(t_1,x_1),(t_1,x_2),(t_2,x_1),(t_2,x_2)$ and $(s_1,y_1),(s_1,y_2),(s_2,y_1),(s_2,y_2)$ are disjoint. Similarly, $\mathbf{E}\square L((t_1,x_1),(t_2,x_2))(L(0,y_2) - L(0,y_1)) = 0$, $\mathbf{E}\square L((t_1,x_1),(t_2,x_2))(L(s_2,0) - L(s_1,0)) = 0$, $\mathbf{E}(L(0,y_2) - L(0,y_1))(L(s_2,0) - L(s_1,0)) = 0$, $\mathbf{E}(L(t_2,0) - L(t_1,0))(L(s_2,0) - L(s_1,0)) = 0$ when $(t_1,t_2)$ and $(s_1,s_2)$ are disjoint, $\mathbf{E}(L(0,x_2) - L(0,x_1))(L(0,y_2) - L(0,y_1)) = 0$ when $(x_1,x_2)$ and $(y_1,y_2)$ are disjoint, $\mathbf{E}\square L((t_1,x_1),(t_2,x_2))L(0,0) = 0$, $\mathbf{E}(L(t_2,0) - L(t_1,0)) L(0,0) = 0$, and $\mathbf{E}(L(0,x_2) - L(0,x_1))L(0,0) = 0$. These properties enable us to define in a standard fashion integrals with respect to $L$.

To recapitulate, suppose $h$ is a Borel function with

$$\int_{\mathbb{R}_+^2} h(s,x)^2 \, d\mathbf{a}(s) \, dF'(x) < \infty. \tag{4.11}$$

If

$$\begin{aligned} h(s,x) = {} & \sum_{i=1}^k a_i \mathbf{1}_{\{s \in (s_{1,i}, s_{2,i}]\}} \mathbf{1}_{\{x \in (x_{1,i}, x_{2,i}]\}} + \sum_{i=1}^l b_i \mathbf{1}_{\{s=0\}} \mathbf{1}_{\{x \in (y_{1,i}, y_{2,i}]\}} \\ & + \sum_{i=1}^m c_i \mathbf{1}_{\{s \in (z_{1,i}, z_{2,i}]\}} \mathbf{1}_{\{x=0\}} + d\mathbf{1}_{\{s=0\}} \mathbf{1}_{\{x=0\}}, \end{aligned} \tag{4.12}$$

where the $a_i$, $b_i$, $c_i$ and $d$ are real numbers, $0 \le x_{1,i} < x_{2,i}$, $0 \le y_{1,i} < y_{2,i}$, $0 \le z_{1,i} < z_{2,i}$, and the sets $(s_{1,i}, s_{2,i}] \times (x_{1,i}, x_{2,i}]$, $(y_{1,i}, y_{2,i}]$, and $(z_{1,i}, z_{2,i}]$ are pairwise disjoint, we set

$$\begin{aligned}\int_{\mathbb{R}_+^2} h(s,x)\, dL(s,x) &= \sum_{i=1}^{k} a_i \square L((s_{1,i}, x_{1,i}), (s_{2,i}, x_{2,i})) \\ &\quad + \sum_{i=1}^{l} b_i (L(0, y_{2,i}) - L(0, y_{1,i})) \\ &\quad + \sum_{i=1}^{m} c_i (L(z_{2,i}, 0) - L(z_{1,i}, 0)) + dL(0,0).\end{aligned} \tag{4.13}$$

Note that by Lemma 4.1

$$\mathbf{E}\Big(\int_{\mathbb{R}_+^2} h(s,x)\, dL(s,x)\Big)^2 = \int_{\mathbb{R}_+^2} h(s,x)^2\, d\mathbf{a}(s)\, dF'(x). \tag{4.14}$$

If $h(s,x)$ is an arbitrary Borel function satisfying (4.11), then there exists a sequence $h^k$ of functions of the form (4.12) such that

$$\lim_{k\to\infty} \int_{\mathbb{R}_+^2} (h(s,x) - h^k(s,x))^2\, d\mathbf{a}(s)\, dF'(x) = 0.$$

This implies by (4.14) that the sequence $\int_{\mathbb{R}_+^2} h^k(s,x)\, dL(s,x)$ is Cauchy in $L_2(\Omega, \mathcal{F}, \mathbf{P})$, so it converges. We define $\int_{\mathbb{R}_+^2} h(s,x)\, dL(s,x)$ as the limit. One can see that the integral is a zero-mean Gaussian random variable and that (4.14) still holds. By polarization, if $g(s,x)$ is another function with $\int_{\mathbb{R}_+^2} g(s,x)^2\, d\mathbf{a}(s)\, dF'(x) < \infty$, then

$$\begin{aligned}\mathbf{E} \int_{\mathbb{R}_+^2} h(s,x)\, dL(s,x) \int_{\mathbb{R}_+^2} g(s,x)\, dL(s,x) \\ = \int_{\mathbb{R}_+^2} h(s,x) g(s,x)\, d\mathbf{a}(s)\, dF'(x).\end{aligned} \tag{4.15}$$

The martingale property asserted in the next lemma is understood with respect to the natural filtration. Recall that Gaussian martingales are locally square integrable.

Lemma 4.2. *If $g(s,x)$ is a Borel function on $\mathbb{R}_+^2$ such that, for all $t \in \mathbb{R}_+$ and $x \in \mathbb{R}_+$,*

$$\int_{\mathbb{R}_+^2} \mathbf{1}_{\{s \le t\}} \mathbf{1}_{\{y \le x\}} g(s,y)^2\, d\mathbf{a}(s)\, dF'(y) < \infty,$$

*then the process* $N = (N(t), t \in \mathbb{R}_+)$ *with*

$$N(t) = \int_{\mathbb{R}_+^2} g(s,x) \mathbf{1}_{\{s+x \le t\}} \, dL(s,x)$$

*is a zero-mean Gaussian martingale with predictable quadratic variation process* $\langle N \rangle = (\langle N \rangle(t), t \in \mathbb{R}_+)$, *where*

$$\langle N \rangle(t) = \int_{\mathbb{R}_+^2} g(s,x)^2 \mathbf{1}_{\{s+x \le t\}} \, d\mathbf{a}(s) \, dF'(x).$$

*The process* $N$ *has a modification with trajectories in* $\mathbb{D}(\mathbb{R}_+, \mathbb{R})$.

PROOF. By construction, finite-dimensional distributions of $N$ are limits of Gaussian distributions, so they are Gaussian too. The covariance of $N$ is given, according to (4.15), by

$$\mathbf{E} N(s) N(t) = \int_{\mathbb{R}_+^2} g(s,x)^2 \mathbf{1}_{\{u+x \le s \wedge t\}} \, d\mathbf{a}(u) \, dF'(x).$$

It follows that $N$ has uncorrelated, hence, independent increments. It is, thus, a zero-mean Gaussian martingale. By the independence-of-increments property, the natural filtration of $N$ is right-continuous; see, for example, Doob [6], Part 2, Chapter VI, Section 8. Consequently, $N$ admits a right-continuous with left-hand limits modification; see, for example, Doob [6], Part 2, Chapter IV, Section 1. $\square$

Let

$$M(t) = \int_{\mathbb{R}_+^2} \mathbf{1}_{\{s+x \le t\}} \, dL(s,x). \tag{4.16}$$

By Lemma 4.2, the process $M = (M(t), t \in \mathbb{R}_+)$ admits a modification which is a Gaussian martingale with trajectories from $\mathbb{D}(\mathbb{R}_+, \mathbb{R})$. We further consider such a modification throughout. The variance of $M(t)$ is given by

$$C(t) = \int_0^t \mathbf{a}(t-s) \, dF'(s). \tag{4.17}$$

We study measurability properties of the introduced processes. Recall that $H$ was defined in (3.17).

LEMMA 4.3. *The processes* $H, G$ *and* $M$ *are separable random elements of* $\mathbb{D}_c(\mathbb{R}_+, \mathbb{R})$ *and the processes* $L, L', U$ *and* $V$ *are separable random elements of* $\mathbb{D}_c(\mathbb{R}_+, \mathbb{D}_c(\mathbb{R}_+, \mathbb{R}))$. *The pair* $(G, M)$ *is Gaussian. The process* $G$ *is a Gaussian semimartingale.*

PROOF. By the hypotheses of Theorem 2.2, the process $Y$ is a separable random element of $\mathbb{D}_c(\mathbb{R}_+,\mathbb{R})$, so its distribution is a tight probability measure. By (3.17), $H$ is obtained from $Y$ by an application of a continuous operator on $\mathbb{D}_c(\mathbb{R}_+,\mathbb{R})$. It follows that the distribution of $H$ is a tight probability measure, hence, $H$ is a separable random element of $\mathbb{D}_c(\mathbb{R}_+,\mathbb{R})$.

Since the process $V$ is bounded on bounded domains [see (4.1)], by (4.5) the process $G$ is of locally bounded variation a.s. By (4.7) and Lebesgue's dominated convergence theorem, the trajectories of $G$ are right-continuous and admit left-hand limits on a set of full probability. Since the integrand in (4.5) is a left-continuous function of $y$, the integral can be interpreted as a Stjeltjes integral, that is, as a limit of Riemann sums. As $(L,V)$ is a Gaussian pair by Lemma 4.1, the definition of $M$ in (4.16) implies that $(G,M)$ is a Gaussian pair.

Hence, the processes $G$ and $M$ are Gaussian semimartingales with paths in $\mathbb{D}(\mathbb{R}_+,\mathbb{R})$. By Liptser and Shiryaev [16], Theorem 4.9.1, their jump times are deterministic, so the ranges of these processes as elements of $\mathbb{D}_c(\mathbb{R}_+,\mathbb{R})$ are separable. Since balls in $\mathbb{D}_c(\mathbb{R}_+,\mathbb{R})$ belong to the Kolmogorov $\sigma$-algebra, the traces of the Kolmogorov and Borel $\sigma$-algebras on a separable set coincide, hence, $G$ and $M$ are separable random elements of $\mathbb{D}_c(\mathbb{R}_+,\mathbb{R})$.

By (4.2) and continuity of the Kiefer process in both variables, $U(t,x)$, as a function of $x$, jumps only when $F$ jumps, so its range is separable, hence, $(U(t,x),x\in\mathbb{R}_+)$ is a separable random element of $\mathbb{D}_c(\mathbb{R}_+,\mathbb{R})$ for each $t$. Next, the map $t\to(U(t,x),x\in\mathbb{R}_+)$ from $\mathbb{R}_+$ to $\mathbb{D}_c(\mathbb{R}_+,\mathbb{R})$ is continuous, so $U$ has a separable range in $\mathbb{D}_c(\mathbb{R}_+,\mathbb{D}_c(\mathbb{R}_+,\mathbb{R}))$. It is, therefore, a separable random element of the latter space. The process $V$ is a separable random element of $\mathbb{D}_c(\mathbb{R}_+,\mathbb{D}_c(\mathbb{R}_+,\mathbb{R}))$ for a similar reason.

The range of $(L(t,x),x\in\mathbb{R}_+)$ for a given $t$ is separable by the fact that the jumps of $(L(t,x),x\in\mathbb{R}_+)$ occur at the times of jumps of $F$, hence, this process is a random element of $\mathbb{D}_c(\mathbb{R}_+,\mathbb{R})$. Since the jumps of $L$ as an element of $\mathbb{D}_c(\mathbb{R}_+,\mathbb{D}_c(\mathbb{R}_+,\mathbb{R}))$ coincide with the jumps of $\mathbf{a}$, it is a separable random element of $\mathbb{D}_c(\mathbb{R}_+,\mathbb{D}_c(\mathbb{R}_+,\mathbb{R}))$.

The assertion of the lemma for $L'$ is obtained analogously. □

In what follows, we always assume the modifications as described in Lemmas 4.2 and 4.3.

We now construct the process $Z=(Z(t),t\in\mathbb{R}_+)$ in the statement of Theorem 2.2. We define

$$Z(t)=G(t)-M(t). \tag{4.18}$$

By Lemma 4.3, $Z$ is a Gaussian semimartingale, so it is a separable random element of $\mathbb{D}_c(\mathbb{R}_+,\mathbb{R})$. In order to verify that its covariance function has the form stated in Theorem 2.2, we find it convenient to approximate this process with Gaussian processes of simpler structure.

For $l \in \mathbb{N}$, let $0 = s_0^l < s_1^l < s_2^l < \cdots$ be a strictly increasing to infinity sequence of real numbers. Following the notation used in the Introduction, we set, for $t \in \mathbb{R}_+$,

$$I_{l,t}(s,x) = \sum_{i=1}^{\infty} \mathbf{1}_{\{s \in (s_{i-1}^l, s_i^l]\}} \mathbf{1}_{\{0 \le x \le t - s_{i-1}^l\}} + \mathbf{1}_{\{s=0\}} \mathbf{1}_{\{0 \le x \le t\}}. \tag{4.19}$$

We also define

$$\begin{aligned} &\int_{\mathbb{R}_+^2} I_{l,t}(s,x)\, dV(s,x) \\ &\quad = \sum_{i=1}^{\infty} (V(s_i^l, t - s_{i-1}^l) - V(s_{i-1}^l, t - s_{i-1}^l)) \mathbf{1}_{\{s_{i-1}^l \le t\}} + V(0,t) \end{aligned} \tag{4.20}$$

and introduce

$$Z_l(t) = -\int_{\mathbb{R}_+^2} I_{l,t}(s,x)\, dV(s,x). \tag{4.21}$$

On recalling (4.1), we see that $(Z_l(t), t \in \mathbb{R}_+)$ is a zero-mean Gaussian process with covariance

$$\begin{aligned} \mathbf{E} Z_l(t) Z_l(s) = {} & \sum_{i=1}^{\infty} (\mathbf{a}(s_i^l) - \mathbf{a}(s_{i-1}^l)) F(t \wedge s - s_{i-1}^l) \\ & \times (1 - F(t \vee s - s_{i-1}^l)) \mathbf{1}_{\{s_{i-1}^l \le t \wedge s\}} \\ & + \mathbf{a}(0) F(t \wedge s)(1 - F(t \vee s)). \end{aligned} \tag{4.22}$$

LEMMA 4.4. *If $\sup_i (s_i^l - s_{i-1}^l) \to 0$ as $l \to \infty$, then the $Z_l(t)$ converge to $Z(t)$ in the mean square sense for each $t$ as $l \to \infty$. The covariance of $Z$ is given by*

$$\mathbf{E} Z(t) Z(s) = \int_0^{t \wedge s} F(t \wedge s - u)(1 - F(t \vee s - u))\, d\mathbf{a}(u).$$

PROOF. By (4.3), (4.13), (4.19) and (4.20),

$$\begin{aligned} \int_{\mathbb{R}_+^2} I_{l,t}(s,x)\, dV(s,x) = {} & \int_{\mathbb{R}_+^2} I_{l,t}(s,x)\, dL(s,x) \\ & - \sum_{i=1}^{\infty} \mathbf{1}_{\{s_{i-1}^l \le t\}} \int_0^{t - s_{i-1}^l} \frac{V(s_i^l, y-) - V(s_{i-1}^l, y-)}{1 - F(y-)}\, dF(y) \\ & - \int_0^t \frac{V(0, y-)}{1 - F(y-)}\, dF(y). \end{aligned} \tag{4.23}$$

Given $\varepsilon > 0$, for all $l$ large enough, $\mathbf{1}_{\{s+x\le t\}}\mathbf{1}_{\{s\ge 0\}}\mathbf{1}_{\{x\ge 0\}} \le I_{l,t}(s,x) \le \mathbf{1}_{\{s+x\le t+\varepsilon\}}\mathbf{1}_{\{s\ge 0\}}\mathbf{1}_{\{x\ge 0\}}$. Therefore,

$$\int_{\mathbb{R}_+^2} (I_{l,t}(s,x) - \mathbf{1}_{\{s+x\le t\}})^2 \, d\mathbf{a}(s)\, dF'(x) \\ \le \int_{\mathbb{R}_+^2} (\mathbf{1}_{\{s+x\le t+\varepsilon\}} - \mathbf{1}_{\{s+x\le t\}})^2 \, d\mathbf{a}(s)\, dF'(x).$$

The right-hand side converges to zero as $\varepsilon \to 0$, so by (4.16) and the definition of the integral

$$M(t) = \operatorname*{l.i.m.}_{l\to\infty} \int_{\mathbb{R}_+^2} I_{l,t}(s,x)\, dL(s,x), \tag{4.24}$$

where l.i.m. stands for mean-square limit. Also

$$\sum_{i=1}^{\infty} \mathbf{1}_{\{s_{i-1}^l \le t\}} \int_0^{t-s_{i-1}^l} \frac{V(s_i^l, y-) - V(s_{i-1}^l, y-)}{1-F(y-)}\, dF(y) + \int_0^t \frac{V(0,y-)}{1-F(y-)}\, dF(y) \\ = \int_0^t \frac{V(s_{i(y)}^l, y-)}{1-F(y-)}\, dF(y),$$

where $s_{i(y)}^l > t - y \ge s_{i(y)-1}^l$. By right continuity of $V(t,x)$ in $t$, (4.5), (4.8) and Lebesgue's dominated convergence theorem,

$$G(t) = \mathbf{P}\text{-}\lim_{l\to\infty} \Bigg( \sum_{i=1}^{\infty} \mathbf{1}_{\{s_{i-1}^l \le t\}} \int_0^{t-s_{i-1}^l} \frac{V(s_i^l, y-) - V(s_{i-1}^l, y-)}{1-F(y-)}\, dF(y) \\ + \int_0^t \frac{V(0,y-)}{1-F(y-)} dF(y) \Bigg), \tag{4.25}$$

where $\mathbf{P}$- lim denotes limit in probability.

By (4.18), (4.21), (4.23)–(4.25), $Z(t) = \mathbf{P}\text{-}\lim_{l\to\infty} Z_l(t)$. Since the $Z_l$ are Gaussian processes, the latter limit holds in the mean square sense too; see, for example, Ibragimov and Rozanov [11], Lemma I.3.1.

By (4.22), we can write

$$\mathbf{E}Z_l(t)Z_l(s) = \int_0^{\infty} F(t\wedge s - r^l(u))(1 - F(t\vee s - r^l(u)))\mathbf{1}_{\{r^l(u)\le t\wedge s\}}\, d\mathbf{a}(u),$$

where

$$r^l(u) = \sum_{i=1}^{\infty} s_{i-1}^l \mathbf{1}_{\{u\in(s_{i-1}^l, s_i^l]\}}.$$

Since $\max_i(s_i^l - s_{i-1}^l) \to 0$ as $l\to\infty$, $r^l(u) \to u$ from the left. By right continuity of $F$, $F(t\wedge s - r^l(u)) \to F(t\wedge s - u)$ and $F(t\vee s - r^l(u)) \to F(t\vee s - u)$ as $l \to \infty$. Also, $\mathbf{1}_{\{r^l(u) \le t\wedge s\}} \to \mathbf{1}_{\{u\le t\wedge s\}}$. Therefore, by Lebesgue's dominated convergence theorem,

$$\lim_{l\to\infty} \mathbf{E} Z_l(t) Z_l(s) = \int_0^{t\wedge s} F(t\wedge s - u)(1 - F(t\vee s - u))\, d\mathbf{a}(u). \qquad \square$$

4.2. *Convergence in distribution.* Let

$$U^n(t,x) = \frac{1}{\sqrt{n}} \sum_{i=1}^{\lfloor nt \rfloor} (\mathbf{1}_{\{\eta_i \le x\}} - F(x)), \tag{4.26}$$

$$L'^n(t,x) = \frac{1}{\sqrt{n}} \sum_{i=1}^{\lfloor nt \rfloor} \left( \mathbf{1}_{\{\eta_i \le x\}} - \int_0^{\eta_i \wedge x} \frac{dF(u)}{1 - F(u-)} \right). \tag{4.27}$$

By (3.1) and (4.26),

$$V^n(t,x) = U^n\left( \frac{A^n(t)}{n}, x \right). \tag{4.28}$$

For what follows, we note that, given arbitrary $T > 0$ and $\delta > 0$,

$$\begin{aligned} &\lim_{\varepsilon\to 0} \limsup_{n\to\infty} \mathbf{P}\left( \sup_{t\in[0,T]} \sup_{x\in\mathbb{R}_+} \left| \int_0^x \frac{U^n(t,u-)}{1-F(u-)} \mathbf{1}_{\{F(u-) > 1-\varepsilon\}}\, dF(u) \right| > \delta \right) \\ &\quad = 0. \end{aligned} \tag{4.29}$$

The limit in (4.29) is analogous to that in equation (3.23) in Krichagina and Puhalskii [15] [see also (3.24) in that paper], so a similar proof applies.

We introduce the processes $L'^n = ((L'^n(t,x), x \in \mathbb{R}_+), t \in \mathbb{R}_+)$, $U^n = ((U^n(t, x), x \in \mathbb{R}_+), t \in \mathbb{R}_+)$, and $L' = ((L'(t,x), x \in \mathbb{R}_+), t \in \mathbb{R}_+)$.

LEMMA 4.5. *As $n \to \infty$, the $(\hat{X}_0^n, \tilde{X}^n, U^n, L'^n)$ converge jointly in distribution in $\mathbb{R} \times \mathbb{D}_c(\mathbb{R}_+,\mathbb{R}) \times \mathbb{D}_c(\mathbb{R}_+, \mathbb{D}_c(\mathbb{R}_+,\mathbb{R}))^2$ to $(\hat{X}_0, \tilde{X}, U, L')$.*

PROOF. Recall that, by Lemma 4.3, $L'$ and $U$ are separable random elements of the space $\mathbb{D}_c(\mathbb{R}_+, \mathbb{D}_c(\mathbb{R}_+,\mathbb{R}))$. The process $\tilde{X}$ is a separable random element of $\mathbb{D}_c(\mathbb{R}_+,\mathbb{R})$ by (3.20) and the assumption that $S$ is a separable random element. The hypotheses of Theorem 2.2 and (2.2c) imply in a standard fashion that the $(\hat{X}_0^n, \tilde{X}^n)$ converge in distribution in $\mathbb{R} \times \mathbb{D}_c(\mathbb{R}_+,\mathbb{R})$ to $(\hat{X}_0, \tilde{X})$.

Since $(\hat{X}_0^n, \tilde{X}^n)$ and $(U^n, L'^n)$ are independent and $(\hat{X}_0, \tilde{X})$ and $(U, L')$ are also independent and are separable random elements of the associated metric spaces, by Theorem D.8 it suffices to establish convergence in distribution in $\mathbb{D}_c(\mathbb{R}_+, \mathbb{D}_c(\mathbb{R}_+,\mathbb{R}))^2$ for $(U^n, L'^n)$. By (4.26) and (4.27) (cf. (3.5)),

$$U^n(t,x) = -\int_0^x \frac{U^n(t,u-)}{1-F(u-)}\, dF(u) + L'^n(t,x). \tag{4.30}$$

Let

$$K^n(t,x) = \frac{1}{\sqrt{n}} \sum_{i=1}^{\lfloor nt \rfloor} (\mathbf{1}_{\{\zeta_i \le x\}} - x), \tag{4.31}$$

where the $\zeta_i$ are independent and uniform on $[0,1]$. By Krichagina and Puhalskii [15], the $K^n = ((K^n(t,x), x \in [0,1]), t \in \mathbb{R}_+)$ converge to $K$ in distribution in $\mathbb{D}(\mathbb{R}_+, \mathbb{D}([0,1],\mathbb{R}))$. Since $K$ is continuous in both variables, it follows that the convergence takes place in $\mathbb{D}_c(\mathbb{R}_+, \mathbb{D}_c([0,1],\mathbb{R}))$ too. (One can apply Corollary D.1.) Since by (4.26), we can assume that $U^n(t,x) = K^n(t, F(x))$, the $U^n$ converge in distribution in $\mathbb{D}_c(\mathbb{R}_+, \mathbb{D}_c(\mathbb{R}_+,\mathbb{R}))$ to the process $U$.

Let, for $\varepsilon \in (0,1)$,

$$L'^{n,\varepsilon}(t,x) = U^n(t,x) + \int_0^x \frac{U^n(t,u-)}{1-F(u-)} \mathbf{1}_{\{F(u-) \le 1-\varepsilon\}}\, dF(u),$$

$$L'^{\varepsilon}(t,x) = U(t,x) + \int_0^x \frac{U(t,u-)}{1-F(u-)} \mathbf{1}_{\{F(u-) \le 1-\varepsilon\}}\, dF(u),$$

$L'^{n,\varepsilon} = ((L'^{n,\varepsilon}(t,x), x \in \mathbb{R}_+), t \in \mathbb{R}_+)$, and $L'^{\varepsilon} = ((L'^{\varepsilon}(t,x), x \in \mathbb{R}_+), t \in \mathbb{R}_+)$. An argument analogous to the one used in the proof of Lemma 4.3 shows that $L'^{\varepsilon}$ is a separable random element of $\mathbb{D}_c(\mathbb{R}_+, \mathbb{D}_c(\mathbb{R}_+,\mathbb{R}))$. The continuous mapping principle (see Theorem D.1) yields the convergence in distribution in $\mathbb{D}_c(\mathbb{R}_+, \mathbb{D}_c(\mathbb{R}_+,\mathbb{R}))^2$ of the $(U^n, L'^{n,\varepsilon})$ to $(U, L'^{\varepsilon})$. Thus, in view of (4.9), (4.30) and Theorem D.10, the result follows by (4.29). □

Lemma 4.6. *As $n \to \infty$, the $(\hat{X}_0^n, \tilde{X}^n, H^n, G^n, L^n, V^n)$ converge in distribution in $\mathbb{R} \times \mathbb{D}_c(\mathbb{R}_+,\mathbb{R})^3 \times \mathbb{D}_c(\mathbb{R}_+, \mathbb{D}_c(\mathbb{R}_+,\mathbb{R}))^2$ to $(\hat{X}, \tilde{X}, H, G, L, V)$.*

Proof. By Lemma 4.3, the processes $H$, $G$, $L$, $V$ and $U$ are separable random elements of the associated function spaces. Since the exogenous arrival process $E^n$ and $(\hat{X}_0^n, \tilde{X}^n, U^n, L'^n)$ are independent, on the one hand, and $Y$ and $(\hat{X}_0, \tilde{X}, U, L')$ are independent, on the other hand, and are separable random elements, Lemma 4.5, the hypotheses of Theorem 2.2, and Theorem D.8 imply that the $(\hat{X}_0^n, \tilde{X}^n, Y^n, U^n, L'^n)$ converge in distribution in $\mathbb{R} \times \mathbb{D}_c(\mathbb{R}_+,\mathbb{R})^2 \times \mathbb{D}_c(\mathbb{R}_+, \mathbb{D}_c(\mathbb{R}_+,\mathbb{R}))^2$ to $(\hat{X}_0, \tilde{X}, Y, U, L')$. Since the random element $(\hat{X}_0, \tilde{X}, Y, U, L')$ is separable, by Theorem 2.1 and Slutsky's lemma (Theorem D.9) the $(\hat{X}_0^n, \tilde{X}^n, Y^n, U^n, L'^n, A^n/n)$ jointly converge in distribution in $\mathbb{R} \times \mathbb{D}_c(\mathbb{R}_+,\mathbb{R})^2 \times \mathbb{D}_c(\mathbb{R}_+, \mathbb{D}_c(\mathbb{R}_+,\mathbb{R}))^2 \times \mathbb{D}_c(\mathbb{R}_+,\mathbb{R})$ to $(\hat{X}_0, \tilde{X}, Y, U, L', \mathbf{a})$. On recalling that $V^n(t,x) = U^n(A^n(t)/n, x)$, $V(t,x) = U(\mathbf{a}(t), x)$, $L^n(t,x) = L'^n(A^n(t)/n, x)$, and $L(t,x) = L'(\mathbf{a}(t), x)$ [see (3.4), (4.1), (4.4), (4.27) and (4.28)], we conclude by the continuous mapping principle (Theorem D.1) that the $(\hat{X}_0^n, \tilde{X}^n, Y^n, V^n, L^n)$ converge in distribution in $\mathbb{R} \times \mathbb{D}_c(\mathbb{R}_+,\mathbb{R})^2 \times \mathbb{D}_c(\mathbb{R}_+, \mathbb{D}_c(\mathbb{R}_+,\mathbb{R}))^2$ to $(\hat{X}_0, \tilde{X}, Y, V, L)$.

Let, for $\varepsilon \in (0,1)$,

$$G^{n,\varepsilon}(t) = \int_0^t \frac{V^n(t-u,u-)}{1-F(u-)} \mathbf{1}_{\{F(u-)\le 1-\varepsilon\}}\, dF(u),$$

$$G^{\varepsilon}(t) = \int_0^t \frac{V(t-u,u-)}{1-F(u-)} \mathbf{1}_{\{F(u-)\le 1-\varepsilon\}}\, dF(u),$$

$G^{n,\varepsilon} = (G^{n,\varepsilon}(t), t \in \mathbb{R}_+)$, and $G^{\varepsilon} = (G^{\varepsilon}(t), t \in \mathbb{R}_+)$. We note that the $G^{\varepsilon}$ are separable random elements of $\mathbb{D}_c(\mathbb{R}_+,\mathbb{R})$. On applying the continuous mapping principle and recalling (3.18), we have that the $(\hat{X}_0^n, \tilde{X}^n, H^n, G^{n,\varepsilon}, L^n, V^n)$ converge in distribution in $\mathbb{R}\times\mathbb{D}_c(\mathbb{R}_+,\mathbb{R})^3\times\mathbb{D}_c(\mathbb{R}_+,\mathbb{D}_c(\mathbb{R}_+,\mathbb{R}))^2$ to $(\hat{X}_0, \tilde{X}, H, G^{\varepsilon}, L, V)$. In view of Theorem D.10, (3.7), (4.5) and (4.8), the assertion of the lemma will follow if for arbitrary $T>0$ and $\delta>0$,

$$\lim_{\varepsilon\to 0}\limsup_{n\to\infty}\mathbf{P}\left(\sup_{t\in[0,T]}\sup_{x\in\mathbb{R}_+}\left|\int_0^x \frac{V^n(t,u-)}{1-F(u-)}\mathbf{1}_{\{F(u-)>1-\varepsilon\}}\, dF(u)\right| > \delta\right) = 0.$$

The latter limit is implied by (4.28) and (4.29). $\square$

Let, for $\varepsilon>0$,

$$M^{n,\varepsilon}(t) = \int_{\mathbb{R}_+^2} \mathbf{1}_{\{\Delta\mathbf{a}(s)\Delta F(x)>\varepsilon\}}\mathbf{1}_{\{s+x\le t\}}\, dL^n(s,x) \tag{4.32}$$

and

$$M^{\varepsilon}(t) = \int_{\mathbb{R}_+^2} \mathbf{1}_{\{\Delta\mathbf{a}(s)\Delta F(x)>\varepsilon\}}\mathbf{1}_{\{s+x\le t\}}\, dL(s,x). \tag{4.33}$$

We note that these two integrals are, in fact, finite sums:

$$M^{n,\varepsilon}(t) = \sum_{\substack{s,x: s+x\le t,\\ \Delta\mathbf{a}(s)\Delta F(x)>\varepsilon}} \square L^n((s-,x-),(s,x)) \tag{4.34}$$

and

$$M^{\varepsilon}(t) = \sum_{\substack{s,x: s+x\le t,\\ \Delta\mathbf{a}(s)\Delta F(x)>\varepsilon}} \square L((s-,x-),(s,x)). \tag{4.35}$$

In particular, $M^{\varepsilon} = (M^{\varepsilon}(t), t\in\mathbb{R}_+)$ is a separable random element of $\mathbb{D}_c(\mathbb{R}_+,\mathbb{R})$ and $M^{n,\varepsilon} = (M^{n,\varepsilon}(t), t\in\mathbb{R}_+)$ is a continuous function of $L^n$ for spaces $\mathbb{D}_c(\mathbb{R}_+,\mathbb{R})$ and $\mathbb{D}_c(\mathbb{R}_+,\mathbb{D}_c(\mathbb{R}_+,\mathbb{R}))$.

LEMMA 4.7. *Given $t_1 < \cdots < t_m \in \mathbb{R}_+$, as $n\to\infty$, the $(\hat{X}_0^n, \tilde{X}^n, H^n, G^n, M^{n,\varepsilon}, M^n(t_1),\dots,M^n(t_m))$ converge in distribution in $\mathbb{R}\times\mathbb{D}_c(\mathbb{R}_+,\mathbb{R})^4\times\mathbb{R}^m$ to $(\hat{X}_0, \tilde{X}, H, G, M^{\varepsilon}, M(t_1),\dots,M(t_m))$.*

PROOF. Since $Z^n(t)=G^n(t)-M^n(t)$ by (3.6) and $Z(t)=G(t)-M(t)$ by (4.18), the continuous mapping principle implies that it suffices to prove that the $(\hat X_0^n,\tilde X^n,H^n,G^n,M^{n,\varepsilon},Z^n(t_1),\dots,Z^n(t_m))$ converge in distribution in $\mathbb{R}\times\mathbb{D}_c(\mathbb{R}_+,\mathbb{R})^4\times\mathbb{R}^m$ to $(\hat X_0,\tilde X,H,G,M^\varepsilon,Z(t_1),\dots,Z(t_m))$.

Let $0=s_0^l<s_1^l<s_2^l<\cdots$ be such that $s_k^l\to\infty$ as $k\to\infty$ and $\sup_{i=1,2,\dots}(s_i^l-s_{i-1}^l)\to 0$ as $l\to\infty$. Recall that the function $I_{l,t}$ and the process $Z_l=(Z_l(t),t\in\mathbb{R}_+)$ are defined by (4.19) and (4.21), respectively. We also set

$$Z_l^n(t)=-\int_{\mathbb{R}_+^2}I_{l,t}(s,x)\,dV^n(s,x), \tag{4.36}$$

where the integral on the right is defined in analogy with (4.20), that is,

$$\int_{\mathbb{R}_+^2}I_{l,t}(s,x)\,dV^n(s,x)=\sum_{i=1}^\infty \mathbf{1}_{\{s_{i-1}^l\le t\}}(V^n(s_i^l,t-s_{i-1}^l)-V^n(s_{i-1}^l,t-s_{i-1}^l))+V^n(0,t).$$

Lemma 4.6, (4.34), (4.35) and the continuous mapping principle yield convergence in distribution in $\mathbb{R}\times\mathbb{D}_c(\mathbb{R}_+,\mathbb{R})^4\times\mathbb{R}^m$ of the $(\hat X_0^n,\tilde X^n,H^n,G^n,M^{n,\varepsilon},Z_l^n(t_1),\dots,Z_l^n(t_m))$ to $(\hat X_0,\tilde X,H,G,M^\varepsilon,Z_l(t_1),\dots,Z_l(t_m))$ as $n\to\infty$. By Lemma 4.4, for arbitrary $t\in\mathbb{R}_+$ and $\gamma>0$,

$$\lim_{l\to\infty}\mathbf{P}(|Z(t)-Z_l(t)|>\gamma)=0.$$

It thus remains to prove that, for arbitrary $t\in\mathbb{R}_+$ and $\gamma>0$,

$$\lim_{l\to\infty}\limsup_{n\to\infty}\mathbf{P}(|Z^n(t)-Z_l^n(t)|>\gamma)=0. \tag{4.37}$$

By (3.1), (3.3) and (4.36), for $s\in\mathbb{R}_+$,

$$\begin{aligned}Z^n(s)&-Z_l^n(s)\\&=\frac{1}{\sqrt n}\sum_{i=1}^{A^n(s)}\sum_{j=1}^\infty \mathbf{1}_{\{s_{j-1}^l\le s\}}\mathbf{1}_{\{\tau_i^n\in(s_{j-1}^l,s_j^l]\}}\\&\qquad\times(\mathbf{1}_{\{\eta_i\in(s-\tau_i^n,s-s_{j-1}^l]\}}-(F(s-s_{j-1}^l)-F(s-\tau_i^n))).\end{aligned}$$

For $k\in\mathbb{N}$, introduce

$$\begin{aligned}Z_{l,k}^n(s)=\frac{1}{\sqrt n}\sum_{i=1}^{A^n(s)\wedge k}\sum_{j=1}^\infty \mathbf{1}_{\{s_{j-1}^l\le s\}}\mathbf{1}_{\{\tau_i^n\in(s_{j-1}^l,s_j^l]\}}\\\times(\mathbf{1}_{\{\eta_i\in(s-\tau_i^n,s-s_{j-1}^l]\}}-(F(s-s_{j-1}^l)-F(s-\tau_i^n))).\end{aligned}$$

Since the $A^n(s)/n$ converge to $\mathbf{a}(s)$ in probability as $n\to\infty$ by Theorem 2.1, (4.37) would follow by

$$\lim_{l\to\infty}\limsup_{n\to\infty}\mathbf{P}(|Z^n_{l,k}(t)|>\gamma)=0. \tag{4.38}$$

Let $\mathcal{F}^n_s$ be complete $\sigma$-algebras generated by the random variables $\tau^n_j\wedge\tau^n_{A^n(s)+1}$ and $\eta_{j\wedge A^n(s)}$, where $j\in\mathbb{N}$. By Brémaud [4], Appendix A3, Theorem 25, the flow $\mathbf{F}^n=(\mathcal{F}^n_s, s\in\mathbb{R}_+)$ is right-continuous, so it is a filtration. By part 4 of Lemma C.1 (see also Lemma 8 in Reed [22] and Lemma 5.2 in Krichagina and Puhalskii [15]), the process $(Z^n_{l,k}(s), s\in\mathbb{R}_+)$ is an $\mathbf{F}^n$-square integrable martingale with predictable quadratic variation process $(\langle Z^n_{l,k}\rangle(s), s\in\mathbb{R}_+)$, where

$$\begin{aligned}\langle Z^n_{l,k}\rangle(s)&=\frac{1}{n}\sum_{i=1}^{A^n(s)\wedge k}\sum_{j=1}^{\infty}\mathbf{1}_{\{s^l_{j-1}\le s\}}\mathbf{1}_{\{\tau^n_i\in(s^l_{j-1},s^l_j]\}}(F(s-s^l_{j-1})-F(s-\tau^n_i))\\&\qquad\times(1-(F(s-s^l_{j-1})-F(s-\tau^n_i)))\\&\le\frac{1}{n}\sum_{j=1}^{\infty}\mathbf{1}_{\{s^l_{j-1}\le s\}}\int_0^s\sum_{i=1}^{A^n(s)}\mathbf{1}_{\{\tau^n_i\in(s^l_{j-1},s^l_j]\}}\mathbf{1}_{\{x\in(s-\tau^n_i,s-s^l_{j-1}]\}}\,dF(x)\\&=\frac{1}{n}\sum_{j=1}^{\infty}\mathbf{1}_{\{s^l_{j-1}\le s\}}\int_{s-s^l_j}^{s-s^l_{j-1}}(A^n(s^l_j\wedge s)-A^n(s-x))\,dF(x).\end{aligned}$$

In view of the compact convergence in probability of the $A^n/n$ to $\mathbf{a}$ as $n\to\infty$ (Theorem 2.1), the latter sum converges in probability as $n\to\infty$ to

$$\begin{aligned}&\sum_{j=1}^{\infty}\mathbf{1}_{\{s^l_{j-1}\le s\}}\int_{s-s^l_j}^{s-s^l_{j-1}}(\mathbf{a}(s^l_j\wedge s)-\mathbf{a}(s-x))\,dF(x)\\&\quad=\sum_{j=1}^{\infty}\mathbf{1}_{\{s^l_{j-1}\le s\}}(F(s-s^l_{j-1})-F((s-s^l_j)^+))\mathbf{a}(s^l_j\wedge s)\\&\qquad-\int_0^s\mathbf{a}(s-x)\,dF(x)\\&\quad=\int_0^s\mathbf{a}(u^l(s-x))\,dF(x)-\int_0^s\mathbf{a}(s-x)\,dF(x),\end{aligned}$$

where $u^l(x)=\sum_{j=1}^{\infty}s^l_j\wedge s\mathbf{1}_{\{s^l_{j-1}\le s\}}\mathbf{1}_{\{x\in[s^l_{j-1},s^l_j)\}}$ for $x\in[0,s]$. Note that $u^l(x)\ge x$ and $u^l(x)\to x$ as $l\to\infty$. Hence, $u^l(s-x)\to s-x$ from the right, so $\mathbf{a}(u^l(s-x))\to\mathbf{a}(s-x)$. By Lebesgue's dominated convergence theorem, $\int_0^s \mathbf{a}(u^l(s-x))\,dF(x) \to \int_0^s \mathbf{a}(s-x)\,dF(x)$ as $l\to\infty$. We conclude that, for arbitrary $\delta>0$,

$$\lim_{l\to\infty}\limsup_{n\to\infty}\mathbf{P}(\langle Z^n_{l,k}\rangle(t)>\delta)=0.$$

Limit (4.38) follows by an application of the Lenglart–Rebolledo inequality; see, for example, Liptser and Shiryayev [16], Theorem 1.9.3. □

LEMMA 4.8. *As $n\to\infty$, the $(\hat X^n_0,\tilde X^n,H^n,G^n,M^{n,\varepsilon},M^n)$ converge in distribution in $\mathbb{R}\times\mathbb{D}_c(\mathbb{R}_+,\mathbb{R})^4\times\mathbb{D}(\mathbb{R}_+,\mathbb{R})$ to $(\hat X_0,\tilde X,H,G,M^\varepsilon,M)$. The latter random element has a tight distribution.*

PROOF. Let us show that the sequence $M^n$ is tight in $\mathbb{D}(\mathbb{R}_+,\mathbb{R})$ (for Skorohod's $J_1$-topology). By Lemma 3.1, the process $M^n_{n^2}$ is a $\mathbf{G}^n$-square integrable martingale with predictable quadratic variation process

$$\begin{aligned}\langle M^n_{n^2}\rangle(t)&=\int_{\mathbb{R}_+^2}\mathbf{1}_{\{s+x\le t\}}\,d\langle L^n_{n^2}\rangle(s,x)\\&=\int_0^t\frac{1}{n}\sum_{i=1}^{A^n(t-x)\wedge n^2}\mathbf{1}_{\{0<x\le\eta_i\}}\frac{1-F(x)}{(1-F(x-))^2}\,dF(x).\end{aligned}$$

Since the $A^n/n$ tend to $\mathbf{a}$ in probability uniformly over compact sets, by the law of large numbers and (4.10), for $T>0$ and $\delta>0$,

$$\lim_{n\to\infty}\mathbf{P}\biggl(\sup_{t\in[0,T]}\biggl|\langle M^n_{n^2}\rangle(t)-\int_0^t\mathbf{1}_{\{x>0\}}\mathbf{a}(t-x)\,dF'(x)\biggr|>\delta\biggr)=0.$$

By Theorem VI.5.17 in Jacod and Shiryaev [12], the sequence $M^n_{n^2}$ is tight in $\mathbb{D}(\mathbb{R}_+,\mathbb{R})$. Since the $A^n/n$ converge in probability to $\mathbf{a}$, the definition of $V^n(t,x)$ in (3.1), Donsker's theorem, Slutsky's lemma (Theorem D.9) and the continuous mapping principle imply that the sequence $\tilde V^n=(V^n(t,0),t\in\mathbb{R}_+)$ converges in distribution in $\mathbb{D}_c(\mathbb{R}_+,\mathbb{R})$ to $(\sqrt{F(0)(1-F(0))}W(\mathbf{a}(t)),t\in\mathbb{R}_+)$, where $(W(t),t\in\mathbb{R}_+)$ is a standard Brownian motion. Thus, this sequence is asymptotically tight in $\mathbb{D}_c(\mathbb{R}_+,\mathbb{R})$ (Theorem D.3). By Theorem D.7, the sequence $(\tilde V^n,M^n_{n^2})$ is asymptotically tight in $\mathbb{D}_c(\mathbb{R}_+,\mathbb{R})\times\mathbb{D}(\mathbb{R}_+,\mathbb{R})$. Since the map $(\mathbf{x},\mathbf{y})\to\mathbf{x}+\mathbf{y}$ is continuous as a map from $\mathbb{D}_c(\mathbb{R}_+,\mathbb{R})\times\mathbb{D}(\mathbb{R}_+,\mathbb{R})$ to $\mathbb{D}(\mathbb{R}_+,\mathbb{R})$, it follows that the sequence $\tilde V^n+M^n_{n^2}$ is asymptotically tight in $\mathbb{D}(\mathbb{R}_+,\mathbb{R})$, see Theorem D.4. By Ulam's theorem, $\tilde V^n+M^n_{n^2}$ is a tight random element of $\mathbb{D}(\mathbb{R}_+,\mathbb{R})$ for each $n$. Therefore, the sequence $\tilde V^n+M^n_{n^2}$ is tight in $\mathbb{D}(\mathbb{R}_+,\mathbb{R})$, van der Vaart and Wellner [23], Problem 1.3.9. Since, for $T>0$, by (3.9)–(3.11),

$$\mathbf{P}\Bigl(\sup_{t\in[0,T]}|M^n(t)-\breve V^n(A^n(t)/n,0)-M^n_{n^2}(t)|>0\Bigr)\le\mathbf{P}(A^n(T)>n^2)$$

and the latter probability tends to zero as $n \to \infty$, it follows that the sequence $M^n$ is tight in $\mathbb{D}(\mathbb{R}_+, \mathbb{R})$ (e.g., by Theorem VI.3.21 in Jacod and Shiryaev [12]).

By Lemma 4.7, the $(\hat{X}_0^n, \tilde{X}^n, H^n, G^n, M^{n,\varepsilon})$ converge in distribution in $\mathbb{R} \times \mathbb{D}_c(\mathbb{R}_+, \mathbb{R})^4$ to $(\hat{X}_0, \tilde{X}, H, G, M^\varepsilon)$, which is a separable, hence, tight, random element (recall the definition of $\tilde{X}$ in the statement of Theorem 2.2, Lemma 4.3, and (4.35)). Therefore, the sequence $(\hat{X}_0^n, \tilde{X}^n, H^n, G^n, M^{n,\varepsilon})$ is asymptotically tight in $\mathbb{R} \times \mathbb{D}_c(\mathbb{R}_+, \mathbb{R})^4$, see Theorem D.3. Since the sequence $M^n$ is tight in $\mathbb{D}(\mathbb{R}_+, \mathbb{R})$, it follows that the sequence $(\hat{X}_0^n, \tilde{X}^n, H^n, G^n, M^{n,\varepsilon}, M^n)$ is asymptotically tight in $\mathbb{R} \times \mathbb{D}_c(\mathbb{R}_+, \mathbb{R})^4 \times \mathbb{D}(\mathbb{R}_+, \mathbb{R})$, see Theorem D.7.

Since the topology of $\mathbb{R} \times \mathbb{D}(\mathbb{R}_+, \mathbb{R})^5$ is coarser than the topology of $\mathbb{R} \times \mathbb{D}_c(\mathbb{R}_+, \mathbb{R})^4 \times \mathbb{D}(\mathbb{R}_+, \mathbb{R})$, the sequence $(\hat{X}_0^n, \tilde{X}^n, H^n, G^n, M^{n,\varepsilon}, M^n)$ is asymptotically tight in $\mathbb{R} \times \mathbb{D}(\mathbb{R}_+, \mathbb{R})^5$. It is thus tight because the $(\hat{X}_0^n, \tilde{X}^n, H^n, G^n, M^{n,\varepsilon}, M^n)$ are random elements of $\mathbb{R} \times \mathbb{D}(\mathbb{R}_+, \mathbb{R})^5$, see van der Vaart and Wellner [23], Problem 1.3.9. Lemma 4.7 and Prohorov's theorem imply that the $(\hat{X}_0^n, \tilde{X}^n, H^n, G^n, M^{n,\varepsilon}, M^n)$ converge in distribution in $\mathbb{R} \times \mathbb{D}(\mathbb{R}_+, \mathbb{R})^5$ to $(\hat{X}_0, \tilde{X}, H, G, M^\varepsilon, M)$. On taking the set of Skorohod-continuous bounded functions as the separating subalgebra in Theorem D.6 we conclude that the sequence $(\hat{X}_0^n, \tilde{X}^n, H^n, G^n, M^{n,\varepsilon}, M^n)$ is asymptotically measurable in $\mathbb{R} \times \mathbb{D}_c(\mathbb{R}_+, \mathbb{R})^4 \times \mathbb{D}(\mathbb{R}_+, \mathbb{R})$. The latter property, coupled with the asymptotic tightness of this sequence in $\mathbb{R} \times \mathbb{D}_c(\mathbb{R}_+, \mathbb{R})^4 \times \mathbb{D}(\mathbb{R}_+, \mathbb{R})$, yields by Theorem D.5 the existence of a subsequence that converges in distribution in $\mathbb{R} \times \mathbb{D}_c(\mathbb{R}_+, \mathbb{R})^4 \times \mathbb{D}(\mathbb{R}_+, \mathbb{R})$ to a random element with a tight probability law. The limit must be the same as in $\mathbb{R} \times \mathbb{D}(\mathbb{R}_+, \mathbb{R})^5$, that is, it is $(\hat{X}_0, \tilde{X}, H, G, M^\varepsilon, M)$. It thus does not depend on a subsequence. □

4.3. *Completion of the proof of Theorem 3.1.* Let

$$M(t) = M_1(t) + M_2(t) \tag{4.39}$$

be the decomposition of the Gaussian martingale $M$ into the sum of a continuous Gaussian martingale $M_1 = (M_1(t), t \in \mathbb{R}_+)$ and a pure-jump Gaussian martingale $M_2 = (M_2(t), t \in \mathbb{R}_+)$ with jumps occurring at the jump times of $C(t)$; see Jacod and Shiryaev [12], Chapter II, Section 4d, Liptser and Shiryayev [16], Chapter 4, Section 9. Formally,

$$\begin{aligned} M_1(t) &= \int_0^t \mathbf{1}_{\{\Delta C(s)=0\}}\, dM(s), \\ M_2(t) &= \int_0^t \mathbf{1}_{\{\Delta C(s)>0\}}\, dM(s). \end{aligned} \tag{4.40}$$

Let also

$$M'^\varepsilon(t) = M(t) - M^\varepsilon(t). \tag{4.41}$$

We show that, for $T > 0$ and $\delta > 0$,

$$\lim_{\varepsilon\to 0} \mathbf{P}\Big(\sup_{t\in[0,T]} |M^\varepsilon(t) - M_2(t)| > \delta\Big) = 0, \tag{4.42a}$$

$$\lim_{\varepsilon\to 0} \mathbf{P}\Big(\sup_{t\in[0,T]} |M'^\varepsilon(t) - M_1(t)| > \delta\Big) = 0. \tag{4.42b}$$

Since by (4.17) $\Delta C(t) = \sum_{s+x=t} \Delta\mathbf{a}(s)\Delta F'(x)$, by (4.16),

$$M_2(t) = \int_{\mathbb{R}_+^2} \mathbf{1}_{\{\Delta\mathbf{a}(s)\Delta F(x)>0\}}\mathbf{1}_{\{s+x\le t\}}\, dL(s,x).$$

Hence, by (4.33) and (4.40),

$$M_2(t) - M^\varepsilon(t) = \int_{\mathbb{R}_+^2} \mathbf{1}_{\{0<\Delta\mathbf{a}(s)\Delta F(x)\le\varepsilon\}}\mathbf{1}_{\{s+x\le t\}}\, dL(s,x).$$

Lemma 4.2 implies that $(M_2(t) - M^\varepsilon(t), t \in \mathbb{R}_+)$ is a Gaussian martingale relative to the natural filtration with predictable quadratic variation process $(\int_{\mathbb{R}_+^2} \mathbf{1}_{\{0<\Delta\mathbf{a}(s)\Delta F(x)\le\varepsilon\}}\mathbf{1}_{\{s+x\le t\}}\, d\mathbf{a}(s)\, dF'(x), t \in \mathbb{R}_+)$. The latter function tends to zero as $\varepsilon \to 0$ for every $t \in \mathbb{R}_+$, so the Lenglart–Rebolledo inequality yields (4.42a). The limit in (4.42b) follows by (4.41), (4.42a) and the identity in (4.39).

We introduce processes $M'^{n,\varepsilon} = (M'^{n,\varepsilon}(t), t \in \mathbb{R}_+)$ by letting $M'^{n,\varepsilon}(t) = M^n(t) - M^{n,\varepsilon}(t)$. Let also $M'^\varepsilon = (M'^\varepsilon(t), t \in \mathbb{R}_+)$. By Lemma 4.8 and the continuous mapping principle, as $n \to \infty$, the $(\hat{X}_0^n, \tilde{X}^n, H^n, G^n, M^{n,\varepsilon}, M'^{n,\varepsilon})$ converge in distribution in $\mathbb{R}\times\mathbb{D}_c(\mathbb{R}_+,\mathbb{R})^4\times\mathbb{D}(\mathbb{R}_+,\mathbb{R})$ to $(\hat{X}_0, \tilde{X}, H, G, M^\varepsilon, M'^\varepsilon)$. The latter random element has a separable range. Therefore, by Theorem D.11,

$$\lim_{n\to\infty} d^*_{\mathrm{BL}_1}((\hat{X}_0^n, \tilde{X}^n, H^n, G^n, M^{n,\varepsilon}, M'^{n,\varepsilon}), (\hat{X}_0, \tilde{X}, H, G, M^\varepsilon, M'^\varepsilon)) = 0,$$

where $d^*_{\mathrm{BL}_1}$ is the distance on the space of mappings from $\Omega$ to $\mathbb{R}\times\mathbb{D}_c(\mathbb{R}_+,\mathbb{R})^4\times\mathbb{D}(\mathbb{R}_+,\mathbb{R})$ as defined in Appendix D. The limits in (4.42a) and (4.42b) imply that $\lim_{\varepsilon\to 0} d^*_{\mathrm{BL}_1}((\hat{X}_0, \tilde{X}, H, G, M^\varepsilon, M'^\varepsilon), (\hat{X}_0, \tilde{X}, H, G, M_2, M_1) = 0$, so

$$\lim_{\varepsilon\to 0}\limsup_{n\to\infty} d^*_{\mathrm{BL}_1}((\hat{X}_0^n, \tilde{X}^n, H^n, G^n, M^{n,\varepsilon}, M'^{n,\varepsilon}), (\hat{X}_0, \tilde{X}, H, G, M_2, M_1)) = 0.$$

Therefore, there exists a sequence $\varepsilon_n \to 0$ such that

$$\lim_{n\to\infty} d^*_{\mathrm{BL}_1}((\hat{X}_0^n, \tilde{X}^n, H^n, G^n, M^{n,\varepsilon_n}, M'^{n,\varepsilon_n}), (\hat{X}_0, \tilde{X}, H, G, M_2, M_1)) = 0,$$

which implies by Theorem D.11 that the $(\hat{X}_0^n, \tilde{X}^n, H^n, G^n, M^{n,\varepsilon_n}, M'^{n,\varepsilon_n})$ converge in distribution in $\mathbb{R}\times\mathbb{D}_c(\mathbb{R}_+,\mathbb{R})^4\times\mathbb{D}(\mathbb{R}_+,\mathbb{R})$ to $(\hat{X}_0, \tilde{X}, H, G, M_2, M_1)$. Since the process $M_1$ has continuous paths and convergence in Skorohod's $J_1$-topology to continuous functions is equivalent to compact convergence, by Corollary D.1 the $(\hat{X}_0^n, \tilde{X}^n, H^n, G^n, M^{n,\varepsilon_n}, M'^{n,\varepsilon_n})$ converge in distribution in $\mathbb{R}\times\mathbb{D}_c(\mathbb{R}_+,\mathbb{R})^5$ to $(\hat{X}_0, \tilde{X}, H, G, M_2, M_1)$. By the continuous mapping principle, the $(\hat{X}_0^n, \tilde{X}^n, H^n, G^n - M^{n,\varepsilon_n} - M'^{n,\varepsilon_n})$ converge in distribution in $\mathbb{R}\times\mathbb{D}_c(\mathbb{R}_+,\mathbb{R})^3$ to $(\hat{X}_0, \tilde{X}, H, G - M_1 - M_2)$. Recalling that $M^n = M^{n,\varepsilon} + M'^{n,\varepsilon}$, $Z^n = G^n - M^n$, $M = M_1 + M_2$, and $Z = G - M$ completes the proof of Theorem 3.1.

## APPENDIX A: EXISTENCE AND UNIQUENESS FOR THE QUEUEING EQUATIONS

We prove existence and uniqueness of solutions to (2.2a)–(2.2c).

LEMMA A.1. *Given a nondecreasing nonnegative integer-valued process $E^n$ with trajectories in $\mathbb{D}(\mathbb{R}_+,\mathbb{R})$, a nonnegative integer-valued random variable $Q_0^n$, and sequences of nonnegative random variables $\{\eta_i\}$ and $\{\tilde{\eta}_i\}$, there exist nonnegative integer-valued processes $Q^n$ and $\tilde{Q}^n$ and nonnegative integer-valued nondecreasing process $A^n$ with trajectories in $\mathbb{D}(\mathbb{R}_+,\mathbb{R})$ such that equations (2.2a)–(2.2c) are satisfied. The process $\tilde{Q}^n$ is specified uniquely. If $\eta_i > 0$ for all $i$, then the processes $A^n$ and $Q^n$ are specified uniquely.*

PROOF. We fix $\omega\in\Omega$ throughout. Obviously, $\tilde{Q}^n$ is specified uniquely by (2.2c). For $Q^n$ and $A^n$, we start with the case where $Q_0^n < n$. First, the existence issue is addressed. Introduce the process $Q^{n,1} = (Q^{n,1}(t), t\in\mathbb{R}_+)$ as the process of the number of customers in the infinite server with the same arrival and service times and the initial number of customers, that is,

$$Q^{n,1}(t) = \tilde{Q}^n(t) + E^n(t) - \int_0^t\int_0^t \mathbf{1}_{\{s+x\le t\}}\, d\sum_{i=1}^{E^n(s)} \mathbf{1}_{\{\eta_i\le x\}}. \tag{A.1}$$

Let $\tau^{n,1} = \inf\{t : Q^{n,1}(t)\ge n\}\le\infty$. Suppose $\tau^{n,1} > 0$. We define $Q^n(t) = Q^{n,1}(t)$ and $A^n(t) = E^n(t)$ for $t < \tau^{n,1}$. Obviously, equations (2.2a) and (2.2b) are satisfied for $t < \tau^{n,1}$.

If $\tau^{n,1} = \infty$, the proof of existence is over. Suppose the contrary. Clearly, $\tau^{n,1}$ is a jump time of $E^n$. We choose $A^n(\tau^{n,1})$ such that

$$\begin{aligned}\sum_{i=E^n(\tau^{n,1}-)+1}^{A^n(\tau^{n,1})} \mathbf{1}_{\{\eta_i>0\}} \\ = n - \left(Q^n(\tau^{n,1}-) + \Delta\tilde{Q}^n(\tau^{n,1}) - \sum_{i=1}^{E^n(\tau^{n,1}-)} \mathbf{1}_{\{\tau_i^n+\eta_i=\tau^{n,1}\}}\right).\end{aligned} \tag{A.2}$$

It exists and is not greater than $E^n(\tau^{n,1})$ due to the facts that $Q^n(\tau^{n,1}-) < n$ and that $Q^{n,1}(\tau^{n,1}) \geq n$, so by (A.1),

$$Q^n(\tau^{n,1}-) + \Delta\tilde{Q}^n(\tau^{n,1}) + \Delta E^n(\tau^{n,1}) - \sum_{i=1}^{E^n(\tau^{n,1}-)} \mathbf{1}_{\{\tau_i^n+\eta_i=\tau^{n,1}\}} - \sum_{i=E^n(\tau^{n,1}-)+1}^{E^n(\tau^{n,1})} \mathbf{1}_{\{\eta_i=0\}} \geq n. \tag{A.3}$$

We also let

$$Q^n(\tau^{n,1}) = n + E^n(\tau^{n,1}) - A^n(\tau^{n,1}). \tag{A.4}$$

Hence, $Q^n(\tau^{n,1}) \geq n$, so (2.2b) holds for $t = \tau^{n,1}$. By (A.2) and (A.4),

$$Q^n(\tau^{n,1}) = Q^n(\tau^{n,1}-) + \Delta E^n(\tau^{n,1}) + \Delta\tilde{Q}^n(\tau^{n,1}) - \sum_{i=1}^{A^n(\tau^{n,1})} \mathbf{1}_{\{\tau_i^n+\eta_i=\tau^{n,1}\}}. \tag{A.5}$$

We obtain that equations (2.2a) and (2.2b) are satisfied on $[0, \tau^{n,1}]$. We next extend the construction past $\tau^{n,1}$. Let us define a nondecreasing nonnegative integer-valued process $A^{n,2} = (A^{n,2}(t), t \in \mathbb{R}_+)$ with right-continuous trajectories admitting left-hand limits and with $A^{n,2}(0) = 0$ by specifying its jumps $\Delta A^{n,2}(t)$ for $t > 0$ as follows:

$$\Delta A^{n,2}(t) = \inf\Biggl\{k \in \mathbb{N} : \sum_{i=A^{n,2}(t-)+1}^{A^{n,2}(t-)+k} \mathbf{1}_{\{\eta_i^{n,2}>0\}} \geq \sum_{i=1}^{A^{n,2}(t-)} \mathbf{1}_{\{\tau_i^{n,2}+\eta_i^{n,2}=t\}} - \Delta\tilde{Q}^{n,2}(t)\Biggr\}, \tag{A.6}$$

where $\tau_i^{n,2} = \inf\{t : A^{n,2}(t) \geq i\}$, $\eta_i^{n,2} = \eta_{i+A^n(\tau^{n,1})}$, and

$$\tilde{Q}^{n,2}(t) = \tilde{Q}^n(t+\tau^{n,1}) + \sum_{i=1}^{A^n(\tau^{n,1})} \mathbf{1}_{\{\tau_i^n+\eta_i>t+\tau^{n,1}\}}. \tag{A.7}$$

It follows by (A.1) and (A.2) on recalling that $Q^n(\tau^{n,1}-) = Q^{n,1}(\tau^{n,1}-)$ and $A^n(\tau^{n,1}-) = E^n(\tau^{n,1}-)$ that $\tilde{Q}^{n,2}(0) = n$ and $\tilde{Q}^{n,2}(t) \leq n$ for $t \geq 0$ as the right-hand side of (A.7) is decreasing in $t$.

In words, $A^{n,2}(t)$ is the number of customers that enter service by time $t$ for an $n$-server queue that always has a nonzero queue length and for which residual service times of customers in service at time 0 are those of customers in the $n$-server queue under consideration at time $\tau^{n,1}$ and the service times of customers entering service after time zero are equal to the service times of customers entering service after time $\tau^{n,1}$ in the queue in question. We note that this process can run away to infinity on finite time. (For instance, if $\tilde{Q}^{n,2}(t)$ has its first jump at $s>0$ and $\eta_i^{n,2}=0$ for all $i$, then $A^{n,2}(t)=\infty$ for $t\ge s$.) We therefore introduce $\tilde{\tau}^{n,1}=\sup\{t: A^{n,2}(t)<\infty\}$ and note that $\tilde{\tau}^{n,1}>0$ by the right continuity of $\tilde{Q}^{n,2}$.

By (A.6), we have for $t<\tilde{\tau}^{n,1}$ that

$$\Delta\sum_{i=1}^{A^{n,2}(t)}\mathbf{1}_{\{\tau_i^{n,2}+\eta_i^{n,2}>t\}}=-\Delta\tilde{Q}^{n,2}(t),$$

so, since $\tilde{Q}^{n,2}(0)=n$,

$$\sum_{i=1}^{A^{n,2}(t)}\mathbf{1}_{\{\tau_i^{n,2}+\eta_i^{n,2}>t\}}+\tilde{Q}^{n,2}(t)=n. \tag{A.8}$$

Introduce

$$\sigma^{n,1}=\inf\{t: A^{n,2}(t)>Q^n(\tau^{n,1})-n+(E^n(t+\tau^{n,1})-E^n(\tau^{n,1}))\}. \tag{A.9}$$

By right continuity of $A^{n,2}$, $\sigma^{n,1}>0$. Also, $\sigma^{n,1}\le\tilde{\tau}^{n,1}$. Intuitively, $\sigma^{n,1}$ is the time when $A^{n,2}(t)$ and $A^n(t+\tau^{n,1})-A^n(\tau^{n,1})$ diverge. To substantiate this, define for $t\in[\tau^{n,1},\tau^{n,1}+\sigma^{n,1})$,

$$\begin{aligned}Q^n(t)&=Q^n(\tau^{n,1})-n+\tilde{Q}^{n,2}(t-\tau^{n,1})+(E^n(t)-E^n(\tau^{n,1}))\\&\quad-\sum_{i=1}^{A^{n,2}(t-\tau^{n,1})}\mathbf{1}_{\{\tau_i^{n,2}+\eta_i^{n,2}\le t-\tau^{n,1}\}}.\end{aligned} \tag{A.10}$$

By (A.8), we can also write

$$Q^n(t)=Q^n(\tau^{n,1})+(E^n(t)-E^n(\tau^{n,1}))-A^{n,2}(t-\tau^{n,1}). \tag{A.11}$$

Therefore, $Q^n(t)\ge n$ by (A.9), and by (A.4),

$$Q^n(t)=n+E^n(t)-A^n(\tau^{n,1})-A^{n,2}(t-\tau^{n,1}). \tag{A.12}$$

On letting

$$A^n(t)=A^n(\tau^{n,1})+A^{n,2}(t-\tau^{n,1}), \tag{A.13}$$

we can see that (2.2b) holds on $[\tau^{n,1},\tau^{n,1}+\sigma^{n,1})$. To obtain (2.2a), we substitute (A.4) in (A.10) and recall the definitions of $\tau_i^{n,2}$, $\eta_i^{n,2}$ and $\tilde{Q}^{n,2}(t)$.

If $\sigma^{n,1}=\infty$, the proof of existence is over. If $\sigma^{n,1}<\infty$, we let

$$\Delta A^n(\tau^{n,1}+\sigma^{n,1})=Q^n((\tau^{n,1}+\sigma^{n,1})-)-n+\Delta E^n(\tau^{n,1}+\sigma^{n,1}) \tag{A.14}$$

and

$$\Delta Q^n(\tau^{n,1}+\sigma^{n,1}) = \Delta E^n(\tau^{n,1}+\sigma^{n,1}) - \sum_{i=1}^{A^{n,2}(\sigma^{n,1}-)+\Delta A^n(\tau^{n,1}+\sigma^{n,1})} \mathbf{1}_{\{\tau_i^{n,2}+\eta_i^{n,2}=\sigma^{n,1}\}} + \Delta\tilde{Q}^{n,2}(\sigma^{n,1}). \tag{A.15}$$

We note that by (A.11) and (A.14),

$$\Delta A^n(\tau^{n,1}+\sigma^{n,1}) = Q^n(\tau^{n,1}) - n + (E^n(\tau^{n,1}+\sigma^{n,1}) - E^n(\tau^{n,1})) - A^{n,2}(\sigma^{n,1}-), \tag{A.16}$$

so by (A.9)

$$0 \le \Delta A^n(\tau^{n,1}+\sigma^{n,1}) < \Delta A^{n,2}(\sigma^{n,1}). \tag{A.17}$$

Let us show that $Q^n(\tau^{n,1}+\sigma^{n,1}) < n$. By (A.6) and the right-hand inequality in (A.17),

$$\sum_{i=A^{n,2}(\sigma^{n,1}-)+1}^{A^{n,2}(\sigma^{n,1}-)+\Delta A^n(\tau^{n,1}+\sigma^{n,1})} \mathbf{1}_{\{\eta_i^{n,2}>0\}} < \sum_{i=1}^{A^{n,2}(\sigma^{n,1}-)} \mathbf{1}_{\{\tau_i^{n,2}+\eta_i^{n,2}=\sigma^{n,1}\}} - \Delta\tilde{Q}^{n,2}(\sigma^{n,1}),$$

so

$$\Delta A^n(\tau^{n,1}+\sigma^{n,1}) < \sum_{i=1}^{A^{n,2}(\sigma^{n,1}-)+\Delta A^n(\tau^{n,1}+\sigma^{n,1})} \mathbf{1}_{\{\tau_i^{n,2}+\eta_i^{n,2}=\sigma^{n,1}\}} - \Delta\tilde{Q}^{n,2}(\sigma^{n,1}), \tag{A.18}$$

and by (A.15) and (A.14),

$$\Delta Q^n(\tau^{n,1}+\sigma^{n,1}) < \Delta E^n(\tau^{n,1}+\sigma^{n,1}) - \Delta A^n(\tau^{n,1}+\sigma^{n,1}) = n - Q^n((\tau^{n,1}+\sigma^{n,1})-). \tag{A.19}$$

Hence, $Q^n(\tau^{n,1}+\sigma^{n,1}) < n$.

Also, $Q^n(\tau^{n,1}+\sigma^{n,1}) \ge 0$ which is shown as follows. By (A.15) and (A.16),

$$\begin{aligned} Q^n(\tau^{n,1}+\sigma^{n,1}) &\ge Q^n((\tau^{n,1}+\sigma^{n,1})-) + \Delta E^n(\tau^{n,1}+\sigma^{n,1}) - \Delta A^n(\tau^{n,1}+\sigma^{n,1}) \\ &\quad - \sum_{i=1}^{A^{n,2}(\sigma^{n,1}-)} \mathbf{1}_{\{\tau_i^{n,2}+\eta_i^{n,2}=\sigma^{n,1}\}} + \Delta\tilde{Q}^{n,2}(\sigma^{n,1}) \\ &= n - \sum_{i=1}^{A^{n,2}(\sigma^{n,1}-)} \mathbf{1}_{\{\tau_i^{n,2}+\eta_i^{n,2}=\sigma^{n,1}\}} + \Delta\tilde{Q}^{n,2}(\sigma^{n,1}). \end{aligned}$$

Since $-\Delta\tilde{Q}^{n,2}(\sigma^{n,1}) \le \tilde{Q}^{n,2}(\sigma^{n,1}-)$, we have by (A.8) that

$$\sum_{i=1}^{A^{n,2}(\sigma^{n,1}-)} \mathbf{1}_{\{\tau_i^{n,2}+\eta_i^{n,2}=\sigma^{n,1}\}} - \Delta\tilde{Q}^{n,2}(\sigma^{n,1}) \le \sum_{i=1}^{A^{n,2}(\sigma^{n,1}-)} \mathbf{1}_{\{\tau_i^{n,2}+\eta_i^{n,2}\le\sigma^{n,1}\}} + \tilde{Q}^{n,2}(\sigma^{n,1}-) = n.$$

The required inequality has been proven.

By (A.14), (2.2b) for $t < \tau^{n,1}+\sigma^{n,1}$, and the inequality $Q^n((\tau^{n,1}+\sigma^{n,1})-) \ge n$, $\Delta A^n(\tau^{n,1}+\sigma^{n,1}) = E^n(\tau^{n,1}+\sigma^{n,1}) - A^n((\tau^{n,1}+\sigma^{n,1})-)$. We thus have (2.2b) for $t = \tau^{n,1}+\sigma^{n,1}$. Equation (2.2a) for $t = \tau^{n,1}+\sigma^{n,1}$ follows by (A.15) and (A.7). Existence has been proven on $[0, \tau^{n,1}+\sigma^{n,1}]$.

At time $\tau^{n,1}+\sigma^{n,1}$, we are in a similar situation to the one we faced at $t = 0$. We can thus define $\tau^{n,2}$ and $\sigma^{n,2}$ and proceed until we get to the desired time $t$. This is bound to happen after a finite number of steps the reason being that in any interval $[\sum_{j=1}^{i-1}(\tau^{n,j}+\sigma^{n,j}), \sum_{j=1}^{i-1}(\tau^{n,j}+\sigma^{n,j})+\tau^{n,i}]$, where $\sigma^{n,1} = 0$ by definition, we have at least one upward jump of $E^n$. Since $E^n$ is finite-valued, there are only finitely many such intervals on any interval $[0, t]$.

We now prove uniqueness. Thus, we suppose that $Q^n$ and $A^n$ satisfy (2.2a) and (2.2b) and that $\eta_i > 0$ for all $i$. Assume for the moment that $\tau^{n,1} > 0$. If $Q^n$ and $A^n$ differed from those defined above on $[0, \tau^{n,1})$, then there would exist $t' < \tau^{n,1}$ with $Q^n(t') \ge n$. We note that $t' > 0$ which is checked as follows. By (2.2a) and (2.2b) for $t = 0$,

$$\begin{aligned} Q^n(0) &= \tilde{Q}^n(0) + (Q^n(0)-n)^+ + A^n(0) - \sum_{i=1}^{A^n(0)} \mathbf{1}_{\{\eta_i=0\}} \\ &\le \tilde{Q}^n(0) + (Q^n(0)-n)^+ + E^n(0) - \sum_{i=1}^{E^n(0)} \mathbf{1}_{\{\eta_i=0\}} \\ &= (Q^n(0)-n)^+ + Q^{n,1}(0) < (Q^n(0)-n)^+ + n, \end{aligned}$$

so $Q^n(0) < n$.

For $\tilde{\tau}^{n,1} = \inf\{t < \tau^{n,1} : Q^n(t) \ge n\} \in (0, \infty)$, we would have that $Q^n(t) = Q^{n,1}(t)$ and $A^n(t) = E^n(t)$ when $t < \tilde{\tau}^{n,1}$. Also, $A^n(\tilde{\tau}^{n,1}) \le E^n(\tilde{\tau}^{n,1})$. Since by (2.2a) and (A.1),

$$Q^n(\tilde{\tau}^{n,1}) = Q^{n,1}(\tilde{\tau}^{n,1}) + \sum_{i=A^n(\tilde{\tau}^{n,1})+1}^{E^n(\tilde{\tau}^{n,1})} \mathbf{1}_{\{\eta_i=0\}} < n + E^n(\tilde{\tau}^{n,1}) - A^n(\tilde{\tau}^{n,1}),$$

we arrive at a contradiction with (2.2b). We have thus proved that a solution to (2.2a) and (2.2b) is specified uniquely on $[0,\tau^{n,1})$.

Next, (2.2a) and (2.2b) imply (A.5), so by (A.3) and the inequality $E^n(\tau^{n,1}) \geq A^{n,1}(\tau^{n,1})$ we have that $Q^n(\tau^{n,1}) \geq n$, and (2.2b) implies (A.4). Also, (A.4) and (A.5) imply (A.2). Since $\eta_i > 0$, the latter condition specifies $A^n(\tau^{n,1})$ uniquely. Uniqueness on $[0,\tau^{n,1}]$ follows.

We now show that $Q^n(t) \geq n$ for $t \in [\tau^{n,1}, \tau^{n,1}+\sigma^{n,1})$. Let $\tilde{\sigma}^{n,1} = \inf\{t > \tau^{n,1} : Q^n(t) < n\}$ and suppose that $\tilde{\sigma}^{n,1} < \tau^{n,1}+\sigma^{n,1}$. By right continuity of $Q^n$, we have that $\tilde{\sigma}^{n,1} > \tau^{n,1}$. By (2.2b), for $t \in [\tau^{n,1}, \tilde{\sigma}^{n,1})$,

$$A^n(t) - A^n(\tau^{n,1}) = E^n(t) - E^n(\tau^{n,1}) + Q^n(\tau^{n,1}) - Q^n(t). \tag{A.20}$$

Since by (2.2a) and (A.7) for $t \geq \tau^{n,1}$,

$$\begin{aligned} Q^n(t) = Q^n(\tau^{n,1}) - n + \tilde{Q}^{n,2}(t-\tau^{n,1}) + (E^n(t) - E^n(\tau^{n,1})) \\ - \sum_{i=A^n(\tau^{n,1})+1}^{A^n(t)} \mathbf{1}_{\{\tau_i^n+\eta_i \leq t\}}, \end{aligned} \tag{A.21}$$

we have that, for $t \in [\tau^{n,1}, \tilde{\sigma}^{n,1})$,

$$\sum_{i=A^n(\tau^{n,1})+1}^{A^n(t)} \mathbf{1}_{\{\tau_i^n+\eta_i > t\}} = n - \tilde{Q}^{n,2}(t-\tau^{n,1}), \tag{A.22}$$

so that

$$\Delta \sum_{i=A^n(\tau^{n,1})+1}^{A^n(t)} \mathbf{1}_{\{\tau_i^n+\eta_i > t\}} = -\Delta\tilde{Q}^{n,2}(t-\tau^{n,1}).$$

Consequently, if $t$ is a jump time of $A^n$, then

$$\sum_{i=A^n(t-)+1}^{A^n(t)} \mathbf{1}_{\{\eta_i > 0\}} = \sum_{i=A^n(\tau^{n,1})+1}^{A^n(t-)} \mathbf{1}_{\{\tau_i^n+\eta_i = t\}} - \Delta\tilde{Q}^{n,2}(t-\tau^{n,1}).$$

It follows by (A.6) and the assumption that $\eta_i > 0$ that $A^n(t) - A^n(\tau^{n,1}) = A^{n,2}(t-\tau^{n,1})$, which implies that $Q^n(t)$ is given by the right-hand side of (A.10) (or the right-hand side of (A.11)) for $t \in [\tau^{n,1}, \tilde{\sigma}^{n,1})$. Besides, $\tilde{\sigma}^{n,1} - \tau^{n,1} \leq \tilde{\tau}^{n,1}$.

We also have by (2.2b) since $Q^n(\tilde{\sigma}^{n,1}) < n$, that $A^n(\tilde{\sigma}^{n,1}) - A^n(\tau^{n,1}) = E^n(\tilde{\sigma}^{n,1}) - E^n(\tau^{n,1}) + Q^n(\tau^{n,1}) - n$, so by (A.21) with $t = \tilde{\sigma}^{n,1}$

$$\begin{aligned} \sum_{i=A^n(\tau^{n,1})+1}^{A^n(\tilde{\sigma}^{n,1})} \mathbf{1}_{\{\tau_i^n+\eta_i > \tilde{\sigma}^{n,1}\}} &= Q^n(\tilde{\sigma}^{n,1}) - \tilde{Q}^{n,2}(\tilde{\sigma}^{n,1} - \tau^{n,1}) \\ &< n - \tilde{Q}^{n,2}(\tilde{\sigma}^{n,1} - \tau^{n,1}). \end{aligned}$$

It follows by (A.22) on noting that $\tilde{\sigma}^{n,1}$ is a jump time of $A^n$ that

$$\sum_{A^n(\tilde{\sigma}^{n,1}-)+1}^{A^n(\tilde{\sigma}^{n,1})} \mathbf{1}_{\{\eta_i>0\}} < \sum_{i=A^n(\tau^{n,1})+1}^{A^n(\tilde{\sigma}^{n,1}-)} \mathbf{1}_{\{\tau_i^n+\eta_i=\tilde{\sigma}^{n,1}\}} - \Delta\tilde{Q}^{n,2}(\tilde{\sigma}^{n,1}-\tau^{n,1}),$$

so by the fact that $A^n(t)-A^n(\tau^{n,1}) = A^{n,2}(t-\tau^{n,1})$ for $t<\tilde{\sigma}^{n,1}$ and (A.6), $\Delta A^n(\tilde{\sigma}^{n,1}) < \Delta A^{n,2}(\tilde{\sigma}^{n,1}-\tau^{n,1})$. Therefore,

$$\sum_{i=A^n(\tau^{n,1})+1}^{A^n(\tilde{\sigma}^{n,1})} \mathbf{1}_{\{\tau_i^n+\eta_i\le\tilde{\sigma}^{n,1}\}} \le \sum_{i=1}^{A^{n,2}(\tilde{\sigma}^{n,1}-\tau^{n,1})} \mathbf{1}_{\{\tau_i^{n,2}+\eta_i^{n,2}\le\tilde{\sigma}^{n,1}-\tau^{n,1}\}}.$$

Since the right-hand side of (A.10) is not less than $n$ at $\tilde{\sigma}^{n,1}$, so is the right-hand side of (A.21), which contradicts the definition of $\tilde{\sigma}^{n,1}$. The obtained contradiction proves that $Q^n(t)\ge n$ on $[\tau^{n,1},\sigma^{n,1})$. Consequently, (A.20) holds for $t$ from this interval. By (A.21), (A.22) holds for those $t$, hence, as we have seen, (A.13) holds. By (A.13) and (A.20), $Q^n(t)$ is given by the right-hand side of (A.11). Thus, $Q^n$ and $A^n$ are specified uniquely on $[0,\tau^{n,1}+\sigma^{n,1})$.

For $t=\tau^{n,1}+\sigma^{n,1}$, we have that (2.2a) implies (A.15) (with possibly different $\Delta A^n(\tau^{n,1}+\sigma^{n,1})$). By (2.2b) on $[0,\tau^{n,1}+\sigma^{n,1}]$,

$$\begin{aligned}\Delta A^n(\tau^{n,1}+\sigma^{n,1}) = \Delta E^n(\tau^{n,1}+\sigma^{n,1}) - (Q^n(\tau^{n,1}+\sigma^{n,1})-n)^+ \\ + (Q^n((\tau^{n,1}+\sigma^{n,1})-)-n).\end{aligned} \tag{A.23}$$

In particular, $\Delta A^n(\tau^{n,1}+\sigma^{n,1})$ is not greater than the right-hand side of (A.14), so by (A.11) it is not greater than the right-hand side of (A.16), and we obtain (A.17) by applying (A.9). A similar argument to the one we used above yields (A.18). Then, in analogy with (A.19) on taking into account (A.23),

$$\begin{aligned}\Delta Q^n(\tau^{n,1}+\sigma^{n,1}) &< \Delta E^n(\tau^{n,1}+\sigma^{n,1}) - \Delta A^n(\tau^{n,1}+\sigma^{n,1}) \\ &= n\vee Q^n(\tau^{n,1}+\sigma^{n,1}) - Q^n((\tau^{n,1}+\sigma^{n,1})-),\end{aligned}$$

so that $Q^n(\tau^{n,1}+\sigma^{n,1}) < n\vee Q^n(\tau^{n,1}+\sigma^{n,1})$. Hence, $Q^n(\tau^{n,1}+\sigma^{n,1})<n$. By (A.23), we have (A.14). Having established uniqueness on $[0,\tau^{n,1}+\sigma^{n,1}]$, we can use this argument repeatedly until a given $t$ is reached.

We have thus established existence and uniqueness for the case where $Q_0^n<n$. If $Q_0^n\ge n$, then we can use an analogous argument to the one employed when $\tau^{n,1}=0$. □

We now provide an example of nonuniqueness in the case where there are zero service times. Consider a single server queue (so $n=1$) with no customers initially. There are two arrivals: the first arrival occurs at $t=1$ and the second at $t=2$. The service times are $\eta_1=2$ and $\eta_2=0$, respectively. Thus, $E^n(t)=\mathbf{1}_{\{t\geq 1\}}+\mathbf{1}_{\{t\geq 2\}}$. The equations for $Q^n(t)$ and $A^n(t)$ are as follows

$$\begin{aligned}
Q^n(t)&=\mathbf{1}_{\{t\geq 1\}}+\mathbf{1}_{\{t\geq 2\}}-\mathbf{1}_{\{A^n(t)\geq 1\}}\mathbf{1}_{\{1+2\leq t\}}-\mathbf{1}_{\{A^n(t)\geq 2\}}\mathbf{1}_{\{2+0\leq t\}},\\
A^n(t)&=\mathbf{1}_{\{t\geq 1\}}+\mathbf{1}_{\{t\geq 2\}}-(Q^n(t)-1)^+.
\end{aligned}$$

They admit two sets of solutions: $Q^n(t)=\mathbf{1}_{\{1\leq t<3\}}, A^n(t)=\mathbf{1}_{\{1\leq t<2\}}+2\mathbf{1}_{\{t\geq 2\}}$ and $Q^n(t)=\mathbf{1}_{\{1\leq t<2\}}+2\mathbf{1}_{\{2\leq t<3\}}, A^n(t)=\mathbf{1}_{\{1\leq t<3\}}+2\mathbf{1}_{\{t\geq 3\}}$.

## APPENDIX B: CONTINUITY PROPERTIES OF CONVOLUTION EQUATIONS

Let $B=(B(x),x\in\mathbb{R}_+)$ be a distribution function on $\mathbb{R}_+$ with $B(0)<1$. Given $T>0$, let $\mathbb{L}_\infty([0,T],\mathbb{R})$ denote the Banach space of $\mathbb{R}$-valued bounded Borel measurable functions on $[0,T]$ which is equipped with the uniform norm $\|\mathbf{x}\|_\infty=\sup_{t\in[0,T]}|\mathbf{x}(t)|$, where $\mathbf{x}=(\mathbf{x}(t),t\in[0,T])$. Consider the equation

$$\mathbf{y}(t)=\mathbf{x}(t)+\int_0^t f(\mathbf{y}(t-s),t-s)\,dB(s), \tag{B.1}$$

where $\mathbf{x}=(\mathbf{x}(t),t\in[0,T])\in\mathbb{L}_\infty([0,T],\mathbb{R})$ and the function $f:\mathbb{R}\times\mathbb{R}_+\to\mathbb{R}$ is Borel measurable.

LEMMA B.1. *If $|f(y,t)|\leq|y|$ and $|f(y_1,t)-f(y_2,t)|\leq|y_1-y_2|$ for all $y$, $y_1$, $y_2$, and $t$ from the domain, then for every $\mathbf{x}\in\mathbb{L}_\infty([0,T],\mathbb{R})$ there exists a unique $\mathbf{y}=(\mathbf{y}(t),t\in[0,T])\in\mathbb{L}_\infty([0,T],\mathbb{R})$ which satisfies equation (B.1). This solution can be obtained by the method of successive approximations and there exists a function $\rho(t)$ which only depends on the function $B$ such that* $\sup_{t\in[0,T]}|\mathbf{y}(t)|\leq\rho(T)\sup_{t\in[0,T]}|\mathbf{x}(t)|$. *The latter bound holds more generally if*

$$|\mathbf{y}(t)|\leq|\mathbf{x}(t)|+\int_0^t|\mathbf{y}(t-s)|\,dB(s).$$

*If, in addition, $\mathbf{x}$ is right-continuous with left-hand limits (respectively, admits limits on the right and limits on the left) and $f(y,t)$ is right-continuous with left-hand limits (respectively, admits limits on the right and limits on the left) in the second argument, then $\mathbf{y}$ is right-continuous with left-hand limits (respectively, admits limits on the right and limits on the left). If $\mathbf{x}$ is continuous and either $f(y,t)$ is continuous in the second argument, $f(0,0)=0$, and $\mathbf{x}(0)=0$, or $f(y,t)$ admits limits on the right and limits on the left in the second argument, and $B$ is continuous on $[0,T]$ with $B(0)=0$, then $\mathbf{y}$ is continuous.*

PROOF. Let $t_0 \in (0,T]$ be such that $B(t_0) < 1$. Define the map $\phi_{\mathbf{x}}$ by

$$\phi_{\mathbf{x}}(\mathbf{z})(t) = \mathbf{x}(t) + \int_0^t f(\mathbf{z}(t-s), t-s)\, dB(s). \tag{B.2}$$

The boundedness condition on $f$ implies that $\phi_{\mathbf{x}}$ is an operator on $\mathbb{L}_\infty([0,t_0],\mathbb{R})$. We show it is a contraction. Since $|\phi_{\mathbf{x}}(\mathbf{z})(t) - \phi_{\mathbf{x}}(\mathbf{z}')(t)| \le \int_0^t |\mathbf{z}(t-s) - \mathbf{z}'(t-s)|\, dB(s)$, we have that $\sup_{t\in[0,t_0]} |\phi_{\mathbf{x}}(\mathbf{z})(t) - \phi_{\mathbf{x}}(\mathbf{z}')(t)| \le B(t_0) \sup_{t\in[0,t_0]} |\mathbf{z}(t) - \mathbf{z}'(t)|$, which proves the claim because $B(t_0) < 1$. Since $\mathbb{L}_\infty([0,t_0],\mathbb{R})$ is a Banach space, the operator $\phi_{\mathbf{x}}$ has a unique fixed point so that the equation $\mathbf{y}(t) = \mathbf{x}(t) + \int_0^t f(\mathbf{y}(t-s), t-s)\, dB(s)$ has a unique solution for $t \in [0,t_0]$, which is obtained by the method of successive approximations.

Next, on introducing $\mathbf{y}_{t_0}(t) = \mathbf{y}(t+t_0)$, $\mathbf{x}_{t_0}(t) = \mathbf{x}(t+t_0) + \int_t^{t+t_0} f(\mathbf{y}(t+t_0-s), t+t_0-s)\, dB(s)$, and $f_{t_0}(u,v) = f(u, t_0+v)$, we can write (B.1) for $t \in [0, T-t_0]$ in the form

$$\mathbf{y}_{t_0}(t) = \mathbf{x}_{t_0}(t) + \int_0^t f_{t_0}(\mathbf{y}_{t_0}(t-s), t-s)\, dB(s). \tag{B.3}$$

The function $\mathbf{x}_{t_0} = (\mathbf{x}_{t_0}(t), t \in [0, T-t_0])$ is uniquely specified by $\mathbf{x}$ and the values of $\mathbf{y}(t)$ for $t \le t_0$. The preceding argument applied to (B.3) shows that given $\mathbf{y}(t)$ for $t \in [0,t_0]$, there exists a unique extension of $\mathbf{y}(t)$ to the interval $[t_0, 2t_0 \wedge T]$ that satisfies the equation. The method of successive approximations converges to this solution. By applying this argument repeatedly, we deduce existence and uniqueness of a solution in $\mathbb{L}_\infty([0,T],\mathbb{R})$. This solution is obtained by the method of successive approximations.

If $\mathbf{x}$ is right-continuous with left-hand limits (respectively, admits limits on the right and limits on the left) and $f(y,t)$ is right-continuous with left-hand limits (respectively, admits limits on the right and limits on the left) in $t$, then the successive approximations starting from the zero function are right-continuous with left-hand limits (respectively, admit limits on the right and limits on the left), so $\mathbf{y}$ is right-continuous with left-hand limits (respectively, admits limits on the right and limits on the left). A similar argument applies in the case where $\mathbf{x}$ is continuous.

Suppose that

$$|\mathbf{y}(t)| \le |\mathbf{x}(t)| + \int_0^t |\mathbf{y}(t-s)|\, dB(s) \qquad \text{for } t \in [0,T].$$

Then

$$\sup_{t\in[0,t_0]} |\mathbf{y}(t)| \le \sup_{t\in[0,t_0]} |\mathbf{x}(t)| + \sup_{t\in[0,t_0]} |\mathbf{y}(t)| B(t_0),$$

so

$$\sup_{t\in[0,t_0]} |\mathbf{y}(t)| \le \frac{1}{1-B(t_0)} \sup_{t\in[0,t_0]} |\mathbf{x}(t)|.$$

Next,

$$\sup_{t\in[t_0,(2t_0)\wedge T]}|\mathbf{y}(t)| \le \sup_{t\in[t_0,(2t_0)\wedge T]}|\mathbf{x}(t)| + \sup_{t\in[0,t_0]}|\mathbf{y}(t)| + \sup_{t\in[t_0,(2t_0)\wedge T]}|\mathbf{y}(t)|(1-B(t_0)).$$

Thus,

$$\sup_{t\in[0,(2t_0)\wedge T]}|\mathbf{y}(t)| \le \sup_{t\in[0,(2t_0))\wedge T]}|\mathbf{x}(t)|\left(\frac{1}{1-B(t_0)}+\frac{1}{(1-B(t_0))^2}\right).$$

It follows that

$$\sup_{t\in[0,T]}|\mathbf{y}(t)| \le \sup_{t\in[0,T]}|\mathbf{x}(t)| \sum_{i=1}^{\lfloor T/t_0\rfloor+1}\frac{1}{(1-B(t_0))^i}. \qquad \square$$

The next lemma is concerned with convergence of solutions to (B.1). Let $\mathbf{y}^n = (\mathbf{y}^n(t), t\in[0,T]) \in \mathbb{L}_\infty([0,T],\mathbb{R})$ solve equations

$$\mathbf{y}^n(t) = \mathbf{x}^n(t) + \int_0^t f^n(\mathbf{y}^n(t-s), t-s)\,dB(s), \tag{B.4}$$

where $\mathbf{x}^n = (\mathbf{x}^n(t), t\in[0,T]) \in \mathbb{L}_\infty([0,T],\mathbb{R})$ and the $f^n$ satisfy the hypotheses on $f$ in Lemma B.1 (in particular, are Lipshitz continuous in the first argument).

Lemma B.2.

1. *If*

$$\lim_{n\to\infty}\sup_{t\in[0,T]}\int_0^t |f^n(\mathbf{y}(t-s),t-s) - f(\mathbf{y}(t-s),t-s)|\,dB(s) = 0$$

*and the* $\mathbf{x}^n$ *converge to* $\mathbf{x}$ *in* $\mathbb{L}_\infty([0,T],\mathbb{R})$, *then the* $\mathbf{y}^n$ *converge to* $\mathbf{y}$ *in* $\mathbb{L}_\infty([0,T],\mathbb{R})$.
2. *Suppose that* $f^n(y,t)\to f(y,t)$ *as* $n\to\infty$ *for every* $y$ *and* $t$ *and that* $f^n(y,t)$ *is monotonically decreasing in* $n$ *and monotonically increasing in* $y$, *that is,* $f^m(y,t)\le f^n(y,t)$ *for* $m\ge n$ *and* $f^n(y_1,t)\le f^n(y_2,t)$ *for* $y_1\le y_2$. *If* $\mathbf{x}^n(t)\to\mathbf{x}(t)$ *for all* $t\in[0,T]$ *and* $\sup_n\sup_{t\in[0,T]}|\mathbf{x}^n(t)|<\infty$, *then* $\mathbf{y}^n(t)\to\mathbf{y}(t)$ *for all* $t\in[0,T]$.

Proof. In order to prove part 1, introduce

$$\tilde{\mathbf{x}}^n(t) = \mathbf{x}^n(t) - \mathbf{x}(t) + \int_0^t (f^n(\mathbf{y}(t-s),t-s) - f(\mathbf{y}(t-s),t-s))\,dB(s).$$

We have,

$$
\begin{aligned}
|\mathbf{y}^n(t) - \mathbf{y}(t)| &\le |\tilde{\mathbf{x}}^n(t)| + \int_0^t |f^n(\mathbf{y}^n(t-s), t-s) - f^n(\mathbf{y}(t-s), t-s)|\, dB(s) \\
&\le |\tilde{\mathbf{x}}^n(t)| + \int_0^t |\mathbf{y}^n(t-s) - \mathbf{y}(t-s)|\, dB(s),
\end{aligned}
$$

so by Lemma B.1 $\sup_{t\in[0,T]} |\mathbf{y}^n(t) - \mathbf{y}(t)| \le \rho(T) \sup_{t\in[0,T]} |\tilde{\mathbf{x}}^n(t)|$. The hypotheses imply that $\lim_{n\to\infty} \sup_{t\in[0,T]} |\tilde{\mathbf{x}}^n(t)| = 0$. Thus, $\lim_{n\to\infty} \sup_{t\in[0,T]} |\mathbf{y}^n(t) - \mathbf{y}(t)| = 0$.

We now prove part 2. Given $m \in \mathbb{N}$, we define

$$\overline{\mathbf{x}}^m(t) = \sup_{n\ge m} \mathbf{x}^n(t), \qquad \underline{\mathbf{x}}^m(t) = \inf_{n\ge m} \mathbf{x}^n(t). \tag{B.5}$$

Let $\overline{\mathbf{y}}^{n,m}$ and $\underline{\mathbf{y}}^{n,m}$ be defined by the respective equations

$$\overline{\mathbf{y}}^{n,m}(t) = \overline{\mathbf{x}}^m(t) + \int_0^t f^n(\overline{\mathbf{y}}^{n,m}(t-s), t-s)\, dB(s), \tag{B.6a}$$

$$\underline{\mathbf{y}}^{n,m}(t) = \underline{\mathbf{x}}^m(t) + \int_0^t f^n(\underline{\mathbf{y}}^{n,m}(t-s), t-s)\, dB(s). \tag{B.6b}$$

Since $f^n(y,t)$ is monotonically increasing in $y$ and is monotonically decreasing in $n$, an application of the method of successive approximations to (B.6a) and (B.6b) with the initial approximations being zero functions shows that the $\overline{\mathbf{y}}^{n,m}(t)$ and $\underline{\mathbf{y}}^{n,m}(t)$ monotonically decrease in $n$ for every $t$ and $m$. Besides, the sequences $\{\sup_{t\in[0,T]} |\overline{\mathbf{y}}^{n,m}(t)|, n \in \mathbb{N}\}$ and $\{\sup_{t\in[0,T]} |\overline{\mathbf{y}}^{n,m}(t)|, n \in \mathbb{N}\}$ are bounded, as it follows by Lemma B.1. Thus, for $t \in [0,T]$, there exist finite limits $\overline{\mathbf{y}}^m(t) = \lim_{n\to\infty} \overline{\mathbf{y}}^{n,m}(t)$ and $\underline{\mathbf{y}}^m(t) = \lim_{n\to\infty} \underline{\mathbf{y}}^{n,m}(t)$. On letting $n \to \infty$ in (B.6a) and (B.6b) and applying Lebesgue's bounded convergence theorem, we obtain that

$$\overline{\mathbf{y}}^m(t) = \overline{\mathbf{x}}^m(t) + \int_0^t f(\overline{\mathbf{y}}^m(t-s), t-s)\, dB(s), \tag{B.7a}$$

$$\underline{\mathbf{y}}^m(t) = \underline{\mathbf{x}}^m(t) + \int_0^t f(\underline{\mathbf{y}}^m(t-s), t-s)\, dB(s). \tag{B.7b}$$

As $m \to \infty$, the $\overline{\mathbf{x}}^m(t)$ monotonically decrease to $\mathbf{x}(t)$ and the $\underline{\mathbf{x}}^m(t)$ monotonically increase to $\mathbf{x}(t)$. An application of the method of successive approximations to (B.7a) and (B.7b) shows that the $\overline{\mathbf{y}}^m(t)$ monotonically decrease and the $\underline{\mathbf{y}}^m(t)$ monotonically increase. The limits satisfy (B.1), so $\overline{\mathbf{y}}^m(t) \to \mathbf{y}(t)$ and $\underline{\mathbf{y}}^m(t) \to \mathbf{y}(t)$, as $m \to \infty$. Also, since $\underline{\mathbf{x}}^m(t) \le \mathbf{x}^n(t) \le \overline{\mathbf{x}}^m(t)$ for $n \ge m$ by (B.5), we have by (B.1), (B.7a) and (B.7b) that $\underline{\mathbf{y}}^{n,m}(t) \le \mathbf{y}^n(t) \le \overline{\mathbf{y}}^{n,m}(t)$ for $n \ge m$. It follows that $\mathbf{y}^n(t) \to \mathbf{y}(t)$ as $n \to \infty$. $\square$

We now study regularity properties of solutions to (B.1) for equations arising in Theorem 2.2. Given $\mathbf{z} = (\mathbf{z}(t), t \in [0,T]) \in \mathbb{L}_\infty([0,T],\mathbb{R})$, consider the condition

$$\lim_{\varepsilon\to 0} \sup_{t\in[0,T]} \int_0^t \mathbf{1}_{\{0<|\mathbf{z}(t-s)|<\varepsilon\}}\, dB(s) = 0. \tag{B.8}$$

Note that (B.8) implies the condition of part 1 of Lemma B.2 for $f^n(y,t) = (y + n\mathbf{z}(t))^+ - n\mathbf{z}(t)^+$ and $f(y,t) = y\mathbf{1}_{\{\mathbf{z}(t)>0\}} + y^+\mathbf{1}_{\{\mathbf{z}(t)=0\}}$.

LEMMA B.3. *If $\mathbf{z}$ admits limits on the right and limits on the left and $B(t)$ is continuous on $[0,T]$ with $B(0) = 0$, then condition (B.8) holds.*

*Assume that $f(y,t) = y\mathbf{1}_{\{\mathbf{z}(t)>0\}} + y^+\mathbf{1}_{\{\mathbf{z}(t)=0\}}$ in (B.1). If both $\mathbf{z}$ and $\mathbf{x}$ are right-continuous with left-hand limits and condition (B.8) holds, then $\mathbf{y}$ is right-continuous with left-hand limits. If $\mathbf{z}$ is right-continuous with left-hand limits, $\mathbf{x}$ is continuous, and $B(t)$ is continuous on $[0,T]$ with $B(0) = 0$, then $\mathbf{y}$ is continuous.*

PROOF. If $\mathbf{z}$ admits limits on the right and limits on the left, then the function $\mathbf{1}_{\{0<|\mathbf{z}(s)|<\varepsilon\}}$ has at most countably many points of discontinuity, which implies that the function $\int_0^t \mathbf{1}_{\{0<|\mathbf{z}(t-s)|<\varepsilon\}}\, dB(s)$ is continuous in $t$ when $B(t)$ is a continuous function with $B(0) = 0$. By Dini's theorem, the monotonic convergence of $\int_0^t \mathbf{1}_{\{0<|\mathbf{z}(t-s)|<\varepsilon\}}\, dB(s)$ to zero as $\varepsilon \to 0$ is uniform in $t \in [0,T]$, as required.

Suppose that $\mathbf{x}$ and $\mathbf{z}$ are right-continuous with left-hand limits and condition (B.8) holds. Let $\hat{\mathbf{y}}^n = (\hat{\mathbf{y}}^n(t), t \in [0,T])$ solve the equation

$$\hat{\mathbf{y}}^n(t) = \mathbf{x}(t) + \int_0^t ((\hat{\mathbf{y}}^n(t-s) + n\mathbf{z}(t-s))^+ - n\mathbf{z}(t-s)^+)\, dB(s).$$

By Lemma B.1, the functions $\hat{\mathbf{y}}^n$ are right-continuous with left-hand limits. By part 1 of Lemma B.2, the $\hat{\mathbf{y}}^n$ converge to $\mathbf{y}$ for the topology of compact convergence. Therefore, the latter function is also right-continuous with left-hand limits. The argument for the case where $\mathbf{z}$ is right-continuous with left-hand limits, $\mathbf{x}$ is continuous, and $B(t)$ is continuous on $[0,T]$ with $B(0) = 0$ is similar. □

## APPENDIX C: MARTINGALE LEMMAS

This section contains results on certain processes being martingales which are instrumental in the proofs above. Let $A = (A(t), t \in \mathbb{R}_+)$ be an integer–valued nonnegative process with nondecreasing trajectories from $\mathbb{D}(\mathbb{R}_+, \mathbb{R})$. We denote $\tau_i = \inf\{t : A(t) \geq i\}, i \in \mathbb{N}$, so the jumps of $A$ occur at the $\tau_i$. Obviously, $A(t) = \sum_{i=1}^{\infty} \mathbf{1}_{\{\tau_i \leq t\}}$. Let $\xi_i, i \in \mathbb{N}$, be nonnegative random variables.

The next lemma is an extension of Lemma 5.2 in Krichagina and Puhalskii [15] (see also Lemma 8 in Reed [22]).

LEMMA C.1. *Let* $\mathbf{A} = (\mathcal{A}_t, t \in \mathbb{R}_+)$ *be a filtration such that the random variables* $\tau_{j\wedge(A(t)+1)}$ *and* $\xi_{j\wedge A(t)}$*, where* $j \in \mathbb{N}$*, are* $\mathcal{A}_t$*-measurable. Let* $\tilde{\mathcal{A}}_i$ *represent the complete* $\sigma$*-algebras generated by the events* $\Lambda \cap \{\tau_i > t\}$*, where* $t \in \mathbb{R}_+$ *and* $\Lambda \in \mathcal{A}_t$*. Suppose that* $(\beta_i(x,y), x \in \mathbb{R}_+, y \in \mathbb{R}_+), i \in \mathbb{N}$*, are real-valued Borel functions such that* $\mathbf{E}\beta_i(\tau_i, \xi_i)^2 < \infty$*. The following assertions hold.*

1. *The* $\tau_i$ *are* $\mathbf{A}$*-predictable stopping times.*
2. *Suppose, in addition, that there exists a nondecreasing sequence of* $\sigma$*-algebras* $\hat{\mathcal{A}}_i, i \in \mathbb{N}$*, such that* $\tilde{\mathcal{A}}_i \subset \hat{\mathcal{A}}_i$ *and such that the random variables* $\tau_j$*, where* $j = 1, 2, \ldots, i$*, and* $\xi_j$*, where* $j = 1, 2, \ldots, i-1$*, are* $\hat{\mathcal{A}}_i$*-measurable. Let, for* $k \in \mathbb{N}$*,*

$$R_k(t) = \sum_{i=1}^{A(t)\wedge k} \beta_i(\tau_i, \xi_i),$$

$$\langle R_k\rangle(t) = \sum_{i=1}^{A(t)\wedge k} \mathbf{E}(\beta_i(\tau_i, \xi_i)^2 | \hat{\mathcal{A}}_i).$$

   *If* $\mathbf{E}(\beta_i(\tau_i, \xi_i) | \hat{\mathcal{A}}_i) = 0$*, then the processes* $R_k = (R_k(t), t \in \mathbb{R}_+)$ *and* $(R_k(t)^2 - \langle R_k\rangle(t), t \in \mathbb{R}_+)$ *are* $\mathbf{A}$*-square-integrable martingales.*
3. *The inclusions* $\tilde{\mathcal{A}}_i \subset \hat{\mathcal{A}}_i$ *hold provided the* $\mathcal{A}_t$ *are the complete* $\sigma$*-algebras generated by the random variables* $\tau_{j\wedge(A(t)+1)}$ *and* $\xi_{j\wedge A(t)}$*, where* $j \in \mathbb{N}$*, and the* $\hat{\mathcal{A}}_i$ *are the complete* $\sigma$*-algebras generated by the random variables* $\tau_j$*, where* $j = 1, 2, \ldots, i$*, and* $\xi_j$*, where* $j = 1, 2, \ldots, i-1$*.*
4. *If the* $\xi_i, i \in \mathbb{N}$*, are independent, and if, for each* $i \in \mathbb{N}$*,* $\xi_i$ *is independent of the* $\tau_j$ *for* $j \le i$*, and* $\mathbf{E}\beta_i(x, \xi_i) = 0$ *for every* $x$*, then, for* $\mathbf{A}$ *and* $\hat{\mathcal{A}}_i$ *defined as in part* 3 *and for* $R_k$ *and* $\langle R_k\rangle$ *defined as in part* 2*, the processes* $R_k$ *are* $\mathbf{A}$*-square-integrable martingales with the predictable quadratic variation processes* $\langle R_k\rangle$*.*

PROOF. The $\tau_i$ are $\mathbf{A}$-stopping times because, for $t \in \mathbb{R}_+$ and $i \in \mathbb{N}$, $\{\tau_i \le t\} = \{\tau_{i\wedge(A(t)+1)} \le t\}$ which follows on noting that $\tau_{A(t)+1} > t$. The latter set belongs to $\mathcal{A}_t$.

The proof of the $\mathbf{A}$-predictability of the $\tau_i$ follows the argument of Krichagina and Puhalskii [15]. The time $\tau_1$ is $\mathbf{A}$-predictable since the times $(1 - 1/l)\tau_1, l = 1, 2, \ldots,$ are stopping times predicting $\tau_1$ (note that $\tau_1$ is $\mathcal{A}_0$-measurable). For $i > 1$, introduce $\sigma_{i,l} = \tau_i - (\tau_i - \tau_{\nu(i)})/l, l \in \mathbb{N}$, where $\tau_{\nu(i)} = \max\{\tau_p : \tau_p < \tau_i\}$ and $\tau_{\nu(i)} = 0$ if no such $p$ exists. Obviously, $\sigma_{i,l} \uparrow \tau_i$ as $l \to \infty$ and $\sigma_{i,l} < \tau_i$ on the set $\{\tau_i > 0\}$. We show that the $\sigma_{i,l}$ are $\mathbf{A}$-stopping time. For $t \in \mathbb{R}_+$, in view of the fact that $\{\tau_i \le t\} \subset \{\sigma_{i,l} \le t\} \subset \{\tau_{\nu(i)} \le t\}$,

$$\{\sigma_{i,l} \le t\} = \{\tau_i \le t\} \cup (\{\sigma_{i,l} \le t\} \cap \{\tau_{\nu(i)} \le t\} \cap \{\tau_{i\wedge(A(t)+1)} > t\}).$$

The first set on the right belongs to $\mathcal{A}_t$, as has been proved. Let $\hat{\tau}_{\nu(i)} = \max\{\tau_{p\wedge(A(t)+1)} : \tau_{p\wedge(A(t)+1)} < \tau_{i\wedge(A(t)+1)}\}$ and $\hat{\sigma}_{i,l} = \tau_{i\wedge(A(t)+1)} - (\tau_{i\wedge(A(t)+1)} - \hat{\tau}_{\nu(i)})/l$. The random variables $\hat{\tau}_{\nu(i)}$ and $\hat{\sigma}_{i,l}$ are $\mathcal{A}_t$-measurable. If $\tau_{\nu(i)} \le t$ and $\tau_{i\wedge(A(t)+1)} > t$, then $\tau_{\nu(i)} = \hat{\tau}_{\nu(i)}$. On the other hand, if $\hat{\tau}_{\nu(i)} \le t$ and $\tau_{i\wedge(A(t)+1)} > t$, then $\tau_{\nu(i)} = \hat{\tau}_{\nu(i)}$. We conclude that

$$\begin{aligned}\{\sigma_{i,l} \le t\} &\cap \{\tau_{\nu(i)} \le t\} \cap \{\tau_{i\wedge(A(t)+1)} > t\} \\ &= \{\hat{\sigma}_{i,l} \le t\} \cap \{\hat{\tau}_{\nu(i)} \le t\} \cap \{\tau_{i\wedge(A(t)+1)} > t\}.\end{aligned}$$

The set on the right is in $\mathcal{A}_t$. Thus, the $\sigma_{i,l}, l \in \mathbb{N}$, are $\mathbf{A}$-stopping times which predict $\tau_i$. Part 1 has been proved.

We prove part 2. As $\{A(t) \ge i\} = \{\tau_{i\wedge(A(t)+1)} \le t\}$, the random variables $A_t$ are $\mathcal{A}_t$-measurable, so the $R_k(t)$ are $\mathcal{A}_t$-measurable. Since $\sup_{t\in\mathbb{R}_+} \mathbf{E}R_k(t)^2 < \infty$, to prove that $R_k$ is an $\mathbf{A}$-square-integrable martingale, it is enough to prove that

$$\mathbf{E} \sum_{i=1}^{A(s\wedge\sigma)\wedge k} \beta_i(\tau_i, \xi_i) = 0 \tag{C.1}$$

for any $\mathbf{A}$-stopping time $\sigma$ (see, e.g., Jacod and Shiryaev [12], I.1.44). We have

$$\mathbf{E} \sum_{i=1}^{A(s\wedge\sigma)\wedge k} \beta_i(\tau_i, \xi_i) = \sum_{i=1}^{k} \mathbf{E}\mathbf{1}_{\{\tau_i \le s\wedge\sigma\}} \beta_i(\tau_i, \xi_i).$$

Since $\{\tau_i > s\wedge\sigma\} = \bigcup(\{\tau_i > u\} \cap \{u \ge s\wedge\sigma\})$, where the union is over all positive rational $u$, we have, by the fact that $s\wedge\sigma$ is an $\mathbf{A}$-stopping time and hence $\{u \ge s\wedge\sigma\} \in \mathcal{A}_u$, and the definition of $\tilde{\mathcal{A}}_i$, that $\{\tau_i > s\wedge\sigma\} \in \tilde{\mathcal{A}}_i$. By the inclusion $\tilde{\mathcal{A}}_i \subset \hat{\mathcal{A}}_i$, $\mathbf{E}\mathbf{1}_{\{\tau_i \le s\wedge\sigma\}} \beta_i(\tau_i, \xi_i) = \mathbf{E}(\mathbf{1}_{\{\tau_i \le s\wedge\sigma\}} \mathbf{E}(\beta_i(\tau_i, \xi_i)|\hat{\mathcal{A}}_i))$. The latter conditional expectation equals zero by hypotheses. Equality (C.1) has been proved.

In order to prove that $(R_k(t)^2 - \langle R_k\rangle(t), t \in \mathbb{R}_+)$ is an $\mathbf{A}$-square-integrable martingale, we show that $\mathbf{E}R_k(s\wedge\sigma)^2 = \mathbf{E}\langle R_k\rangle(s\wedge\sigma)$ for any $\mathbf{A}$-stopping time $\sigma$. The definition implies that

$$\begin{aligned}R_k(s\wedge\sigma)^2 &= \sum_{i=1}^{k} \mathbf{1}_{\{\tau_i \le s\wedge\sigma\}} \beta_i(\tau_i, \xi_i)^2 \\ &\quad + 2\sum_{i=1}^{k} \sum_{j=i+1}^{k} \mathbf{1}_{\{\tau_i \le s\wedge\sigma\}} \mathbf{1}_{\{\tau_j \le s\wedge\sigma\}} \beta_i(\tau_i, \xi_i)\beta_j(\tau_j, \xi_j).\end{aligned}$$

Next, for $i<j$, in analogy with an earlier argument, using the facts that $\tau_i$ and $\xi_i$ are $\hat{\mathcal{A}}_j$-measurable and that $\tau_j$ is $\hat{\mathcal{A}}_j$-measurable,

$$\begin{aligned}&\mathbf{E}\mathbf{1}_{\{\tau_i\le s\wedge\sigma\}}\mathbf{1}_{\{\tau_j\le s\wedge\sigma\}}\beta_i(\tau_i,\xi_i)\beta_j(\tau_j,\xi_j)\\&\qquad=\mathbf{E}(\mathbf{1}_{\{\tau_i\le s\wedge\sigma\}}\mathbf{1}_{\{\tau_j\le s\wedge\sigma\}}\beta_i(\tau_i,\xi_i)\mathbf{E}(\beta_j(\tau_j,\xi_j)|\hat{\mathcal{A}}_j))=0.\end{aligned}$$

Therefore,

$$\mathbf{E}R_k(s\wedge\sigma)^2=\mathbf{E}\sum_{i=1}^k\mathbf{1}_{\{\tau_i\le s\wedge\sigma\}}\mathbf{E}(\beta_i(\tau_i,\xi_i)^2|\hat{\mathcal{A}}_i)=\mathbf{E}\langle R_k\rangle(s\wedge\sigma).$$

Part 2 has been proved.

By Brémaud [4], Appendix A3, Theorem 25, the flow $\mathcal{A}_t, t\in\mathbb{R}_+$, is right-continuous, so we prove part 3 by showing that each generator of $\tilde{\mathcal{A}}_i$ is in $\hat{\mathcal{A}}_i$. Since $\{\tau_i>\ t\}=\{A(t)+1\le i\}$, we have that, if $j>i$, then $\{\tau_{j\wedge(A(t)+1)}\le u\}\cap\{\tau_i>\ t\}=\{\tau_{i\wedge(A(t)+1)}\le u\}\cap\{\tau_i>t\}$ which event is seen to be in $\hat{\mathcal{A}}_i$, and $\{\tau_{j\wedge(A(t)+1)}\le u\}\cap\{\tau_i>t\}\in\hat{\mathcal{A}}_i$ when $j\le i$. Similarly, if $j\ge i$, then $\{\xi_{j\wedge A(t)}\le u\}\cap\{\tau_i>t\}=\{\xi_{(i-1)\wedge A(t)}\le u\}\cap\{\tau_i>t\}\in\hat{\mathcal{A}}_i$, and if $j<i$, then $\{\xi_{j\wedge A(t)}\le u\}\cap\{\tau_i>t\}\in\hat{\mathcal{A}}_i$. Part 3 has been proved.

We prove part 4. Under the hypotheses, $\mathbf{E}(\beta_i(\tau_i,\xi_i)|\hat{\mathcal{A}}_i)=\mathbf{E}(\beta_i(x,\xi_i)|\hat{\mathcal{A}}_i)|_{x=\tau_i}=0$. By parts 2 and 3, it remains to prove that $\langle R_k\rangle$ is $\mathbf{A}$-predictable. As $\langle R_k\rangle(t)=\sum_{i=1}^k\mathbf{1}_{\{\tau_i\le t\}}\mathbf{E}\beta_i(x,\xi_i)^2|_{x=\tau_i}$, the $\tau_i$ are $\mathbf{A}$-predictable, and the $\mathbf{E}\beta_i(x,\xi_i)^2|_{x=\tau_i}$ are $\mathcal{A}_{\tau_i-}$-measurable, it follows that $\langle R_k\rangle$ is $\mathbf{A}$-predictable, Dellacherie [5], Subsection 5 of Section 1 of Chapter 4. □

For $k\in\mathbb{N}$, introduce the two-parameter processes $L_k=(L_k(t,x),t\in\mathbb{R}_+,x\in\mathbb{R}_+)$ with

$$L_k(t,x)=\sum_{i=1}^{A(t)\wedge k}\left(\mathbf{1}_{\{0<\xi_i\le x\}}-\int_0^{\xi_i\wedge x}\mathbf{1}_{\{u>0\}}\frac{dF_i(u)}{1-F_i(u-)}\right),\tag{C.2}$$

where $F_i$ denotes the distribution function of $\xi_i$. Define also

$$\langle L_k\rangle(t,x)=\sum_{i=1}^{A(t)\wedge k}\int_0^{\xi_i\wedge x}\mathbf{1}_{\{u>0\}}\frac{1-F_i(u)}{(1-F_i(u-))^2}\,dF_i(u).\tag{C.3}$$

Let, for $t\in\mathbb{R}_+$ and $x\in\mathbb{R}_+$, complete $\sigma$-algebras $\hat{\mathcal{F}}_{t,x}$ be generated by the random variables $\tau_{i\wedge(A(t)+1)}$, $\xi_{i\wedge A(t)}$, and $\mathbf{1}_{\{\tau_i\le s\}}$, where $i\in\mathbb{N}$, $s\in\mathbb{R}_+$, and $s\le t+x$, and by $\mathbf{1}_{\{\tau_i\le s\}}\mathbf{1}_{\{\xi_i\le y\}}$, where $i\in\mathbb{N}$, $s\in\mathbb{R}_+$, $y\in\mathbb{R}_+$, and $s+y\le t+x$. We denote $\mathcal{F}_{t,x}=\bigcap_{\varepsilon>0,\delta>0}\hat{\mathcal{F}}_{t+\varepsilon,x+\delta}$ (see Figure 1) and note that both $L_k(t,x)$ and $\langle L_k\rangle(t,x)$ are $\mathcal{F}_{t,x}$-measurable. (As a matter of fact, these random variables are $\mathcal{F}_{t,0}$-measurable.) In the next lemma, we define conditional probabilities as being equal to zero when the conditioning events are of probability zero.

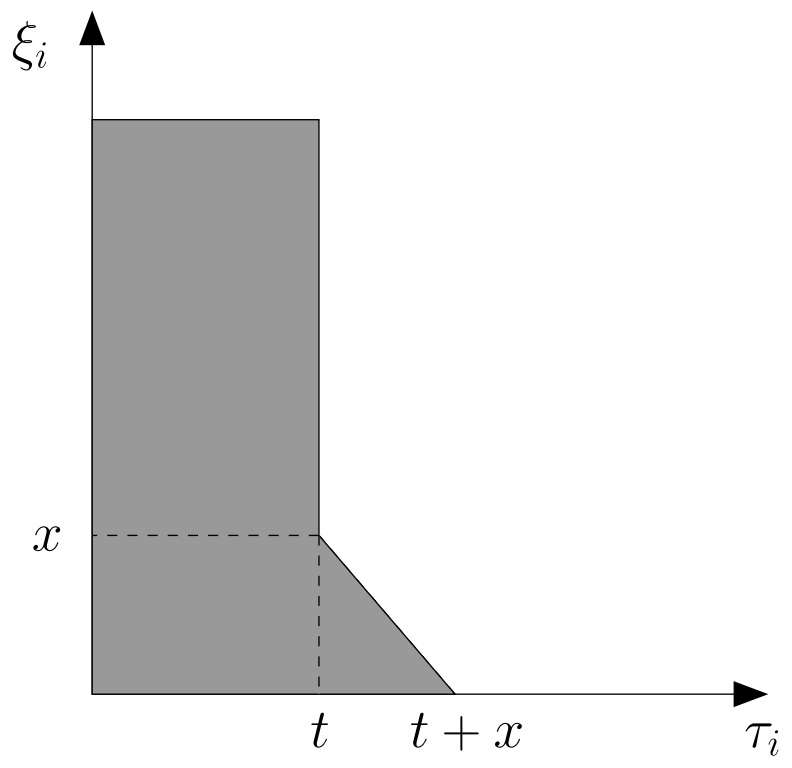

FIG. 1. *The $\sigma$-algebra $\mathcal{F}_{t,x}$.*

LEMMA C.2. *Suppose that, given $i \in \mathbb{N}$, $s \in \mathbb{R}_+$, and $x \in \mathbb{R}_+$, the random variable $\xi_i$ is independent of the random variables $\xi_j$ for $j \neq i$, of the $\tau_j$ for $j \leq i$, and of the $\tau_j \wedge (s+x)$ for $j > i$, conditioned on the event $\{\tau_i \geq s, \xi_i > x\}$. Suppose also that each $\xi_i$, for $i \in \mathbb{N}$, is independent of the $\tau_j$ for $j \leq i$. Then for $s \leq t$ and $x \leq y$,*

$$\mathbf{E}(\Box L_k((s,x),(t,y))|\mathcal{F}_{s,x}) = 0$$

*and*

$$\mathbf{E}((\Box L_k((s,x),(t,y)))^2|\mathcal{F}_{s,x}) = \Box\langle L_k\rangle((s,x),(t,y)).$$

PROOF. Assuming $x$ and $y$ are held fixed, we apply Lemma C.1 with $\mathcal{A}_t = \mathcal{F}_{t,x}$ and $\hat{\mathcal{A}}_i$ being the complete $\sigma$-algebra generated by $\xi_j$ for $j < i$, by $\tau_j$ for $j \leq i$, by $\mathbf{1}_{\{\tau_j \leq s\}}\mathbf{1}_{\{\tau_i > t\}}$ for $j > i$, $t \in \mathbb{R}_+$, and $s \leq t+x$, and by $\mathbf{1}_{\{\tau_j \leq s\}}\mathbf{1}_{\{\tau_i > t\}}\mathbf{1}_{\{0<\xi_j \leq u\}}$ for $j \geq i$, $s \in \mathbb{R}_+$, $t \in \mathbb{R}_+$, $u \in \mathbb{R}_+$, and $s+u \leq t+x$.

Let us check that $\tilde{\mathcal{A}}_i \subset \hat{\mathcal{A}}_i$. Pick $j \in \mathbb{N}$. The inclusions $\{\tau_{j\wedge(A(t)+1)} \leq u\} \cap \{\tau_i > t\} \in \hat{\mathcal{A}}_i$ and $\{\xi_{j\wedge A(t)} \leq u\} \cap \{\tau_i > t\} \in \hat{\mathcal{A}}_i$ follow by part 3 of Lemma C.1. For $t \in \mathbb{R}_+$ and $s \leq t+x$, $\{\tau_j \leq s\} \cap \{\tau_i > t\} \in \hat{\mathcal{A}}_i$ by definition. Similarly, $\{\tau_j \leq s\} \cap \{0 < \xi_j \leq u\} \cap \{\tau_i > t\} \in \hat{\mathcal{A}}_i$, where $i \in \mathbb{N}$, $s \in \mathbb{R}_+$, $u \in \mathbb{R}_+$, and $s+u \leq t+x$. Thus, if $\Lambda \in \mathcal{A}_t$, then $\Lambda \cap \{\tau_i > t+\varepsilon\} \in \hat{\mathcal{A}}_i$ for all $\varepsilon > 0$, hence, $\Lambda \cap \{\tau_i > t\} \in \hat{\mathcal{A}}_i$.

Introduce

$$\beta_i(v) = \mathbf{1}_{\{x<v\leq y\}} - \int_{v\wedge x}^{v\wedge y} \mathbf{1}_{\{u>0\}} \frac{dF_i(u)}{1-F_i(u-)}, \tag{C.4}$$

so that

$$L_k(t,y) - L_k(t,x) = \sum_{i=1}^{A(t)\wedge k} \beta_i(\xi_i). \tag{C.5}$$

We show that

$$\mathbf{E}(\beta_i(\xi_i)|\hat{\mathcal{A}}_i) = 0. \tag{C.6}$$

We have for natural numbers $i < k_1 < k_2 < \cdots < k_p$, $i < n_1 < n_2 < \cdots < n_m$, nonnegative real numbers $\tilde{u}$, $\tilde{t}$, $\tilde{s}$, $t_1, t_2, \ldots, t_p$, $s_1, s_2, \ldots, s_p$, $t'_1, t'_2, \ldots, t'_m$, $s'_1, s'_2, \ldots, s'_m$, and $u_1, \ldots, u_m$ such that $s_j \le t_j + x$, $\tilde{s} + \tilde{u} \le \tilde{t} + x$, and $s'_j + u_j \le t'_j + x$, and for Borel bounded functions $f_1, \ldots, f_l$, $g_1, \ldots, g_i$, $h_1^{(1)}, \ldots, h_p^{(1)}$, $h_1^{(2)}, \ldots, h_m^{(2)}$, and $f$, on noting that by (C.4) $\beta_i(v) = \mathbf{1}_{\{v>x\}}\beta_i(v)$ and that the condition $\tilde{s} + \tilde{u} \le \tilde{t} + x$ implies the identity $\mathbf{1}_{\{\xi_i > x\}}\mathbf{1}_{\{\tilde{t} < \tau_i \le \tilde{s}\}}\mathbf{1}_{\{0<\xi_i \le \tilde{u}\}} = 0$,

$$\begin{aligned}
&\mathbf{E}\Bigg(\prod_{j=1}^{i-1} f_j(\xi_j) \prod_{j=1}^{i} g_j(\tau_j) \prod_{j=1}^{p} h_j^{(1)}(\mathbf{1}_{\{\tau_{k_j} \le s_j\}}\mathbf{1}_{\{\tau_i > t_j\}}) \\
&\quad \times \prod_{j=1}^{m} h_j^{(2)}(\mathbf{1}_{\{\tau_{n_j} \le s'_j\}}\mathbf{1}_{\{\tau_i > t'_j\}}\mathbf{1}_{\{0<\xi_{n_j} \le u_j\}}) \\
&\qquad \times f(\mathbf{1}_{\{\tilde{t} < \tau_i \le \tilde{s}\}}\mathbf{1}_{\{0<\xi_i \le \tilde{u}\}})\beta_i(\xi_i)\Bigg) \\
&= \mathbf{E}\Bigg(\prod_{j=1}^{i-1} f_j(\xi_j) \prod_{j=1}^{i} g_j(\tau_j) \prod_{j=1}^{p} (h_j^{(1)}(\mathbf{1}_{\{\tau_{k_j} \le s_j\}})\mathbf{1}_{\{\tau_i > t_j\}} + h_j^{(1)}(0)\mathbf{1}_{\{\tau_i \le t_j\}}) \\
&\quad \times \prod_{j=1}^{m} (h_j^{(2)}(\mathbf{1}_{\{\tau_{n_j} \le s'_j\}}\mathbf{1}_{\{0<\xi_{n_j} \le u_j\}})\mathbf{1}_{\{\tau_i > t'_j\}} + h_j^{(2)}(0)\mathbf{1}_{\{\tau_i \le t'_j\}}) \\
&\qquad \times f(0)\mathbf{1}_{\{\xi_i > x\}}\beta_i(\xi_i)\Bigg)
\end{aligned} \tag{C.7}$$

$$\begin{aligned}
&= f(0) \sum_{\substack{J \subset \{1,2,\ldots,p\},\\ J' \subset \{1,2,\ldots,m\}}} \mathbf{E}\Bigg(\prod_{j=1}^{i-1} f_j(\xi_j) \prod_{j=1}^{i} g_j(\tau_j) \prod_{j \in J} h_j^{(1)}(\mathbf{1}_{\{\tau_{k_j} \le s_j\}})\mathbf{1}_{\{\tau_i > t_j\}} \\
&\quad \times \prod_{j \in \{1,2,\ldots,p\} \setminus J} h_j^{(1)}(0)\mathbf{1}_{\{\tau_i \le t_j\}} \\
&\quad \times \prod_{j \in J'} h_j^{(2)}(\mathbf{1}_{\{\tau_{n_j} \le s'_j\}}\mathbf{1}_{\{0<\xi_{n_j} \le u_j\}})\mathbf{1}_{\{\tau_i > t'_j\}}
\end{aligned}$$

$$\times \prod_{j\in\{1,2,\dots,m\}\setminus J'} h_j^{(2)}(0)\mathbf{1}_{\{\tau_i\le t'_j\}}\mathbf{1}_{\{\xi_i>x\}}\beta_i(\xi_i)\Bigg).$$

We show that a generic term in the sum on the rightmost side of (C.7) equals zero. For given $J$ and $J'$, let $\hat t=\max_{j\in J}t_j\vee\max_{j\in J'}t'_j$ and $\breve t=\min_{j\in\{1,2,\dots,p\}\setminus J}t_j\wedge\min_{j\in\{1,2,\dots,m\}\setminus J'}t'_j$. Since $s_j\le\hat t+x$ and $s'_j\le\hat t+x$, we have

$$\begin{aligned}
&\mathbf{E}\Bigg(\prod_{j=1}^{i-1}f_j(\xi_j)\prod_{j=1}^{i}g_j(\tau_j)\prod_{j\in J}h_j^{(1)}(\mathbf{1}_{\{\tau_{k_j}\le s_j\}})\mathbf{1}_{\{\tau_i>t_j\}}\prod_{j\in\{1,2,\dots,p\}\setminus J}h_j^{(1)}(0)\mathbf{1}_{\{\tau_i\le t_j\}}\\
&\quad\times\prod_{j\in J'}h_j^{(2)}(\mathbf{1}_{\{\tau_{n_j}\le s'_j\}}\mathbf{1}_{\{0<\xi_{n_j}\le u_j\}})\mathbf{1}_{\{\tau_i>t'_j\}}\\
&\qquad\qquad\times\prod_{j\in\{1,2,\dots,m\}\setminus J'}h_j^{(2)}(0)\mathbf{1}_{\{\tau_i\le t'_j\}}\mathbf{1}_{\{\xi_i>x\}}\beta_i(\xi_i)\Bigg)\\
&\quad=\prod_{j\in\{1,2,\dots,p\}\setminus J}h_j^{(1)}(0)\prod_{j\in\{1,2,\dots,m\}\setminus J'}h_j^{(2)}(0)\\
&\qquad\times\mathbf{E}\Bigg(\prod_{j=1}^{i-1}f_j(\xi_j)\prod_{j=1}^{i}g_j(\tau_j)\prod_{j\in J}h_j^{(1)}(\mathbf{1}_{\{\tau_{k_j}\wedge(\hat t+x)\le s_j\}})\\
&\qquad\quad\times\prod_{j\in J'}h_j^{(2)}(\mathbf{1}_{\{\tau_{n_j}\wedge(\hat t+x)\le s'_j\}}\mathbf{1}_{\{0<\xi_{n_j}\le u_j\}})\\
&\qquad\qquad\times\mathbf{1}_{\{\hat t<\tau_i\le\breve t\}}\beta_i(\xi_i)\Big|\tau_i\ge\hat t,\xi_i>x\Bigg)\mathbf{P}(\tau_i\ge\hat t,\xi_i>x).
\end{aligned}$$

The independence hypotheses imply that

$$\begin{aligned}
&\mathbf{E}\Bigg(\prod_{j=1}^{i-1}f_j(\xi_j)\prod_{j=1}^{i}g_j(\tau_j)\prod_{j\in J}h_j^{(1)}(\mathbf{1}_{\{\tau_{k_j}\wedge(\hat t+x)\le s_j\}})\\
&\quad\times\prod_{j\in J'}h_j^{(2)}(\mathbf{1}_{\{\tau_{n_j}\wedge(\hat t+x)\le s'_j\}}\mathbf{1}_{\{0<\xi_{n_j}\le u_j\}})\\
&\qquad\times\mathbf{1}_{\{\hat t<\tau_i\le\breve t\}}\beta_i(\xi_i)\Big|\tau_i\ge\hat t,\xi_i>x\Bigg)\\
&\quad=\mathbf{E}\Bigg(\prod_{j=1}^{i-1}f_j(\xi_j)\prod_{j=1}^{i}g_j(\tau_j)\prod_{j\in J}h_j^{(1)}(\mathbf{1}_{\{\tau_{k_j}\wedge(\hat t+x)\le s_j\}})
\end{aligned}$$

$$\times \prod_{j \in J'} h_j^{(2)}(\mathbf{1}_{\{\tau_{n_j} \wedge (\hat{t}+x) \le s'_j\}} \mathbf{1}_{\{0<\xi_{n_j} \le u_j\}}) \mathbf{1}_{\{\hat{t}<\tau_i \le \check{t}\}} \Big| \tau_i \ge \hat{t}, \xi_i > x \Bigg)$$

$$\times \mathbf{E}(\beta_i(\xi_i)|\xi_i > x).$$

By (C.4), $\mathbf{E}(\beta_i(\xi_i)|\xi_i > x) = 0$. We conclude that the leftmost side of (C.7) equals zero, which establishes (C.6).

We now show that

$$\mathbf{E}(\beta_i(\xi_i)^2|\hat{\mathcal{A}}_i) = \int_{\xi_i \wedge x}^{\xi_i \wedge y} \mathbf{1}_{\{u>0\}} \frac{1 - F_i(u)}{(1 - F_i(u-))^2} \, dF_i(u). \tag{C.8}$$

Denoting

$$\overline{\beta}_i(\xi_i) = \beta_i(\xi_i)^2 - \int_{\xi_i \wedge x}^{\xi_i \wedge y} \mathbf{1}_{\{u>0\}} \frac{1 - F_i(u)}{(1 - F_i(u-))^2} \, dF_i(u)$$

and replicating the argument that established (C.6), we can see that (C.8) follows provided $\mathbf{E}(\overline{\beta}_i(\xi_i)|\xi_i > x) = 0$, which is a consequence of (C.4).

The required properties now follow by part 2 of Lemma C.1. $\square$

Let $\hat{\mathcal{M}}_t$ denote the complete $\sigma$-algebra generated by the events $\{\tau_i \le s\} \cap \{\xi_i \le x\}$ and $\{\tau_i \le s\}$ where $i \in \mathbb{N}$, $s \in \mathbb{R}_+$, $x \in \mathbb{R}_+$, and $s + x \le t$, and let $\mathcal{M}_t = \bigcap_{\varepsilon>0} \hat{\mathcal{M}}_{t+\varepsilon}$ (see Figure 2). The flow $\mathbf{M} = (\mathcal{M}_t, t \in \mathbb{R}_+)$ is a filtration. Note that $L_k(t,x)$ and $\langle L_k \rangle(t,x)$ are $\mathcal{M}_{t+x}$-measurable and $\mathcal{F}_{t,x} \supset \mathcal{M}_{t+x}$. Define, for $k \in \mathbb{N}$ and $t \in \mathbb{R}_+$,

$$M_k(t) = \int_{\mathbb{R}_+^2} \mathbf{1}_{\{s+x \le t\}} \, dL_k(s,x) \tag{C.9}$$

and

$$\langle M_k \rangle(t) = \int_{\mathbb{R}_+^2} \mathbf{1}_{\{s+x \le t\}} \, d\langle L_k \rangle(s,x). \tag{C.10}$$

LEMMA C.3. *Under the independence hypotheses of Lemma C.2, the processes $M_k = (M_k(t), t \in \mathbb{R}_+)$ are $\mathbf{M}$-square integrable martingales with predictable quadratic variation processes $\langle M_k \rangle = (\langle M_k \rangle(t), t \in \mathbb{R}_+)$.*

PROOF. Since $M_k$ and $\langle M_k \rangle$ are right-continuous, it suffices to prove that $M_k$ is a square integrable martingale relative to the flow $\hat{\mathbf{M}} = (\hat{\mathcal{M}}_t, t \in \mathbb{R}_+)$ with predictable quadratic variation process $\langle M_k \rangle$.

Fix $s < t$. We first show that

$$\mathbf{E}((L_k(0,t) - L_k(0,s))|\hat{\mathcal{M}}_s) = 0, \tag{C.11}$$

$$\mathbf{E}((L_k(0,t) - L_k(0,s))^2|\hat{\mathcal{M}}_s) = \mathbf{E}(\langle L_k \rangle(0,t) - \langle L_k \rangle(0,s)|\hat{\mathcal{M}}_s). \tag{C.12}$$

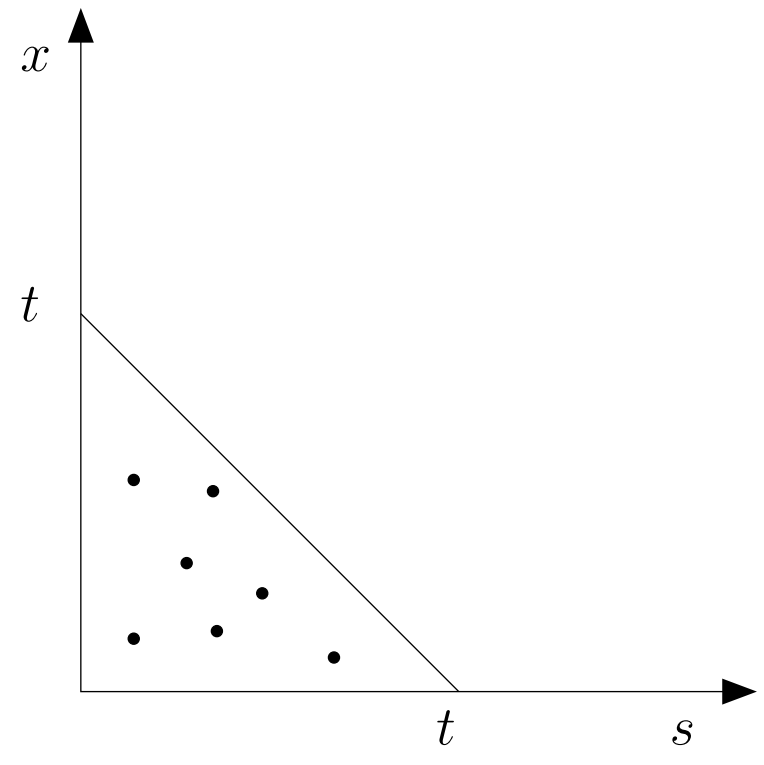

FIG. 2. *The $\sigma$-algebra $\mathcal{M}_t$.*

The argument is similar to the one used in the proof of Lemma C.2. Denote

$$\zeta_i = \mathbf{1}_{\{s<\xi_i\le t\}} - \int_{\xi_i\wedge s}^{\xi_i\wedge t} \mathbf{1}_{\{u>0\}} \frac{dF_i(u)}{1-F_i(u-)}, \tag{C.13}$$

$$\overline{\zeta}_i = \zeta_i^2 - \int_{\xi_i\wedge s}^{\xi_i\wedge t} \mathbf{1}_{\{u>0\}} \frac{1-F_i(u)}{(1-F_i(u-))^2}\, dF_i(u). \tag{C.14}$$

Let $\hat{\mathcal{M}}_s^j$ denote the $\sigma$-algebra generated by the $\xi_k$ for $k\le j$ and by $\hat{\mathcal{M}}_s$ when $j\in\mathbb{N}$ and let $\hat{\mathcal{M}}_s^0=\hat{\mathcal{M}}_s$. Obviously, $\zeta_j$ is $\hat{\mathcal{M}}_s^j$-measurable. We prove that

$$\mathbf{E}(\zeta_i|\hat{\mathcal{M}}_s^{i-1}) = 0, \tag{C.15}$$

$$\mathbf{E}(\overline{\zeta}_i|\hat{\mathcal{M}}_s^{i-1}) = 0. \tag{C.16}$$

We have for distinct natural numbers $k_1,k_2,\dots,k_p$ and $n_1,n_2,\dots,n_m$, none of which equals $i$, nonnegative real numbers $\tilde{t}$, $\tilde{s}$, $\tilde{x}$, $s_1,s_2,\dots,s_p$, $x_1,x_2,\dots,x_p$, and $t_1,t_2,\dots,t_m$, such that $\tilde{s}+\tilde{x}\le s$, $\tilde{t}\le s$, $s_j+x_j\le s$, and $t_j\le s$, and for Borel bounded functions $f_1$, $f_2$, $g_1,\dots,g_k$, $h_1^{(1)},\dots,h_p^{(1)}$, and $h_1^{(2)},\dots,h_m^{(2)}$,

$$\begin{aligned}
&\mathbf{E}\Bigg( f_1(\mathbf{1}_{\{\tau_i\le\tilde{s}\}}\mathbf{1}_{\{\xi_i\le\tilde{x}\}})f_2(\mathbf{1}_{\{\tau_i\le\tilde{t}\}})\\
&\qquad\times\prod_{j=1}^{i-1} g_j(\xi_j)\prod_{j=1}^{p} h_j^{(1)}(\mathbf{1}_{\{\tau_{k_j}\le s_j\}}\mathbf{1}_{\{\xi_{k_j}\le x_j\}})\prod_{j=1}^{m} h_j^{(2)}(\mathbf{1}_{\{\tau_{n_j}\le t_j\}})\zeta_i\Bigg)\\
&\quad= f_1(0)\mathbf{E}\Bigg( f_2(\mathbf{1}_{\{\tau_i\le\tilde{t}\}})\prod_{j=1}^{i-1} g_j(\xi_j)\prod_{j=1}^{p} h_j^{(1)}(\mathbf{1}_{\{\tau_{k_j}\wedge s\le s_j\}}\mathbf{1}_{\{\xi_{k_j}\le x_j\}})\\
&\qquad\qquad\times\prod_{j=1}^{m} h_j^{(2)}(\mathbf{1}_{\{\tau_{n_j}\wedge s\le t_j\}})\mathbf{1}_{\{\xi_i>s\}}\zeta_i\Bigg)
\end{aligned}$$

$$= f_1(0)\mathbf{E}\Bigg( f_2(\mathbf{1}_{\{\tau_i\le\tilde t\}})\prod_{j=1}^{i-1} g_j(\xi_j)\prod_{j=1}^{p} h_j^{(1)}(\mathbf{1}_{\{\tau_{k_j}\wedge s\le s_j\}}\mathbf{1}_{\{\xi_{k_j}\le x_j\}}) \times\prod_{j=1}^{m} h_j^{(2)}(\mathbf{1}_{\{\tau_{n_j}\wedge s\le t_j\}})\Big|\xi_i>s\Bigg) \times\mathbf{E}(\zeta_i|\xi_i>s)\mathbf{P}(\xi_i>s).$$

Since $\mathbf{E}(\zeta_i|\xi_i>s)=0$ by (C.13), (C.15) follows. Equality (C.16) is obtained analogously.

By (C.2) and (C.13), $L_k(0,t)-L_k(0,s)=\sum_{i=1}^k \mathbf{1}_{\{\tau_i=0\}}\zeta_i$. Since $\{\tau_i=0\}\in\hat{\mathcal{M}}_s$, by (C.15)

$$\mathbf{E}((L_k(0,t)-L_k(0,s))|\hat{\mathcal{M}}_s)=\sum_{i=1}^k \mathbf{1}_{\{\tau_i=0\}}\mathbf{E}(\zeta_i|\hat{\mathcal{M}}_s)=0.$$

By (C.15), if $j<i$, then $\mathbf{E}(\zeta_j\zeta_i|\hat{\mathcal{M}}_s)=\mathbf{E}(\zeta_j\mathbf{E}(\zeta_i|\hat{\mathcal{M}}_s^{i-1}))=0$. We thus obtain on recalling (C.3), (C.14) and (C.16),

$$\begin{aligned}&\mathbf{E}((L_k(0,t)-L_k(0,s))^2|\hat{\mathcal{M}}_s)\\&\quad=\sum_{i=1}^k \mathbf{1}_{\{\tau_i=0\}}\mathbf{E}(\zeta_i^2|\hat{\mathcal{M}}_s)+2\sum_{j<i}\mathbf{1}_{\{\tau_j=0\}}\mathbf{1}_{\{\tau_i=0\}}\mathbf{E}(\zeta_j\zeta_i|\hat{\mathcal{M}}_s)\\&\quad=\mathbf{E}(\langle L_k\rangle(0,t)-\langle L_k\rangle(0,s)|\hat{\mathcal{M}}_s).\end{aligned}$$

This completes the proof of (C.11) and (C.12).

By (C.2) and (C.9),

$$M_k(t)=\sum_{i=1}^{A(t)\wedge k}\left(\mathbf{1}_{\{0<\xi_i\le t-\tau_i\}}-\int_0^{\xi_i\wedge(t-\tau_i)}\mathbf{1}_{\{u>0\}}\frac{dF_i(u)}{1-F_i(u-)}\right), \tag{C.17}$$

which implies that $M_k(t)$ is $\hat{\mathcal{M}}_t$-measurable. We prove that $M_k$ is a martingale.

Let $0=s_0^l<s_1^l<s_2^l<\cdots<s_l^l=t$ be such that $\max_i(s_i^l-s_{i-1}^l)\to 0$ as $l\to\infty$ and $s=s_m^l$ for some $m$. Define (see Figure 3)

$$J_{l,s}(u,x)=\sum_{i=1}^m \mathbf{1}_{\{u\in(s_{i-1}^l,s_i^l]\}}\mathbf{1}_{\{x\le s-s_{i-1}^l\}}+\mathbf{1}_{\{u=0\}}\mathbf{1}_{\{x\le s\}}, \tag{C.18a}$$

$$J_{l,t}(u,x)=\sum_{i=1}^l \mathbf{1}_{\{u\in(s_{i-1}^l,s_i^l]\}}\mathbf{1}_{\{x\le t-s_{i-1}^l\}}+\mathbf{1}_{\{u=0\}}\mathbf{1}_{\{x\le t\}}, \tag{C.18b}$$

$$M_{k,l}(s)=\int_{\mathbb{R}_+^2} J_{l,s}(u,x)\,dL_k(u,x), \tag{C.18c}$$

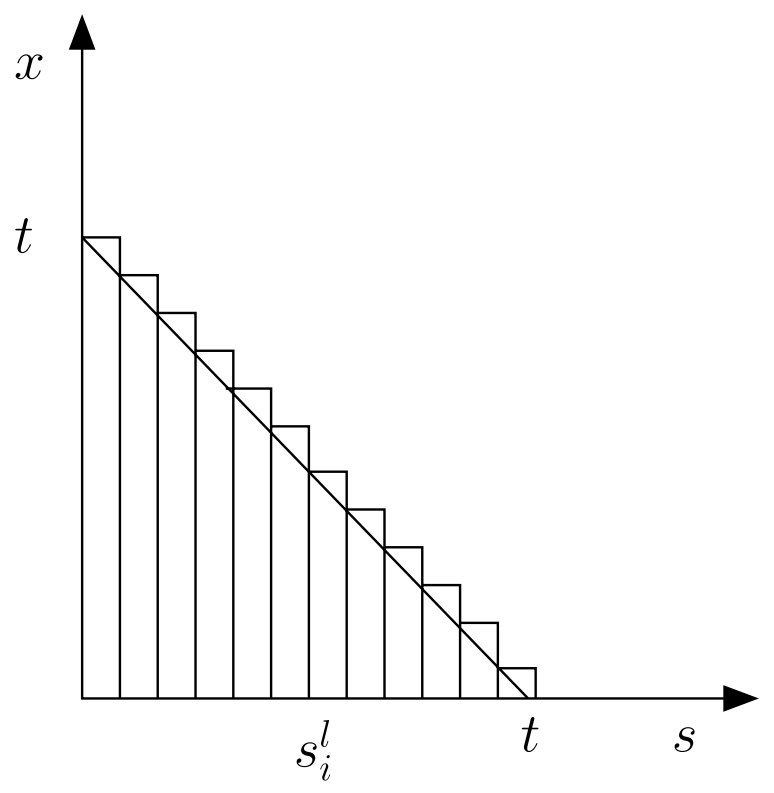

FIG. 3. *Approximating martingales.*

$$M_{k,l}(t) = \int_{\mathbb{R}_+^2} J_{l,t}(u,x)\, dL_k(u,x). \tag{C.18d}$$

Note that

$$M_{k,l}(s) = \sum_{i=1}^{m} \square L_k((s_{i-1}^l, 0), (s_i^l, s - s_{i-1}^l)) + L_k(0,s), \tag{C.19a}$$

$$M_{k,l}(t) = \sum_{i=1}^{l} \square L_k((s_{i-1}^l, 0), (s_i^l, t - s_{i-1}^l)) + L_k(0,t). \tag{C.19b}$$

Since $J_{l,s}(u,x) \to \mathbf{1}_{\{u+x\le s\}}$ and $J_{l,t}(u,x) \to \mathbf{1}_{\{u+x\le t\}}$ as $l\to\infty$, it follows by (C.9), (C.18c) and (C.18d) that $M_{k,l}(s) \to M_k(s)$ and $M_{k,l}(t) \to M_k(t)$ as $l \to \infty$. By Lemma C.2 and (C.19a),

$$\mathbf{E} M_{k,l}(s)^2 = \sum_{i=1}^{m} \mathbf{E} \square \langle L_k \rangle ((s_{i-1}^l, 0), (s_i^l, s - s_{i-1}^l)) \le \mathbf{E}\langle L_k\rangle(s,s) < \infty.$$

Similarly, $\mathbf{E} M_{k,l}(t) \le \mathbf{E}\langle L_k\rangle(t,t) < \infty$. Hence,

$$\lim_{l\to\infty} \mathbf{E}(M_{k,l}(s)|\hat{\mathcal{M}}_s) = M_k(s), \tag{C.20a}$$

$$\lim_{l\to\infty} \mathbf{E}(M_{k,l}(t)|\hat{\mathcal{M}}_s) = \mathbf{E}(M_k(t)|\hat{\mathcal{M}}_s). \tag{C.20b}$$

Lemma C.2 implies that $\mathbf{E}(\square L_k((s_{i-1}^l, 0), (s_i^l, t - s_i^l))|\mathcal{F}_{s_{i-1}^l, s - s_{i-1}^l}) = \mathbf{E}(\square L_k((s_{i-1}^l, 0), (s_i^l, s - s_i^l))|\mathcal{F}_{s_{i-1}^l, s - s_{i-1}^l})$ for $i = 1, 2, \ldots, m$ and that $\mathbf{E}(\square L_k((s_{i-1}^l, 0), (s_i^l, t - s_i^l))|\mathcal{F}_{s_{i-1}^l, 0}) = 0$ for $i = m+1, \ldots, l$. Since $\hat{\mathcal{M}}_s \subset (\bigcap_{i=1,2,\ldots,m} \mathcal{F}_{s_{i-1}^l, s - s_{i-1}^l}) \cap (\bigcap_{i=m+1,\ldots,l} \mathcal{F}_{s_{i-1}^l, 0})$, by (C.11), (C.19a) and (C.19b),

$\mathbf{E}(M_{k,l}(t)|\hat{\mathcal{M}}_s) = \mathbf{E}(M_{k,l}(s)|\hat{\mathcal{M}}_s)$. The martingale property of $M_k$ follows by (C.20a) and (C.20b).

We now compute the predictable quadratic variation process. By (C.17), recalling (C.19b),

$$M_k(t) - M_{k,l}(t)$$
$$(C.21)\qquad = -\sum_{i=1}^{A(t)\wedge k}\sum_{j=1}^{l}\mathbf{1}_{\{\tau_i\in(s_{j-1}^l,s_j^l]\}}\biggl(\mathbf{1}_{\{\xi_i\in(t-\tau_i,t-s_{j-1}^l]\}} - \int_{\xi_i\wedge(t-\tau_i)}^{\xi_i\wedge(t-s_{j-1}^l)}\mathbf{1}_{\{u>0\}}\frac{dF_i(u)}{1-F_i(u-)}\biggr).$$

We apply part 4 of Lemma C.1 with $\beta_i(x,y) = -\sum_{j=1}^{l}\mathbf{1}_{\{x\in(s_{j-1}^l,s_j^l]\}}(\mathbf{1}_{\{y\in(t-x,t-s_{j-1}^l]\}} - \int_{y\wedge(t-x)}^{y\wedge(t-s_{j-1}^l)}\mathbf{1}_{\{u>0\}}\,dF_i(u)/(1-F_i(u-)))$. According to the lemma,

$$\begin{aligned}
&\mathbf{E}(M_k(t) - M_{k,l}(t))^2\\
&\quad = \mathbf{E}\sum_{i=1}^{A(t)\wedge k}\sum_{j=1}^{l}\mathbf{1}_{\{\tau_i\in(s_{j-1}^l,s_j^l]\}}\int_{t-\tau_i}^{t-s_{j-1}^l}\mathbf{1}_{\{u>0\}}\frac{1-F_i(u)}{1-F_i(u-)}\,dF_i(u)\\
&\quad \le \sum_{i=1}^{k}\mathbf{E}\sum_{j=1}^{l}\mathbf{1}_{\{\tau_i\in(s_{j-1}^l,s_j^l]\}}(F_i(t-s_{j-1}^l) - F_i(t-\tau_i))\\
&\quad = \sum_{i=1}^{k}\mathbf{E}(F_i(t-s_{j(\tau_i)-1}^l) - F_i(t-\tau_i)),
\end{aligned}$$

where $j(\tau_i)$ is defined by the requirement that $\tau_i \in (s_{j(\tau_i)-1}^l, s_{j(\tau_i)}^l]$. Since $t - s_{j(\tau_i)-1}^l$ converges from the right to $t-\tau_i$ as $l\to\infty$, $F(t-s_{j(\tau_i)-1}^l)$ converges to $F(t-\tau_i)$. By bounded convergence,

$$(C.22)\qquad \lim_{l\to\infty}\mathbf{E}(M_k(t) - M_{k,l}(t))^2 = 0.$$

Similarly,

$$(C.23)\qquad \lim_{l\to\infty}\mathbf{E}(M_k(s) - M_{k,l}(s))^2 = 0.$$

By (C.19a) and (C.19b),

$$(C.24)\qquad \begin{aligned}
M_{k,l}(t) - M_{k,l}(s) &= \sum_{i=1}^{l}\square L_k((s_{i-1}^l, s - s_{i-1}^l\wedge s),(s_i^l, t - s_{i-1}^l))\\
&\quad + (L_k(0,t) - L_k(0,s)).
\end{aligned}$$

By Lemma C.2, for $i < j$,

$$\begin{aligned}&\mathbf{E}(\square L_k((s^l_{i-1}, s - s^l_{i-1} \wedge s), (s^l_i, t - s^l_{i-1}))\\ &\quad\times \square L_k((s^l_{j-1}, s - s^l_{j-1} \wedge s), (s^l_j, t - s^l_{j-1}))|\mathcal{F}_{s^l_{j-1}, s - s^l_{j-1}\wedge s})\\ &\quad= \square L_k((s^l_{i-1}, s - s^l_{i-1} \wedge s), (s^l_i, t - s^l_{i-1}))\\ &\qquad\times \mathbf{E}(\square L_k((s^l_{j-1}, s - s^l_{j-1} \wedge s), (s^l_j, t - s^l_{j-1}))|\mathcal{F}_{s^l_{j-1}, s - s^l_{j-1}\wedge s}) = 0\end{aligned}$$

and

$$\begin{aligned}&\mathbf{E}((\square L_k((s^l_{i-1}, s - s^l_{i-1} \wedge s), (s^l_i, t - s^l_{i-1})))^2|\mathcal{F}_{s^l_{j-1}, s - s^l_{j-1}\wedge s})\\ &\quad= \square\langle L_k\rangle((s^l_{i-1}, s - s^l_{i-1} \wedge s), (s^l_i, t - s^l_{i-1})),\end{aligned}$$

so, due to the fact that $\hat{\mathcal{M}}_s \subset \mathcal{F}_{s^l_{j-1}, s - s^l_{j-1}\wedge s}$,

$$\begin{aligned}&\mathbf{E}(\square L_k((s^l_{i-1}, s - s^l_{i-1} \wedge s), (s^l_i, t - s^l_{i-1}))\\ &\quad\times \square L_k((s^l_{j-1}, s - s^l_{j-1} \wedge s), (s^l_j, t - s^l_{j-1}))|\hat{\mathcal{M}}_s) = 0\end{aligned} \tag{C.25}$$

and

$$\begin{aligned}&\mathbf{E}((\square L_k((s^l_{i-1}, s - s^l_{i-1} \wedge s), (s^l_i, t - s^l_{i-1})))^2|\hat{\mathcal{M}}_s)\\ &\quad= \mathbf{E}(\square\langle L_k\rangle((s^l_{i-1}, s - s^l_{i-1} \wedge s), (s^l_i, t - s^l_{i-1}))|\hat{\mathcal{M}}_s).\end{aligned} \tag{C.26}$$

Since $L_k(0,t)$ and $L_k(0,s)$ are $\mathcal{F}_{0,0}$-measurable,

$$\begin{aligned}&\mathbf{E}(\square L_k((s^l_{i-1}, s - s^l_{i-1} \wedge s), (s^l_i, t - s^l_{i-1}))(L_k(0,t) - L_k(0,s))|\mathcal{F}_{s^l_{i-1}, s - s^l_{i-1}\wedge s})\\ &\quad= (L_k(0,t) - L_k(0,s))\\ &\qquad\times \mathbf{E}(\square L_k((s^l_{i-1}, s - s^l_{i-1} \wedge s), (s^l_i, t - s^l_{i-1}))|\mathcal{F}_{s^l_{i-1}, s - s^l_{i-1}\wedge s}) = 0,\end{aligned}$$

so

$$\begin{aligned}&\mathbf{E}(\square L_k((s^l_{i-1}, s - s^l_{i-1} \wedge s), (s^l_i, t - s^l_{i-1}))\\ &\quad\times (L_k(0,t) - L_k(0,s))|\hat{\mathcal{M}}_s) = 0.\end{aligned} \tag{C.27}$$

Putting together (C.12), (C.18a), (C.18b), (C.24)–(C.27) yields

$$\begin{aligned}&\mathbf{E}((M_{k,l}(t) - M_{k,l}(s))^2|\hat{\mathcal{M}}_s)\\ &\quad= \mathbf{E}\left(\sum_{i=1}^{l} \square\langle L_k\rangle((s^l_{i-1}, s - s^l_{i-1} \wedge s), (s^l_i, t - s^l_{i-1}))\right.\\ &\qquad\left. + (\langle L_k\rangle(0,t) - \langle L_k\rangle(0,s))\middle|\hat{\mathcal{M}}_s\right)\end{aligned} \tag{C.28}$$

$$= \mathbf{E}\biggl(\int_{\mathbb{R}_+^2} (J_{l,t}(u,x) - J_{l,s}(u,x))\, d\langle L_k\rangle(u,x) \bigg| \hat{\mathcal{M}}_s\biggr).$$

Since $J_{l,s}(u,x) \to \mathbf{1}_{\{u+x\le s\}}$ and $J_{l,t}(u,x) \to \mathbf{1}_{\{u+x\le t\}}$ as $l \to \infty$,

$$\begin{aligned}
&\lim_{l\to\infty} \int_{\mathbb{R}_+^2} J_{l,s}(u,x)\, d\langle L_k\rangle(u,x) = \int_{\mathbb{R}_+^2} \mathbf{1}_{\{u+x\le s\}}\, d\langle L_k\rangle(u,x),\\
&\lim_{l\to\infty} \int_{\mathbb{R}_+^2} J_{l,t}(u,x)\, d\langle L_k\rangle(u,x) = \int_{\mathbb{R}_+^2} \mathbf{1}_{\{u+x\le t\}}\, d\langle L_k\rangle(u,x),
\end{aligned} \tag{C.29}$$

so by bounded convergence, a.s.,

$$\begin{aligned}
&\lim_{l\to\infty} \mathbf{E}\biggl(\int_{\mathbb{R}_+^2} (J_{l,t}(u,x) - J_{l,s}(u,x))\, d\langle L_k\rangle(u,x) \bigg| \hat{\mathcal{M}}_s\biggr)\\
&\qquad = \mathbf{E}\biggl(\int_{\mathbb{R}_+^2} \mathbf{1}_{\{s<u+x\le t\}}\, d\langle L_k\rangle(u,x) \bigg| \hat{\mathcal{M}}_s\biggr).
\end{aligned} \tag{C.30}$$

Putting together (C.22), (C.23), (C.28) and (C.30) and recalling the definition of $\langle M_k\rangle$ in (C.10) obtains the equality

$$\mathbf{E}((M_k(t) - M_k(s))^2 | \hat{\mathcal{M}}_s) = \mathbf{E}((\langle M_k\rangle(t) - \langle M_k\rangle(s)) | \hat{\mathcal{M}}_s) \qquad \text{a.s.}$$

It remains to verify that the process $\langle M_k\rangle$ is $\hat{\mathbf{M}}$-predictable. By (C.3),

$$\langle L_k\rangle(t,x) = \sum_{i=1}^{k} \hat{L}_i(t,x),$$

where

$$\hat{L}_i(t,x) = \mathbf{1}_{\{\tau_i\le t\}} \int_0^{\xi_i\wedge x} \mathbf{1}_{\{u>0\}} \frac{1 - F_i(u)}{(1 - F_i(u-))^2}\, dF_i(u). \tag{C.31}$$

Denote

$$\hat{M}_i(t) = \int_{\mathbb{R}_+^2} \mathbf{1}_{\{u+x\le t\}}\, d\hat{L}_i(u,x). \tag{C.32}$$

Since $\langle M_k\rangle = \sum_{i=1}^k \hat{M}_i$, it suffices to prove that each process $\hat{M}_i = (\hat{M}_i(t), t \in \mathbb{R}_+)$ is $\hat{\mathbf{M}}$-predictable. It is $\hat{\mathbf{M}}$-adapted which follows by the representation

$$\hat{M}_i(t) = \mathbf{1}_{\{\tau_i\le t\}} \int_0^{\xi_i\wedge(t-\tau_i)} \mathbf{1}_{\{u>0\}} \frac{1 - F_i(u)}{(1 - F_i(u-))^2}\, dF_i(u).$$

By (C.31) and (C.32),

$$\hat{M}_i(t) = \mathbf{1}_{\{\tau_i\le t\}} \int_0^{\xi_i\wedge(t-\tau_i)} \mathbf{1}_{\{x>0\}} \frac{1 - F_i(x)}{(1 - F_i(x-))^2}\, dF_i(x). \tag{C.33}$$

Consider the decomposition $\hat{M}_i = \hat{M}_i^c + \hat{M}_i^d$, where $\hat{M}_i^c$ is a continuous-path adapted process and $\hat{M}_i^d$ is a pure-jump adapted process. By being continuous and adapted, $\hat{M}_i^c$ is $\hat{\mathbf{M}}$-predictable. The process $\hat{M}_i^d$ is of the form

$$\hat{M}_i^d(t) = \sum_{u>0} \mathbf{1}_{\{\tau_i+u\le t\}} \mathbf{1}_{\{u\le \xi_i\}} \frac{1-F_i(u)}{(1-F_i(u-))^2} \Delta F_i(u), \tag{C.34}$$

where the sum is taken over all positive times of jumps of $F_i$. Note that $\tau_i$ is an $\hat{\mathbf{M}}$-stopping time by the definition of $\hat{\mathcal{M}}_t$. Since $u>0$ in (C.34), $\tau_i+u$ is an $\hat{\mathbf{M}}$-predictable stopping time. We show that $\mathbf{1}_{\{u\le\xi_i\}}$ is $\hat{\mathcal{M}}_{(\tau_i+u)-}$-measurable. Note that $\{\xi_i<u\} = \bigcup_{t\in\mathbb{Q}_+}(\{\tau_i+\xi_i<t\}\cap\{\tau_i+u>t\})$. By the representation $\{\tau_i+\xi_i<t\} = \bigcup_{s,x\in\mathbb{Q}_+:s+x<t}(\{\tau_i\le s\}\cap\{\xi_i\le x\})$, we have that $\{\tau_i+\xi_i<t\}\in\hat{\mathcal{M}}_t$, hence, $\{\xi_i<u\}\in\hat{\mathcal{M}}_{(\tau_i+u)-}$. By Dellacherie [5], IV.1.5, each summand on the right of (C.34) is an $\hat{\mathbf{M}}$-predictable process. It follows that $\hat{M}_i^d$ is too. □

## APPENDIX D: HOFFMANN–JØRGENSEN'S CONVERGENCE IN DISTRIBUTION

In this section, we state the properties of Hoffmann–Jørgensen's convergence in distribution envoked in the proofs of Section 4. We recall the definition stated in the introductory part of the paper. Let $(\Omega,\mathcal{F},\mathbf{P})$ be a probability space. Given a real-valued function $\xi$ on $\Omega$, the outer expectation $\mathbf{E}^*\xi$ of $\xi$ is defined as the infimum of $\mathbf{E}\zeta$ over all random variables $\zeta$ on $(\Omega,\mathcal{F},\mathbf{P})$ such that $\zeta\ge\xi$ a.s. and $\mathbf{E}\zeta$ exists. Let $S$ be a metric space made into a measurable space by endowing it with the Borel $\sigma$-algebra $\mathcal{B}(S)$. Given a sequence $X_n$ of maps from $\Omega$ to $S$ and a measurable map $X$ from $(\Omega,\mathcal{F})$ to $(S,\mathcal{B}(S))$, we say that the $X_n$ converge to $X$ in distribution if

$$\lim_{n\to\infty} \mathbf{E}^* f(X_n) = \mathbf{E} f(X) \tag{D.1}$$

for all bounded continuous functions $f$ on $S$. We also define $\mathbf{E}_* f(X_n) = -\mathbf{E}^*(-f(X_n))$.

The next result is the continuous mapping principle, which is Theorem 1.3.6 in van der Vaart and Wellner [23].

THEOREM D.1. *Let $S'$ be a metric space and function $f: S\to S'$ be Borel and continuous at all points of a set $S_0\subset S$. If the $X_n$ converge in distribution to $X$ in $S$ and $X$ assumes values in $S_0$, then the $f(X_n)$ converge in distribution to $f(X)$ in $S'$.*

Let $\tau$ denote the topology on $S$.

COROLLARY D.1. *Suppose that $\tau'$ is a metric topology on $S$ which is coarser than $\tau$ and that open balls in $(S,\tau)$ are Borel sets in $(S,\tau')$. If the $X_n$ converge in distribution to $X$ for topology $\tau'$, a subset $S_0$ of $S$ is such that convergence to the elements of $S_0$ in $\tau'$ implies convergence in $\tau$, and $X$ assumes values in $S_0$, then the $X_n$ converge in distribution to $X$ for topology $\tau$.*

The next theorem relaxes the requirement of the separability of the range of $X$ in the extended continuous mapping principle of van der Vaart and Wellner [23], cf., Theorem 1.11.1 in [23]. The proof is modelled on that of Lemma 3.1.13 in Puhalskii [20].

THEOREM D.2. *Let $f_n$ and $f$ be maps from $S$ to a metric space $S'$. Suppose that $f$ is measurable and that the $f_n(x_n)$ converge to $f(x)$ whenever the $x_n$ converge to $x$ and $x \in S_0$, where $S_0$ is a Borel subset of $S$. If the $X_n$ converge in distribution to $X$ and $X$ assumes values in $S_0$, then the $f_n(X_n)$ converge in distribution to $f(X)$.*

PROOF. Note that $f(X)$ is a random element of $S'$, so convergence in distribution of the $f(X_n)$ to $f(X)$ is well defined. Let $g$ be a continuous bounded function from $S'$ to $\mathbb{R}$. Let $B_\delta(x)$ denote the closed ball of radius $\delta$ about $x \in S$ and introduce

$$h_k(x) = \inf_{\delta>0} \sup_{y\in B_\delta(x)} \sup_{n\geq k} g(f_n(y)).$$

The $h_k$ are bounded upper semicontinuous functions on $S$ and $h_k(x) \geq g(f_n(x))$ for $n \geq k$, so, for arbitrary $k$, envoking the Portmanteau theorem (Theorem 1.3.4 in van der Vaart and Wellner [23])

$$\limsup_{n\to\infty} \mathbf{E}^* g(f(X_n)) \leq \limsup_{n\to\infty} \mathbf{E}^* h_k(X_n) \leq \mathbf{E} h_k(X).$$

The convergence hypotheses on the $f_n$ and the continuity of $g$ imply that the $h_k(x)$ monotonically decrease to $g(f(x))$ as $k \to \infty$ for all $x \in S_0$. By bounded convergence, $\lim_k \mathbf{E} h_k(X) = \mathbf{E} g(f(X))$, and we conclude that

$$\limsup_{n\to\infty} \mathbf{E}^* g(f(X_n)) \leq \mathbf{E} g(f(X)).$$

Applying this inequality to the function $-g$, we have that

$$\liminf_{n\to\infty} \mathbf{E}_* g(f(X_n)) \geq \mathbf{E} g(f(X)).$$

Since $\mathbf{E}^* g(f(X_n)) \geq \mathbf{E}_* g(f(X_n))$, the proof is over. $\square$

We now review the extensions of Prohorov's theorem. We say that the sequence $X_n$ is asymptotically measurable if

$$\lim_{n\to\infty} (\mathbf{E}^* f(X_n) - \mathbf{E}_* f(X_n)) = 0 \tag{D.2}$$

for all bounded continuous functions on $S$.

Define $\mathbf{P}^*(A) = \mathbf{E}^* \mathbf{1}_A$ for $A \subset \Omega$. We say that the sequence $X_n$ is asymptotically tight if for every $\varepsilon > 0$ there exists a compact set $K$ such that

$$\limsup_{n\to\infty} \mathbf{P}^*(X_n \in S \setminus K^\delta) \le \varepsilon$$

for all $\delta > 0$, where $K^\delta = \{x \in S : m(x, K) < \delta\}$, $m$ being the metric on $S$.

The next four results are copied from van der Vaart and Wellner [23].

THEOREM D.3. *If the sequence $X_n$ converges in distribution to a random element with a tight probability law, then it is asymptotically tight.*

THEOREM D.4. *If the sequence $X_n$ is asymptotically tight and $f : S \to S'$ is a continuous mapping, then the sequence $f(X_n)$ is asymptotically tight.*

THEOREM D.5. *If the sequence $X_n$ is asymptotically tight and asymptotically measurable, then it has a subsequence which converges in distribution to a random element with a tight probability law.*

THEOREM D.6. *Suppose that the sequence $X_n$ is asymptotically tight and (D.2) holds for all functions $f$ from a subalgebra of the algebra of bounded continuous functions which separates points in $S$. Then the sequence $X_n$ is asymptotically measurable.*

The next group of results is concerned with joint convergence. Let $Y_n$ be a sequence of maps from $\Omega$ to a metric space $S'$. The following is Lemma 1.4.3 from van der Vaart and Wellner [23].

THEOREM D.7. *Each of the sequences $X_n$ and $Y_n$ is asymptotically tight if and only if the sequence $(X_n, Y_n)$ is asymptotically tight in $S \times S'$.*

We use the following corollary of Example 1.4.6 in van der Vaart and Wellner [23].

THEOREM D.8. *Suppose that the $X_n$ converge in distribution in $S$ to a separable random element $X$ and the $Y_n$ converge in distribution in $S'$ to a separable random element $Y$. If $X_n$ and $Y_n$ are independent for each $n$, then the $(X_n, Y_n)$ converge in distribution in $S \times S'$ to $(X, Y)$, where $X$ and $Y$ are independent.*

The next result is usually referred to as Slutsky's lemma; see Example 1.4.7 in van der Vaart and Wellner [23].

THEOREM D.9. *If the sequence $X_n$ converges in distribution in $S$ to a separable random element $X$ and the sequence $Y_n$ converges in distribution in $S'$ to a constant element $c \in S'$, then the sequence $(X_n, Y_n)$ converges in distribution in $S \times S'$ to $(X, c)$.*

The following result is in a similar vein.

THEOREM D.10. *Suppose that $X_n^\varepsilon$, where $\varepsilon > 0$, converge in distribution in $S$ to random elements $X^\varepsilon$ as $n \to \infty$ and that the $X^\varepsilon$ converge in distribution in $S$ to a random element $X$ as $\varepsilon \to 0$. If*

$$\lim_{\varepsilon \to 0} \limsup_{n \to \infty} \mathbf{P}^*(m(X_n, X_n^\varepsilon) > \delta) = 0$$

*for arbitrary $\delta > 0$, then the $X_n$ converge in distribution in $S$ to $X$.*

PROOF. As it follows by the Portmanteau theorem (Theorem 1.3.4 in van der Vaart and Wellner [23]), it suffices to prove that (D.1) holds for any real-valued bounded uniformly continuous function $f$ on $S$. We have

$$|\mathbf{E}^* f(X_n) - \mathbf{E} f(X)| \leq |\mathbf{E}^* f(X_n) - \mathbf{E}^* f(X_n^\varepsilon)| + |\mathbf{E}^* f(X_n^\varepsilon) - \mathbf{E} f(X^\varepsilon)| \\ + |\mathbf{E} f(X^\varepsilon) - \mathbf{E} f(X)|.$$

The hypotheses imply that

$$\limsup_{\varepsilon \to 0} \limsup_{n \to \infty} |\mathbf{E}^* f(X_n) - \mathbf{E} f(X)| \leq \limsup_{\varepsilon \to 0} \limsup_{n \to \infty} |\mathbf{E}^* f(X_n) - \mathbf{E}^* f(X_n^\varepsilon)|.$$

By the triangle inequality for the outer expectation, $|\mathbf{E}^* f(X_n) - \mathbf{E}^* f(X_n^\varepsilon)| \leq \mathbf{E}^* |f(X_n) - f(X_n^\varepsilon)|$. Given arbitrary $\gamma > 0$, let $\delta > 0$ be such that $|f(x) - f(y)| \leq \gamma$ if $m(x, y) < \delta$. Then $|f(X_n) - f(X_n^\varepsilon)| \leq \gamma + 2 \sup_{x \in S} |f(x)| \mathbf{1}_{\{m(X_n, X_n^\varepsilon) > \delta\}}$. We conclude that

$$\limsup_{\varepsilon \to 0} \limsup_{n \to \infty} |\mathbf{E}^* f(X_n) - \mathbf{E}^* f(X_n^\varepsilon)| \\ \leq \gamma + 2 \sup_{x \in S} |f(x)| \limsup_{\varepsilon \to 0} \limsup_{n \to \infty} \mathbf{P}^*(m(X_n, X_n^\varepsilon) > \delta) = \gamma.$$

Convergence (D.1) follows. □

Convergence in distribution to separable random elements can be metrized. Let $X$ and $Y$ be mappings from $\Omega$ into $S$. Suppose that $Y$ is measurable. We define

$$d_{\mathrm{BL}}^*(X, Y) = \sup_{f \in \mathrm{BL}_1} |\mathbf{E}^* f(X) - \mathbf{E} f(Y)|,$$

where $\mathrm{BL}_1$ denotes the set of real-valued functions on $S$ that are bounded in absolute value by 1 and admit a Lipshitz constant of 1. The next theorem appears in Dudley [7]; see also van der Vaart and Wellner [23], Section 12.

THEOREM D.11. *The sequence $X_n$ converges in distribution in $S$ to a separable random element $X$ if and only if* $\lim_{n\to\infty} d^*_{\mathrm{BL}_1}(X_n, X) = 0$.

**Acknowledgments.** The authors would like to thank the referees for their careful reading of the manuscript and helpful suggestions. In particular, the authors followed one of the referee's advice to incorporate general initial conditions in Theorem 2.2.

Department of Mathematical
and Statistical Sciences
University of Colorado Denver
P.O. Box 173364, Campus Box 170
Denver, Colorado 80217-3364
USA
and
IITP, Moscow
Bolshoy Karetny per. 19
Moscow 127994
Russia
E-mail: anatolii.puhalskii@ucdenver.edu

Leonard N. Stern School of Business
New York University
Kaufman Management Center
44 West 4th Street, KMC 8-79
New York, New York 10012
USA
E-mail: jreed@stern.nyu.edu